\documentclass[10pt,a4paper]{article}

\usepackage[margin=1.25in, a4paper]{geometry}
\usepackage{amsmath, bm}
\usepackage[utf8]{inputenc}
\usepackage{amsmath,mathtools,amsthm, enumitem, hyperref}
\usepackage{amssymb}
\usepackage{amsthm}
\usepackage{svg}
\usepackage{float}
\usepackage{pinlabel}
\usepackage{caption}
\usepackage{array}
\usepackage{esvect}
\usepackage{xfrac, bigints}
\usepackage[normalem]{ulem}
\usepackage{stmaryrd}
\usepackage[nottoc]{tocbibind}
\usepackage{enumitem}
\usepackage{graphicx,import}
\usepackage{mathrsfs}
\usepackage[english]{babel}

\usepackage{svg}
\usepackage{float}
\theoremstyle{plain}
\newtheorem{theorem}{Theorem}[section]
\newtheorem{lem}[theorem]{Lemma}
\newtheorem{prop}[theorem]{Proposition}
\newtheorem{coro}[theorem]{Corollary}

\theoremstyle{definition}
\newtheorem{defn}[theorem]{Definition}
\newtheorem{ex}[theorem]{Example}
\newtheorem{rem}[theorem]{Remark}
\newtheorem{ass}[theorem]{Assumption}

\numberwithin{figure}{section}
\numberwithin{equation}{section}

\newcommand{\la}{\langle}
\newcommand{\ra}{\rangle}
\newcommand{\bR}{\mathbb{R}}

\newcommand{\bZ}{\mathbb{Z}}

\newcommand{\pa}{\partial}
\newcommand{\rt}{\rightarrow}

\newcommand{\mancac}{\mathbf{Man}^{\mathbf{c,ac}}}
\newcommand{\manc}{\mathbf{Man}^{\mathbf{c}}}
\newcommand{\ssubset}{\subset\joinrel\subset}
\newcommand{\tu}{\tilde{u}}
\newcommand{\ep}{\varepsilon}

\newcommand{\cL}{\mathcal{L}}
\newcommand{\rel}{\mathrm{rel}}

\newcommand{\chL}{\check{L}}

\newcommand{\cF}{\mathcal{F}}

\newcommand{\cK}{\mathcal{K}}

\newcommand{\loc}{\mathrm{loc}}

\newcommand{\bM}{\mathbb{M}}
\newcommand{\rh}{\rho_1}
\newcommand{\rg}{\rho_2}

\newcommand{\bla}{\boldsymbol{\lambda}}
\newcommand{\lah}{\lambda_1}
\newcommand{\lag}{\lambda_2}
\newcommand{\ze}{\mathbf{0}}

\newcommand{\ind}{\operatorname{index}}

\newcommand{\m}{\mathrm{m}}

\newcommand{\brh}{\boldsymbol{\rho}}

\newcommand{\hu}{\hat{u}}

\newcommand{\cB}{\mathcal{B}}

\newcommand{\cV}{\mathcal{V}}

\newcommand{\we}{{\boldsymbol{\rho}^{\boldsymbol{\lambda}}}}

\newcommand{\cT}{\mathcal{T}}
\newcommand{\ord}{\mathrm{o}}
\newcommand{\Sym}{\mathrm{Sym}}
\newcommand{\ma}{\mathrm{a}}
\newcommand{\Ker}{\mathrm{Ker}}
\newcommand{\Coker}{\mathrm{Coker}}

\newcommand{\cP}{\mathcal{P}}
\newcommand{\rj}{\rho_j}

\newcommand{\cQ}{\mathcal{Q}}
\newcommand{\tsig}{\tilde{\sigma}}

\newcommand{\be}{\boldsymbol{\eta}}
\newcommand{\uV}{\underline{V}}

\newcommand{\de}{\rho\pa_{\rho}}

\newcommand{\ttil}{\tilde{t}}
\newcommand{\vol}{\textrm{vol}}
\newcommand{\sC}{\mathscr{C}}
\newcommand{\ttp}{{t'\partial_{t'}}}
\newcommand{\Ima}{\textrm{Im}}
\newcommand{\ttt}{{\tilde{t}\partial_{\tilde{t}}}}

\usepackage[T1]{fontenc}
\usepackage{imakeidx}
\makeindex
\date{}
\author{Spandan Ghosh 
	\footnote{ Mathematical Institute, University of Oxford, Oxford OX2 6GG. \\
	E-mail: ghoshs@maths.ox.ac.uk}
}

\title{Fredholm theorems for parabolic operators over a b-fibration}

\begin{document}
\maketitle
\begin{abstract}
	We introduce a framework for singularity analysis in geometric flows, using the language of manifolds with corners (and a-corners). We define the notion of a mixed fibration over a time interval, and study heat type operators of the form $\pa_v-\cL$ for $\pa_v$ being a time-like vector field with prescribed decay at a-boundaries and $\cL$ being a second-order relatively elliptic operator. We consider mappings between weighted Sobolev spaces, and prove Fredholm theorems and change-of-index formulae for these operators. We develop this theory with a view towards studying singularity formation as well as flowing out of singularities in flows such as mean curvature flow and Ricci flow.
\end{abstract}
\tableofcontents
\section{Introduction}
The analysis of singularities in geometric flows is important in order to understand long-term behaviour of the flow and its geometric/topological consequences like Perelman's proof of the Poincar\'{e} conjecture 
\cite{Per02, Per03a, Per03b}.
 One of the salient tools to describe singularity formation is the notion of Type I/Type II blow-ups. These notions provide important classifications of singularities that may form, for example we know that any Type I singularity of mean curvature flow must be modelled on a shrinker due to Huisken \cite{Hui90} and a subsequence along Type II rescalings must converge to an ancient solution due to Hamilton. However, simply knowing the Type I/II classification does not give a complete geometric picture of singularity formation. This is because taking a rescaling and then a subsequential limit only looks at the local picture of the singularity and information about the global geometry is lost.

In this paper, we work with a relatively new way to model and understand singularities of a flow, via the language of \textit{manifolds with corners}. This viewpoint was first articulated by the school of Melrose, for example \cite{Mel90}, and more recently by Joyce \cite[\S 5]{Joy16}, and was first applied to a geometric flow in the paper by Lira--Mazzeo--Pluda--S\'{a}ez \cite{LMPS23}, where they studied  the network flow problem in the plane and proved the existence of a regular network flow out of a smooth irregular network. There is also the more interesting problem of modelling singularity formation at a final time. This paper systematically puts both problems on the same framework as studying heat-type operators over a certain notion of fibration. 
 
 \paragraph*{\textbf{Fredholm theorems}}
 Although the equation of the flow is a \textit{nonlinear} parabolic partial differential equation, one needs to first study its linearisation, which is a linear second-order parabolic p.d.e. In this paper, we work with certain weighted Sobolev spaces and prove Fredholm theorems for the operator between these spaces. These results are analogous to that of Loya and Loya--Melrose \cite{LM03, Loy04, Loy05}.  However, we do not study Schwartz kernels as in loc. cit. Instead, we use direct Fourier analytic methods to prove our theorems, inspired from work of Agmon--Nirenberg \cite{AN63}. We also use methods inspired from work of Mazzeo \cite{Maz91a,Maz91b} to upgrade our Fredholm theory from spaces with weak regularity to spaces with stronger regularity which are better suited for nonlinear applications.

\paragraph{Outline of the paper}
We begin with background material in \S \ref{sec_background}, recalling the notions of manifolds with corners and a-corners and the parabolic interior estimates we need. We also introduce the notion of a mixed fibration which is central to our setup. In \S \ref{sec_applications}, we outline geometric problems where our theorems may be applied. In \S \ref{sec_heatoperators}, we explain our setup and prove a-priori estimates for our global operator (Theorem \ref{thm_degenerateapriori}), obtained by adding up interior estimates. Combined with a functional analytic lemma, this reduces semi-Fredholmness of the operator to proving $L^2$ estimates near each a-boundary (Theorem \ref{thm_fredholmlocal}). Then we introduce model problems at a-boundaries which control Fredholmness of the mapping, and relate these to notions of normal operators considered by Melrose and Loya--Melrose. In \S \ref{sec_models}, we reduce these model problems to checking the injectivity of certain normal operators and degenerate normal operators at the a-boundaries and a-corners, giving necessary and sufficient conditions for semi-Fredholmness. In \S \ref{sec_fredholm}, we prove the main Fredholm theorem (Theorem \ref{thm_mainthmfred}) by combining semi-Fredholmness of the global operator and its conjugate. Finally, we upgrade the regularity of our solutions in \S \ref{sec_optimalsobolev},write down a wall-crossing formula at non-vanishing a-boundaries in \S \ref{subsec_coi} (Theorem \ref{prop_polyhomogeneous}), and outline the extension of these result to vector-bundle valued operators in \S\ref{subsubsec_vecval}.

\paragraph{\textbf{Limitations and future work}}
The theory discussed in this paper is only the beginning of a systematic study of singularity formation using manifolds with corners and a-corners. We outline some directions for future work and possible limitations of some aspects:
\begin{itemize}
	\item The assumptions on the geometry of the space-time manifold $X$ in \S \ref{sec_heatoperators} are mainly for application to our intended problems as mentioned in \S \ref{sec_applications}. One would need to work with more general spaces as well as admitting higher codimensional corners for other problems. Similarly, the definition of a-smoothness is a natural weakening of the usual notion of smoothness which applies well to out intended applications, but one may need weaker notions of regularity for other problems. One may also need to allow for more general notions of fibration and more degenerate heat operators than those in \S \ref{sec_heatoperators}.
	
	\item It may happen that for a specific geometric problem there are no appropriate rates where the operator is Fredholm. Then, one has to further refine the analysis by \textit{working at a critical rate}. However, this is a feature, not a bug: for the problem to be well-behaved at a critical rate, we usually require that a certain projection of the leading order error term has to vanish, and this is what often gives us the characterising equation for the scale parameters in a gluing ansatz.
	\item We work mainly with $L^2$-Sobolev spaces in this paper, which makes for an easier (and slightly ad-hoc) approach to proving Fredholm theorems. However, ideally one would like to understand the Schwartz kernel of the heat operator on a blow-up resolution of $  X\times  X$ following Melrose, and this would allow one to directly prove H\"{o}lder estimates and $L^p$ estimates for all $1<p<\infty$.
\end{itemize}

\paragraph*{Acknowledgements}
The author is grateful to his D.Phil.\ supervisor Dominic Joyce for his guidance throughout the project, and for introducing him to the ideas in this paper. This work was funded by a Clarendon Scholarship from the University of Oxford and a Kingsgate Graduate Scholarship from Lincoln College.

\paragraph{AI disclosure} This paper was mostly written up by July 2026, without usage of LLMs. During the last stages of writing, ChatGPT Edu was used for proofreading and consistency checks.

\subsection*{Notation} We collect here some of the non-standard notations used in this paper for easy reference.
\begin{table}[h]
	\centering 
	\begin{tabular}{ p{0.2\linewidth} p{0.7\linewidth} }
		$[0,\infty),\; \llbracket0,\infty)$ & the topological space $[0,\infty)$ endowed with the standard smooth structure and the standard a-smooth structure respectively
		\\
		
		$\bR^s_k$ & the model space $[0,\infty)^k\times \bR^{s-k}$  
		\\
		
		$\bR^{k,s}$ & the model space $\llbracket0,\infty)^k\times \bR^{s-k}$
		
		\\
		$\bR^{k,s}_l$ & the model space $\llbracket0,\infty)^k\times [0,\infty)^l\times \bR^{s-k-l}$
		\\
		$\mancac$ & the category of manifolds with corners and a-corners 
		\\
	
		$\pa_v$ & a time-like vector field on $X$, i.e.\ $\pa_v(t)>0$ away from the a-boundary of $X$
		\\
        $\cV_t$ &{a local choice of non-degenerate time-like vector field}
        \\
        $\rho_1,\dots,\rho_N$ &{boundary-defining functions for the a-boundaries of $X$}
        \\
        $\brh^{\be}$ &{the weighted product of boundary-defining functions $\Pi_{i=1}^N\rho_i^{\eta_i}$}
        \\
		$C^\infty_0(X)$ & {the smooth functions that are compactly supported away from the a-boundaries of $X$}
		\\ $t$ & {the global time function on $X$}
		\\ $t'$ & {the local time function on $X$ near an a-corner, i.e.\ $t'=\pm t+c$}
		\\ $\pa^{\ord}X$ & {the ordinary boundary of $X$}
		\\ $\pa^{\ord}X^-, \pa^{\ord}X^+$ & {the ordinary boundaries of $X$ at initial time and at final time}
		\\ $\cT_\Xi$ & {the p.d.e.\ on model space obtained by freezing coefficients at the boundary}
		\\ $\cT_{\Xi,\deg}$ & {the p.d.e.\ on model space for vanishing  or vanishing corners}
        \\$N_M(\cT)(\xi)$&{the family of normal operators at an a-corner $M$ dependent on a parameter $\xi\in \bR^k$}
        \\$N_{\Xi,\deg}(\cT_{\Xi,\deg})(\xi)$ &{the family of degenerate normal operators at $\Xi$ dependent on a parameter $\xi\in \bR$ or $\xi\in \bR^2$}
		 
	\end{tabular}
\end{table}

\newpage
\section{Background}\label{sec_background}

\subsection{Manifolds with corners and a-corners}
In this subsection, we define the notions of manifolds with corners following Melrose \cite{Mel96} (see also \cite{Mel93}), and its generalization to manifolds with corners and a-corners following \cite{Joy16}, where more details can be found.

\begin{rem}
 The main reason for working with a-smooth structures is that many geometric problems admit solutions for which the usual notion of smoothness is too strong to be a natural regularity condition. For example, one may study the deformation theory of asymptotically conical shrinkers in $\bR^n$ by considering immersed manifolds with boundary inside the compactification of $\bR^n$ at infinity, satisfying a nonlinear p.d.e. Although the compactification of $\bR^n$ may be endowed with a natural smooth structure, the immersed submanifold will not be naturally smooth in general, but it will be a-smooth.
\end{rem}
\subsubsection{Manifolds with corners}
A manifold with corners is locally modelled on the space $\mathbb{R}^n_k:=[0,\infty)^k\times\bR^{n-k}$, with coordinates $(x_1,\dots,x_n)$ with $x_1,\dots,x_k$ in $[0,\infty)$ and $x_{k+1},\dots,x_n$ in $\bR$. Here the coordinates $\{x_1,\dots,x_k\}$ are also \textit{boundary-defining functions} for the boundaries $\{x_1=0,\dots,x_k=0\}$. This space has its smooth structure endowed from the inclusion $\bR^n_k\subset \bR^n$. We can define the notions of tangent bundle, cotangent bundle, metric, connections, etc. in a similar way as in ordinary manifolds. These structures may also be obtained by restricting them via inclusion in an ordinary smooth manifold.

We define the $b$-tangent bundle of a manifold with corners, which is a vector bundle whose sections are vector fields tangent to all the boundaries. This leads to associated notions of $b$-cotangent bundles, $b$-connections, $b$-metric and associated $b$-Levi-Civita connection. 

 There is an important class of mappings between manifolds with corners, called \textit{$b$-maps} as defined by Melrose \cite[\S 1.12]{Mel96}. The collection of manifolds with corners forms a category $\manc$ with morphisms being smooth $b$-maps.
\subsubsection{Manifolds with corners and a-corners}

A manifold with corners and a-corners is locally modelled on spaces $\mathbb{R}^{k,n}_l = \llbracket0,\infty)^k\times [0,\infty)^l\times \mathbb{R}^{n-k-l}$. We have coordinates $(x_1,\dots,x_n)\in \bR^{k,n}_l$ with $x_1,\dots,x_k$ in $\llbracket 0,\infty)$, $x_{k+1},\dots,x_{k+l}$ in $[0,\infty)$, and $x_{k+l+1},\dots, x_{n}$ in $\mathbb{R}$. Here, the coordinates $(x_1,\dots,x_k)$ are \textit{boundary-defining functions} for the \textit{a-boundaries} $\{x_1=0,\dots,x_k=0\}$, while the coordinates $(x_{k+1},\dots, x_{k+l})$ are boundary-defining functions for the \textit{ordinary boundaries} $\{x_{k+1}=0,\dots,x_{k+l}=0\}$. Near the a-boundaries, the smooth structure resembles that of a manifold with a-corners, while near the ordinary boundaries, the smooth structure resembles that of a manifold with corners.

First, we define the notion of an a-smooth function, following Definition 3.18 of \cite{Joy16}.

\begin{defn}\label{defn_asmoothfunc}
	Let $U\subset \mathbb{R}^{k,n}_l$ be an open set and $f: U\rightarrow \mathbb{R}$ be a continuous function. It is said to have a \textit{mixed derivative} or \textit{$\m$-derivative}, denoted $^{\m}\partial f$, if the function $^{\m}\partial f:U\rightarrow \mathbb{R}^n$, written as $^{\m}\partial f= (^{\m}\partial_1 f,\dots ,^{\m}\partial f_n)$ for $^{\m}\partial_if : U\rightarrow\mathbb{R}$ defined by
	\[
	^{\m}\partial_i f (x_1,\dots, x_n)= \begin{cases}
		0, & x_i=0, \: i=1,\dots,k \\
		x_i\partial_{x_i}f &  x_i>0,\: i=1,\dots,k \\
		\partial_{x_i}f &  i = k+1,\dots,n	
	\end{cases}
	\]
	(where we take one-sided derivatives for $x_i=0, i = k+1,\dots, k+l$) is a continuous function on $U$ for each $i=1,\dots, n$. We can iterate m-derivatives, if they exist, to get maps $^{\m}\pa^lf: U \rt \otimes^l\bR^n$ by taking m-derivatives of the components. Then,
	\begin{itemize}
		\item we say that $f$ is \textit{roughly differentiable} or \textit{r-differentiable} if $^{\m}\pa f$ exists as a continuous function, and \textit{$r$-smooth} if all the m-derivatives exist for $l=0,1,\dots$.
		\item we say that $f$ is \textit{a-differentiable}, if it is r-differentiable and we have that for any compact $S\subseteq U$ and $i=1,\dots, k$ , there are positive constants $C,\alpha$ such that \[
		|^{\m}\pa_if(x)| \leq Cx_i^\alpha \; \;\:\; \textrm{for all} \: (x_1,\dots,x_n)\in S.
		\]
		\item we say that $f$ is \textit{a-smooth} if the m-derivatives $^{\m}\pa ^lf$ are all a-differentiable for $l=0,1,\dots$.
	\end{itemize}	
\end{defn}

We note that a-smooth functions are closed under addition and multiplication, and the quotient of an a-smooth function by a non-vanishing a-smooth function is a-smooth. Also, we note that there exist functions which are r-smooth but not a-smooth, for example $f: \llbracket0,1/2) \rt \bR$ given by $f(x)= \log(x)^{-1}$, since $|^{\m}\pa ^lf| = O(|\log(x)|^{-{l-1}})$ which does not decay fast enough.

Now, we define a-smooth maps between subsets of the model spaces, similar to Definition 3.19 in \cite{Joy16} (however, there are some key differences regarding the exponents that one allows in parts (a),(c) below -- we allow mapping an a-boundary to an ordinary boundary and disallow mapping an ordinary boundary to an a-boundary, which is to ensure that a-smooth maps are closed under composition). Note that this notion is analogous to the $b$-maps of Melrose.

\begin{defn}\label{defn_asmoothmap}
	Let $U\subseteq \mathbb{R}^{k,m}_l$, $V\subseteq \mathbb{R}^{p,n}_q$ be open subsets, and $f=(f_1,\dots, f_n): U\rt V$ be a continuous map, so that $f_i:U\rt \llbracket 0,\infty)$ for $i=1,\dots p$, $f_i : U\rt [0,\infty)$ for $i = p+1,\dots, p+q$, and $f_i : U \rt \bR$ for $i=p+q+1,\dots, n$. Then,
	\begin{enumerate}
		\setcounter{enumi}{0}
		\item we say $f$ is \textit{a-smooth} if $f_j$ is a-smooth for $j=p+q+1,\dots, n$, and every $u = (x_1,\dots, x_m)\in U$ has an open neighbourhood $\tilde{U}\subseteq U$ such that for each $j=1,\dots, p$, we have either:
		\begin{enumerate}
			\item we can uniquely write \[f_j(\tilde{x}_1,\dots,\tilde{x}_m) = F_j(\tilde{x}_1,\dots,\tilde{x}_m)\tilde{x}_1^{a_{1,j}}\cdots \tilde{x}_{k}^{a_{k,j}} \]
			for all $(\tilde{x}_1,\dots,\tilde{x}_m)\in \tilde{U}$, where $\tilde{F}_j$ is a positive a-smooth function, and $a_{1,j},\dots, a_{k,j} \in [0,\infty)$ with $a_{i,j}= 0$ if $x_i\neq 0$; or
			\item $f_j|_{\tilde{U}}=0$,
		\end{enumerate}  
		and for each $j=p+1,\dots, p+q$, we have either:
		\begin{enumerate}
			\setcounter{enumii}{2}
			\item we can uniquely write \[f_j(\tilde{x}_1,\dots,\tilde{x}_m) = F_j(\tilde{x}_1,\dots,\tilde{x}_m)\tilde{x}_{1}^{b_{1,j}}\cdots \tilde{x}_{k}^{b_{k,j}}\cdot\tilde{x}_{k+1}^{b_{k+1,j}}\cdots \tilde{x}_{k+l}^{b_{k+l,j}} \]
			for all $(\tilde{x}_1,\dots,\tilde{x}_m)
			\in \tilde{U}$, where $F_j$ is a positive a-smooth function, $b_{1,j},\dots,b_{k,j} \in [0,\infty)$, and $b_{k+1,j},\dots, b_{k+l,j} \in \bZ_{\geq 0}$ with $b_{i,j}= 0$ if $x_i\neq 0$; or
			\item $f_j|_{\tilde{U}}=0$,
		\end{enumerate}  
		\item we say $f$ is \textit{interior} if it is a-smooth and cases (b) or (d) in the definition do not occur.
		\item we say $f$ is \textit{$b$-normal} if it is interior, and in cases (a) and (c) in the definition, for each $i=1,\dots, k+l$ we have $a_{i,
			j}>0$ or $b_{i,j}>0$ for at most one $j=1,\dots, p+q$ (taking $a_{i,j}=0$ for $j=p+1,\dots,p+q$).
		\item we say $f$ is an \textit{a-diffeomorphism} if it is a-smooth with an a-smooth inverse.
	\end{enumerate} 
\end{defn}
We note that a-smooth functions defined as above contain identities and are closed under composition. Now, we can define manifolds with corners and a-corners:
\begin{defn}
	A \textit{manifold with corners and a-corners} is defined to be a second-countable Hausdorff topological space, with charts to the model spaces $\bR^{k,m}_l$, such that the transition maps are a-diffeomorphisms. We denote the category of manifolds with corners and a-corners, with morphisms a-smooth maps, by $\mancac$. Note that we have a fully faithful inclusion of categories $\manc\subset \mancac$.
\end{defn}

Now, we recall the notion of ordinary boundaries and a-boundaries of a manifold with corners and a-corners:

\begin{defn}\label{defn_acorners}
	Let $X$ be a manifold with corners and a-corners. Then, it has a natural stratification, $X= \sqcup_{k=0}^{\dim X} S^k(X)$, where $S^k(X)$ are open manifolds without corners of codimension $k$. A \textit{local $k$-corner component $\gamma$ of $X$ at a point $x$} is a choice of connected component $\gamma \in S^k(X)$ such that $x$ is contained in the closure of $\gamma$. The case $k=1$ corresponds to a \textit{local boundary component}. As sets, we define the \textit{boundary} and \textit{codimension $k$ corners} of $X$ by:
	\[
	\pa X := \{ (x,\beta): \beta \textrm{ is a local boundary component of $X$ at $x$}\},
	\]
	\[
	M_k(X):= \{(x,\beta): \beta \textrm{ is a local $k$-corner component of $X$ at $x$}\}.
	\]
	The boundary $\pa X$ of $X$ admits a natural partition, $\pa X = \pa^{\textrm{o}}X \sqcup \pa^{\textrm{a}}X$, where $\pa^{\textrm{o}}$ is the \textit{ordinary boundary} of $X$ and $\pa^{\textrm{a}}X$ is the \textit{a-boundary} of $X$. These are defined as follows: in a chart $(U,\phi)$ for $X$ with $U\subset \bR^{p,n}_q$, the interior of $\pa X$ consists of points of the form $\phi(x_1,\dots,x_n)$ with exactly one of $\{x_1,\dots,x_{p+q}\}$ being zero, while the interior of $\pa^\textrm{o}X$ consists of $\phi(x_1,\dots,x_n)$ with exactly one of $\{x_{p+1},\dots,x_{p+q}\}$ being zero, and the interior of $\pa^{\textrm{a}}X$ consists of $\phi(x_1,\dots,x_n)$ with exactly one of $\{x_1,\dots,x_p\}$ being zero. We similarly define the \textit{codimension $k$ a-corners} of $X$ as $M^{\ma}_k(X)\subset M_k(X)$ such that in a chart $U\in \bR^{p,n}_q$, the interior consists of $\phi(x_1,\dots,x_n)$ with exactly $k$ of $\{x_1,\dots,x_p\}$ being zero, i.e.\ these are the corners which are locally the intersection of a-boundaries. We say that the boundaries of $X$ are \textit{faces} if the inclusion $\iota_X|_C: C\rt X$ of each connected component $C$ of the boundary $\pa X$ is injective. 
\end{defn}

\begin{ex}
	Consider the example of a manifold with corners and a-corners given by $X = [-1,1]\times \llbracket-1,1\rrbracket$. This has boundary the disjoint union of four pieces, given by $\pa X = (\{-1\}\times \llbracket-1,1\rrbracket) \sqcup (\{1\}\times\llbracket-1,1\rrbracket) \sqcup ([-1,1]\times \{-1\} )\sqcup ([-1,1]\times \{1\})$. The a-boundary is given by $([-1,1]\times \{-1\} )\sqcup ([-1,1]\times \{1\})$, while the ordinary boundary is given by $(\{-1\}\times \llbracket-1,1\rrbracket) \sqcup (\{1\}\times\llbracket-1,1\rrbracket)$. The codimension $2$ corner is the disjoint union of $4$ points, while there are no codimension $2$ a-corners, and all the boundaries are faces.
\end{ex}

\subsubsection{Vector bundles, connections and metrics}
We define the notion of a vector bundle on elements of $\mancac$ as in Definition 4.4 of \cite{Joy16}:
\begin{defn}
	Let $X$ be an element of $\mancac$. A \textit{vector bundle of rank $k$} over $X$ is an object $E\in \mancac$ with an a-smooth map $\pi : E\rt X$ such that each fibre $E_x:=\pi^{-1}(x)$ has the structure of a vector space, and there is a trivialising cover $\{U_i\}_{i\in I}$ of $X$ such that $\pi^{-1}(U_i) \cong U_i\times \bR^k$ and $\pi|_{\pi^{-1}(U_i)}$ is the projection $ U_i\times \bR^k \rt U_i$. A \textit{section} of $E$ is an a-smooth map $s: X\rt E$ such that $\pi\circ s = \mathrm{id}_X$, and we denote the vector space of a-smooth sections by $\Gamma^\infty(E)$.
\end{defn}
Now, we define the notion of the $b$-tangent bundle and the mixed tangent bundle:
\begin{defn}
	Let $X$ be an $m$-dimensional manifold with corners and a-corners. We define the \textit{mixed tangent bundle} or \textit{${\m}$-tangent bundle} of $X$ as follows: it is a rank $m$ vector bundle over $X$, whose sections on a local chart $U\subset \bR^{k,m}_l$ are spanned by $\{\{x_i\frac{\partial}{\partial x_i}\}_{i=1}^{k}\sqcup \{\frac{\partial}{\partial x_i}\}_{i=k+1}^m\}$, for $\{x_i\}_{i=1}^k$ being the boundary-defining coordinates for the a-boundary. We denote it by $^{\m}TX$.
\end{defn}
\begin{defn}
	Let $X$ be an $m$-dimensional manifold with corners and a-corners. We define the \textit{$b$-tangent bundle} of $X$ as follows: it is a rank $m$ vector bundle over $X$, whose sections on a local chart $U\subset \bR^{k,m}_l$ are spanned by $\{\{x_i\frac{\pa}{\pa x_i}\}_{i=1}^{k+l} \sqcup \{\frac{\pa}{\pa x_i}\}_{i=k+l+1}^m\}$, for $\{x_i\}_{i=1}^k$ being the boundary-defining functions for the a-boundaries and $\{x_i\}_{i=k+1}^{k+l}$ being the boundary-defining functions for the ordinary boundaries. We denote it by $^bTX$. Note that there is a natural inclusion $\cV^b=\Gamma^\infty(^bTX) \hookrightarrow \cV^\m=\Gamma^\infty({^{\m}TX})$.
\end{defn}

We can define the $b$-cotangent bundle and m-cotangent bundle of $X$, which are duals to the corresponding tangent bundles. As in the case of ordinary vector bundles, we can form new bundles by taking duals, tensor products, and direct sums. Now, we define the notion of an m-metric on $X$:
\begin{defn}
	An \textit{${\m}$-metric} $g$ on $X$ is a section of $\mathrm{Sym}^2(^{\m}T^*X)$, such that it is positive-definite at each point of $X$.
\end{defn}
\begin{rem}
	Note that if we consider the interior of $X$ as an ordinary manifold $X^\circ$, the m-metric restricts to an ordinary Riemannian metric, such that the a-boundaries are `at infinity' in $X^\circ$, while the ordinary boundaries are at a `finite distance'.
\end{rem}
We define the notion of $b$-connections and m-connections on a vector bundle:
\begin{defn}
	A \textit{$b$-connection} on a vector bundle $E\rt X$ is an $\bR$-linear map $^{b}\nabla : \Gamma^\infty(E)\rt\Gamma^\infty(E\otimes {^{b}}T^*X)$, satisfying \[
	^{b}\nabla(c\cdot e) = c\cdot {^{b}\nabla e} + e\otimes {^{b}d c},
	\]
	for $e\in \Gamma^\infty(E)$ and $c\in C^\infty(X)$. Similarly, an \textit{${\m}$-connection} is an $\bR$-linear map $^{\m}\nabla : \Gamma^\infty(E)\rt\Gamma^\infty(E\otimes {^{\m}T}^*X)$, satisfying \[
	^{\m}\nabla(c\cdot e) = c\cdot {^{\m}\nabla e} + e\otimes {^{\m}d c}.
	\]
\end{defn}
Given an m-metric, we can define a natural m-connection on $^{\m}TX$ analogously to the Levi-Civita connection, as in Proposition 5.2$.$ of \cite{Joy16}. This restricts to the ordinary Levi-Civita connection on $X^\circ$ with respect to the Riemannian metric on $X^\circ$ given by the restriction of $g$ on $X^\circ$. 

Now, we recall the notion of derivatives of interior mappings, following Definition 4.6 of \cite{Joy16}:
\begin{defn}
	Given an interior a-smooth map $f:X\rt Y$, it induces a well-defined mapping $^\m df: {^\m TX}\rt f^*{^\m TY}$ called the \textit{mixed derivative} or m-derivative, and a mapping $^bdf: {^bTX}\rt {f^*(^bTY)}$, called the \textit{$b$-derivative}.
\end{defn}

\begin{rem}
	The $\m$-derivative of a map $g:\llbracket0,\infty)\rt \llbracket0,\infty)$ is given by mapping $r\pa_r\mapsto {^\m\pa}\log(g(r))\cdot t\pa_t$, so the coefficient is \textbf{not} taking the $\m$-derivative $^\m \pa g$ as a real-valued function in Definition \ref{defn_asmoothfunc}. 
\end{rem}
\begin{defn}
	Given a boundary $\Sigma \in M_1(X)$, a \textit{boundary-defining function} is an a-smooth map $\rho: X\rt \llbracket0,\infty)$ in case $\Sigma$ is an a-boundary or $\rho: X\rt [0,\infty)$ in case $\Sigma$ is an ordinary boundary, such that $\rho^{-1}(0)=\Sigma$ and ${^\m}d \rho|_{\Sigma} \neq 0$.
\end{defn}

Now, we define the important notion of a $b$-fibration:
\begin{defn}
	Let $f:X\rightarrow Y$ be an a-smooth, interior map between manifolds with corners and a-corners. It is called a \textit{$b$-submersion} if the induced $b$-derivative is surjective, and a \textit{$b$-immersion} if the induced $b$-derivative is injective. It is called is a \textit{$b$-fibration} if it is a $b$-submersion and $b$-normal.
\end{defn}
We have the following important result from \S 2.7 of \cite{Mel96}:
\begin{prop}\label{prop_bfibr}
Given a $b$-fibration $f:X\rt Y$, the kernel $\ker({^bd f})$ is a well-defined sub-bundle of $^bTX$, called the \textit{relative $b$-tangent bundle} and denoted $^bT_\rel X$.
\end{prop}

\subsubsection{Relative bundles over an m-fibration}

For our purposes, we require a refinement of a $b$-fibration. We will often have to work with a \textit{time variable} which is locally given by a mapping $t: [0,\infty)\times \llbracket0,\infty)\rt [0,\infty)$, $t(x,y)=xy$. This is an interior map, so the $b$-derivative sends:\[
^bd t(x\pa_x)=t\pa_t,\;\; ^bd t(y\pa_y)=t\pa_t,
\]
while the m-derivative sends\[
^\m d t(\pa_x)=y\pa_t\;\; ^\m d t(y\pa_y)=xy\pa_t.
\]
Clearly, this is a $b$-fibration as the $b$-derivative is surjective. However, the m-derivative is surjective only away from the a-boundary $y=0$. If we define a notion of `m-fibration' by declaring the m-derivative to be surjective everywhere, we would be disallowing any maps from an a-boundary to an ordinary boundary, which seems too strong. Instead, we define the following:
\begin{defn}\label{defn_mfibration}
	A \textit{mixed fibration} or \textit{$\mathrm{m}$-fibration} is a $b$-fibration $f:X\rt Y$ such that the m-derivative $^\m d f$ is surjective away from the a-boundaries of $X$. Equivalently, it is a $b$-fibration such that in the local charts of Definition \ref{defn_asmoothmap}, for each $j=p+1,\dots,p+q$, the unordered collection of exponents $\{b_{k+1},\dots,b_{k+l}\}$ is either $\{0,\dots,0\}$ or $\{1,0,\dots,0\}$.
\end{defn}
	For an m-fibration a similar result to \ref{prop_bfibr} holds:
\begin{prop}\label{prop_mfibre}
	Given an m-fibration $f:X\rt Y$, the relative tangent bundle over the interior of $X$, $^b T_\rel X^\circ \subset {^\m T X^\circ}$ extends over the boundary to give a vector sub-bundle of $^\m TX$ called the \textit{relative $\m$-tangent bundle}, denoted $^\m T_\rel X$.  
\end{prop}
\begin{proof}	
	This follows by a similar calculation as for that of showing $^bT_\rel$ is a sub-bundle, see equation 2.4.10 in \S 2.4 and \S 2.7 of \cite{Mel96}.
	In particular, we can write down non-vanishing local sections spanning the $b$-tangent bundle, and show that each of them provides a non-vanishing section of the relative m-tangent bundle. We introduce local coordinates $(x_1,\dots,x_n)$ on $X$ and $(y_1,\dots,y_m)$ on $Y$ such that\[
	y_j=x_1^{a^j_{1}}\cdots x_{p+q}^{a^j_{p+q}} \;\forall j\in\{1,\dots,k+l\},\; \; y_j=x_{p+q-k-l+j} \;\forall j\in\{k+l+1,\dots,m\}.
	\]
    By surjectivity of $^bd$ and $b$-normality, it follows that for each $i=1,\dots,p+q$, $a_i^j>0$ for exactly one $j$.
	Now the m-tangent bundle on $X$ is spanned in these coordinates by $\{x_i\pa_{x_i}\}_{i=1}^k\sqcup\{\pa_{x_i}\}_{i=k+1}^n$. The vector fields spanning $^bT_\rel X^\circ$ over the interior are of the form $a_i^jx_k\pa_{x_k}-a_i^kx_j\pa_{x_j}$ for $i\in\{1,\dots, k+l\}$ and $j<k$ or $\pa_{x_j}$ for $j\in\{p+q+1,\dots, p+q+m-k-l\}$. Clearly, these sections extend to non-vanishing sections of $^\m T_\rel X$ except for the first case when $j,k\in\{p+1,\dots,p+q\}$. In this case, by Definition \ref{defn_mfibration}, either $a_i^k=0$ or $a_i^j=0$, so we can write the vector field as $a_i^sx_s\pa_{x_s}$, for $s=j\textrm{ or }k$. Over the interior, this is the same as the span $\la\pa_{x_s} \ra$, so it extends to a non-vanishing section of $^\m T_\rel X$.
	\end{proof}

\begin{ex}
	Consider the map at the beginning $f: [0,\infty)\times\llbracket0,\infty)\rt [0,\infty)$ mapping $(x,y)\rt t=xy$. Then $^\m TX \cong \la \pa_x, y\pa_y\ra$ and \[
	^\m d f (\pa_x)=y\pa_t,\;\; ^\m d f (y\pa_y)=xy\pa_t.
	\]
	So, $\ker(^\m d f)=\{a\pa_x+by\pa_y: y(a+bx)=0\}$. This is a $1$-dimensional subspace for $y\neq0$, but a $2$-dimensional subspace for $y=0$, and thus the kernel fails to be a sub-bundle of $^\m TX$. However, over the interior we can write $^\m T_\rel X^\circ=\la x\pa_x-y\pa_y\ra$, and $x\pa_x-y\pa_y\in \Gamma^\infty(^\m TX)$ is a non-vanishing section because the component of $y\pa_y$ is non-vanishing. Thus it extends over the boundary to give a well-defined sub-bundle of $^\m TX$.
\end{ex}

\begin{ex}
 Consider the mapping, $f:[0,\infty)^2\rt [0,\infty)$ that maps $(x,y)\mapsto t=xy$. This is a $b$-fibration; however, the m-tangent bundle is $^\m TX\cong \la\pa_x,\pa_y\ra$ and the m-differential maps:
	\[
	{^\m d f}(\pa_x) = y\pa_t,\; 	{^\m d f}(\pa_y) = x\pa_t.
	\]
	Therefore, the image of $^\m d f$ is $\{(ay+bx)\pa_t\}_{a,b\in \bR}$, and so, $\ker(^\m d f)(x,y)=\{a\pa_x+b\pa_y: ay+bx=0\}$. This is a subspace of dimension $1$ for $(x,y)\neq(0,0)$ but a subspace of dimension $2$ for $(x,y)=(0,0)$, so the dimension of the kernel jumps at the origin and $\ker(^md f)$ fails to be a sub-bundle of $^\m TX$. In fact, we cannot extend $^\m T_\rel X^\circ$ over $X$, so it is not an m-fibration.
\end{ex}

As a result, we have an injection of vector bundles, \[
0\rt {^\m} T_\rel X \rt {^\m} TX ,
\]
and taking duals gives us a surjection,\[
^\m T^* X \rt {^\m} T^*_\rel X \rt 0.
\]
\begin{rem}
    The above exact sequence does not extend to $\cdots\rt f^*TY\rt 0$, because as we have noted, the m-derivative may fail to be surjective at a-boundaries.
\end{rem}

\begin{defn}
	Given an m-fibration $f:X\rt Y$ we define the \textit{relative differential} as a map \[
	^\m d_\mathrm{rel} : C^\infty(X) \rt \Gamma^{\infty}(^\m T^*_\mathrm{rel}X),
	\]
	by composing the $\m$-differential $^\m d$ with the dual to the inclusion map $^\m T_\rel X\rt {^\m T}X$.
\end{defn}

Then we can define \textit{relative m-connections} on a vector bundle $E\downarrow X$ over an m-fibration as an $\bR$-linear map $^{\m}\nabla_\rel : \Gamma^\infty(E)\rt\Gamma^\infty(E\otimes {^{\m}T_\rel^*X})$, satisfying \[
^{\m}\nabla_\rel(c\cdot e) = c\cdot {^{\m}\nabla_\rel e} + e\otimes {^{\m}d_\rel c},
\]
for $e\in \Gamma^\infty(E)$ and $c\in C^\infty(X)$. Similarly, we can define a \textit{relative m-metric} as section of $\mathrm{Sym}^2(^{\m}T^*_\rel X)$, such that it is positive-definite at each point of $X$.

\subsubsection{Partial differential operators}
Let $E,F\downarrow X$ be vector bundles, $\nabla^E$ be an m-connection on $E$ and $\nabla$ be an m-connection on $^\m TX$.
\begin{defn}
	A \textit{linear partial differential operator} of order $l \in \bZ_{\geq 0}$ is an $\bR$-linear mapping,\[
	\cP : \Gamma^\infty(E)\rt \Gamma^\infty(F),
	\]
	which can be written as \[
	\cP= \sum_{j=0}^la_j(\nabla^E)^j,
	\] for $\nabla^E$ being a connection on $E$, and $\nabla^E$ for higher orders is defined by combining $\nabla^E$ with $\nabla$, and $a_j\in \Gamma^\infty(E^*\otimes (\bigotimes^j {^\m TX})\otimes F)$. A \textit{nonlinear partial differential operator} of order $l$ is a mapping \[
	\cP: \Gamma^\infty(E) \rt \Gamma^\infty(F)
	\]
	which can be written as $\cP = \cQ(\sigma,\nabla^E\sigma,\dots,(\nabla^E)^l\sigma)$, for \[
	\cQ:\bigoplus_{j=0}^l (E\otimes \Sym^j({^\m}T^*X)) \rt F
	\]
	being a smooth mapping between vector bundles (i.e. $\cQ: E'\rt F$ is smooth and maps $\cQ|_{x}: E'_{x}\rt F_x$ are smooth).
\end{defn}

\begin{defn}
	The \textit{principal symbol} of a linear partial differential operator of order $l$ is given for each $x\in X, \xi \in {^\m T^* X}$ by \[
	\sigma_{\cP}(x,\xi):E_x\rt F_x, \;\; \sigma_\cP(x,\xi)(e_x):= a_l|_x\cdot(e_x\otimes^l\xi).
	\] 
\end{defn}
\begin{defn}
	We say that $\cP$ is \textit{elliptic} if the principal symbol of its linearisation is an isomorphism for $x\in X$, $\xi \neq 0$.
\end{defn}
For our purposes, we define a strong notion of ellipticity following Chapter 4 of the book by Topping \cite{Top06}:
\begin{defn}
	Let $\cP:\Gamma^\infty(E)\rt \Gamma^\infty(E)$ be an elliptic partial differential operator between sections of the same vector bundle. Then, it is defined to be \textit{strongly elliptic} if there exists an inner product on $E$ such that \[
	\la \sigma_{\cP}(x,\xi) v,v\ra_E\geq \lambda|\xi|^2|v|_E^2, 
	\]
	for $\lambda>0$ and for all $\xi \in {^\m T^*_{x} X}, v\in E_x$.
\end{defn}

Now, we define the restrictions of nonlinear differential operators at a-corners:
\begin{defn}
	Given an a-corner $\Sigma \in M^{\ma}_k(X)$ (as in Definition \ref{defn_acorners}), we can define the \textit{restriction of $\cP$ at $\Sigma$} as follows: it is a map \[
	\cP|_{\Sigma}: \Gamma^\infty(E|_{\Sigma}) \rt \Gamma^\infty(F|_{\Sigma}),
	\]
	defined as $\cP|_{\Sigma}(\sigma):=\cP(\tsig)|_{\Sigma}$, where $\tsig \in \Gamma^\infty(E|_{U})$ is a section over a neighbourhood $U$ of $\Sigma$ such that $\tsig|_{\Sigma}=\sigma$. 
\end{defn}
\begin{rem}
	It is necessary that $\Sigma$ is an a-corner for the restriction to be well-defined, otherwise the restriction $\cP(\tsig)|_\Sigma$ would depend on the choice of extension $\tsig$.
\end{rem}
Then, we have that the restriction of a nonlinear partial differential operator of order $l$ at an a-boundary $\Sigma$ 
is also a differential operator of order $l$. If $\cP$ is elliptic, so is $\cP|_{\Sigma}$.
\subsubsection{Relative differential operators}\label{subsubsec_reloper}
Now let $E,F\downarrow X$ be vector bundles as before, $f:X\rt Y$ be an m-fibration, and $\nabla_\rel^E$ be a relative connection on $E$ with respect to the m-fibration, and $\nabla_\rel$ be a relative m-connection on $^\m T_\rel X$. Then, as before we can define linear \textit{relative partial differential operators} as in \S 5.4 of \cite{Joy16}, given by a mapping \[
\cP : \Gamma^\infty(E)\rt \Gamma^\infty(F),
\]
which can be written as \[
\cP= \sum_{j=0}^la_j(\nabla_\rel^E)^j,
\] with $\nabla_\rel^E$ for higher orders defined by combining $\nabla^E$ with $\nabla_\rel$, and $a_j\in \Gamma^\infty(E^*\otimes (\bigotimes^j {^\m T_\rel X})\otimes F)$.
Then exactly as before, we can define \textit{nonlinear relative differential operators} and \textit{(strongly) elliptic relative differential operators}, by replacing instances of $^\m TX$ with $^\m T_\rel X$. 
\subsection{Interior estimates and Sobolev embeddings}
In this section, we recall some parabolic interior estimates and embedding theorems we will need, mainly following Krylov's books \cite{Kry08,Kry96}. 
\subsubsection{Parabolic Sobolev spaces}
We consider our model space as $\bR^{n+1}_0:=\{(x,t):x\in\bR^{n}, t\in [0,\infty)\}$, where $x$ represents the space variable and $t$ represents the time variable. 
\begin{defn}
	For a differentiable function $u$ having continuous derivatives of the form $\pa^k_{{i_1\dots i_k}}u$ for $(i_1,\dots, i_k)\in \{1,\dots, n\}$, we denote by $\nabla^ku$ the tensor given by \[
	\nabla^ku:= \sum_{1\leq i_1\leq \dots i_k\leq n} (\pa^k_{i_1\dots i_k}u)dx^{i_1}\dots dx^{i_k}.\]
	Similarly, for a function $u$ having continuous derivatives of the form $\pa_t^l\pa^k_{i_1\dots i_k}u$, we denote by $\pa_t^l\nabla^ku$ the tensor given by \[
	\pa_t^l\nabla^ku:=\sum_{1\leq i_1\leq \dots i_k\leq n} (\pa_t^l\pa^k_{i_1\dots i_k}u)dx^{i_1}\dots dx^{i_k}.
	\]
\end{defn}
Now, we define the $L^2_{k+2r,r}$-norms on our model space.
\begin{defn}
	Let $k,r$ be non-negative integers. For a locally integrable function $u$ having locally integrable weak derivatives of the form $\pa_t^l\pa^m_{i_1\dots i_m}u$ for all $l\leq r$ and $m+2l\leq k+2r$, we define \[
	\|u\|^2_{L^2_{k+2r,r}}:= \sum_{i\leq r}\sum_{j\leq k+2r-2i}\|\pa_t^i\nabla^ju\|^2_{L^2}.
	\]
\end{defn}

We define the parabolic Sobolev spaces on $\bR^{n+1}_{0}$ as follows:
\begin{defn}
	The Sobolev space $L^2_{k+2r,r}$ is defined to be the space of locally integrable functions $u$ having locally integrable weak derivatives of the form $\pa_t^l\pa^m_{i_1\dots i_m}u$ for all $l\leq r$ and $m+2l\leq k+2r$ and having a finite $L^2_{k+2r,r}$-norm.
\end{defn}
\subsubsection{Parabolic H\"{o}lder spaces}\label{subsec_holdereuclidean}
Let us denote the parabolic distance function $\rho(z_1,z_2)$ between two points $z_1= (x_1,t_1)$ and $z_2 = (x_2,t_2)$ as $\rho(z_1,z_2):= \sqrt{|x_1-x_2|^2+|t_1-t_2|}$. Then, we define the parabolic H\"{o}lder spaces:
\begin{defn}
	Let $Q\subset \bR^{n+1}_0$ be a domain. Then, we denote the norm \[
	[u]_{\delta,\delta/2;Q} := \sup_{z_1\neq z_2, \; z_i\in Q} \frac{|u(z_1)-u(z_2)|}{\rho^\delta(z_1,z_2)}.
	\] 
	We define the H\"{o}lder space, for non-negative integers $k,r\geq 0$ as the space of functions having continuous $C^{k+2r,r}$-derivatives and finite norm\[
	|u|_{C^{k+2r+\delta, r+\delta/2}(Q)} := |u|_{C^{k+2r,r}(Q)}+\sum_{\alpha \leq r,\;2\alpha+|\beta|=k+2r} [\pa_t^\alpha\nabla^\beta u]_{\delta,\delta/2; Q}.
	\]
\end{defn}
\subsubsection{Traces in parabolic Sobolev spaces}

Consider the space $\bR^{n+1}_{\geq0}$. We note that there exists a bounded linear mapping, \[
T: L^2_1(\bR^{n+1}_{\geq0}) \rt L^2(\bR^n), 
\]
such that for $u\in L^2_1\cap C^0$, we have $T(u)=u|_{t=0}$. This mapping is called the trace, see for example \S 5.5 of Evans' book \cite{Eva10}. We can extend the trace map to parabolic Sobolev spaces as in \S 2.5 of \cite{Kry08}, to define a bounded linear map\[
T: L^2_{k+2r,r}(\bR^{n+1}_{\geq 0})\rt L^2_{k+2r-1}(\bR^n).
\]
We also obtain a surjectivity result, following Krylov:
\begin{lem}
	Let $u\in L^2_1(\bR^n)$. Then, there exists a function $v\in L^2_{2,1}(\bR^{n+1}_{\geq0})$ such that $\|v\|_{L^2_{2,1}}\leq C\|u\|_{L^2_1}$ for a constant $C$ independent of $u$, and $T(v)=u$. 
\end{lem}
\begin{proof}
	We define $v$ to be the solution to the heat equation \[
	(\pa_t-\Delta+1)v=0, \; v|_{t=0}=u.
	\]
	Then, taking Fourier transforms in space, we get for each $\xi\in \bR^n$,\[
	(\pa_t+|\xi|^2+1)\hat{v}(\xi)(t)=0,\; \hat{v}(\xi)(0)=\hat{u}(\xi).  
	\]
	This has unique solution $\hat{v}(\xi)(t)= e^{-(1+|\xi|^2)t}\hat{u}(\xi)$, from which we get 
	\begin{equation*}\begin{split}
			\|\nabla_x^2v\|_{L^2} &= \||\xi|^2e^{-(1+|\xi|^2)t}\hat{u}(\xi)\|_{L^2} = \int_{\bR^n}|\xi|^4\hat{u}^2(\xi) \int_0^\infty e^{-2(1+|\xi|^2)t}dtd\xi \\
			&= \int_{\bR^n}\frac{|\xi|^4}{2(1+|\xi|^2)}\hat{u}^2(\xi)d\xi \leq \frac{1}{2}\|\nabla_x u\|^2_{L^2}.
		\end{split}
	\end{equation*}
	Similarly, we can show that $\|v\|_{L^2_{2,1}}\leq C\|u\|_{L^2_1}$.
\end{proof}
Now, we can extend the above to domains of the form $\Omega\times[0,\ep)$, to get:
\begin{lem}
	Let $\Omega \subset \bR^n$ be a bounded domain with smooth boundary. Let $u\in L^2_1(\Omega)$. Then, there exists a function $v\in L^2_{2,1}(\Omega\times[0,\ep))$ such that $\|v\|_{L^2_{2,1}}\leq C\|u\|_{L^2_{1}}$, for a constant $C$ depending only on $\Omega$, and $T(v)=u$.
\end{lem}
\begin{proof}
	By Stein's extension theorem, we can first extend $u$ to a function $\tu$ on all of $\bR^n$. Then, we can use the previous lemma to get a function $v$ on $\bR^{n+1}_{\geq0}$, then restrict it to $\Omega\times[0,\ep)$.
\end{proof}

Therefore, we obtain:
\begin{lem}\label{lem_extensionoftrace}
	Let $\Omega \subset \bR^n$ be a bounded domain with smooth boundary. Let $u\in L^2_{k+2r-1}(\Omega)$. Then, there exists $v\in L^2_{k+2r,r}(\Omega\times[0,\ep))$ and a constant $C$ depending only on $\Omega$ such that \[
	\|v\|_{L^2_{k+2r,r}}\leq C\|u\|_{L^2_{k+2r-1}},\; v|_{t=0}=u.
	\]
\end{lem}
\begin{proof}
	Let us construct $v$ as a solution to the heat equation as before. Then, the estimates follow since the derivatives $\nabla^{k+2r-j}_xv$ for $2\leq j\leq k+2r$ satisfy the same heat equation, with initial data $\nabla^{k+2r-j}_xu$. 
\end{proof}

\subsubsection{Interior estimates}\label{subsec_interior}
Let $Lu := a^{ij}\partial^2_{ij}u + b^i \partial_i u +cu$ be a second order, elliptic partial differential operator defined on a domain $\Omega \subset \mathbb{R}^n$ such that \[
\kappa |\xi|^2 \geq \sum_{i,j}a^{ij} \xi_i\xi_j \geq \kappa^{-1}|\xi|^2
\]  
for a fixed constant $\kappa$ for all $\xi \in \mathbb{R}^n$, and the $C^{k+2r-2+\delta,r-1+\delta/2}$-norms of $a^{ij}, b^{i},c$ are bounded in $\Omega$ by a constant $\Lambda$. We will require two kinds of domains for our interior estimates:
\begin{enumerate}
	\item Let $B_1 \ssubset B_2$ be two concentric Euclidean balls in $\Omega$ with radii $r_2>r_1>0$, and $I_i:= (-\ep_i,0)$ for $i=1,2$ be two open intervals, such that $\ep_1<\ep_2$. We denote $Q_i:=B_{i}\times I_i$ for $i=1,2$, so that $Q_1\subset Q_2$. These will be used for estimates away from the ordinary boundaries at initial time.
	\item Let $B_1\ssubset B_2$ be two open balls as before and $I_i:= [0,\ep_i)$ for $i=1,2$ be two half-open intervals, such that $\ep_1<\ep_2$. We denote $P_i:=B_{i}\times I_i$ for $i=1,2$, such that $P_1\subset P_2$. These will be used for estimates near the ordinary boundaries at initial time.
\end{enumerate} 
On each of these, we state the following $L^2$-Sobolev and H\"{o}lder estimates we will require:
\begin{theorem}\label{thm_ffestimate}
	Let $u \in L^2_{k+2r,r}(Q_2)$ and let us denote $h:= (\pa_t - L)u$. Then, we have that
	\begin{equation}\label{eq_ffestimate}
		\|u|_{Q_1}\|_{L^2_{k+2r,r} } \leq C\|h|_{Q_2}\|_{L^2_{k+2r-2,r-1}} + C'\|u|_{Q_2}\|_{L^2}.
	\end{equation}
	If $u\in C^{k+2r+\delta, r+\delta/2}(Q_2)$, we have that\begin{equation}\label{eq_ffestimateholder}
		\|u|_{Q_1}\|_{C^{k+2r+\delta,r+\delta/2} } \leq C\|h|_{Q_2}\|_{C^{k+2r-2+\delta,r-1+\delta/2}} + C'\|u|_{Q_2}\|_{C^{0}}.
	\end{equation} 
	In each case we have constants $C = C(\kappa,\Lambda)$ and $C' = C'(\kappa,\Lambda,B_1,B_2,\ep_1,\ep_2)$. 
\end{theorem}
\begin{proof}[Proof sketch]
	The case $r=1$ for \eqref{eq_ffestimate} follows from the interior estimate in \cite{Kry08} Theorem 5.2.5, along with an estimate over an intermediate ball $Q_{R'}$ and Lemma \ref{lem_extensionoftrace}. The case $r>1$ follows by induction over an intermediate sequence of balls. For example, for $r=2$, we choose an intermediate ball $Q_1\subset Q_1'\subset Q_2$, so that we have estimates of the type \eqref{eq_ffestimate} for $Q_1\subset Q_1'$ with constants $C_1,C_1'$ and similarly for $Q_1'\subset Q_2$ with $C_2,C_2'$ then we can estimate one higher derivative $\pa_tu$ as\[
	\|\pa_tu|_{Q_1}\|_{L^2_{k+2,1}}\leq C_1\|(\pa_t-L)\pa_tu|_{Q_1'}\|_{L^2_{k,0}}+C_1'\|\pa_tu|_{Q_1'}\|_{L^2}\leq C_3\|\pa_t h|_{Q_2}\|_{L^2_{k,0}} +C'_3\|u|_{Q_2}\|_{L^2_{k+2,1}},
	\]
	where $C_3,C_3'$ depends on the earlier constants,
	thereby combining this with the estimate for $r=1$ proves the estimate \eqref{eq_ffestimate} for $r=2$. The H\"{o}lder estimates follow similarly from \cite[Theorem 8.12.1]{Kry96}.
\end{proof}

\begin{theorem}\label{thm_bfestimate}
	Let $u \in L^2_{k+2r,r}(P_2)$ and let us denote $h:= (\pa_t - L)u$. Then, we have that
	\begin{equation}\label{eq_bfestimate}
		\|u|_{P_1}\|_{L^2_{k+2r,r}} \leq C\|h|_{P_2}\|_{L^2_{k+2r-2,r-1}} + C'( \|u|_{P_2\cap\{t=0\}}\|_{L^2_{k+2r-1}} +\|u|_{P_2}\|_{L^2}).
	\end{equation}
If $u\in C^{k+2r+\delta, r+\delta/2}(P_2)$, we have that\begin{equation}\label{eq_bfestimateholder}
	\|u|_{P_1}\|_{C^{k+2r+\delta,r+\delta/2} } \leq C\|h|_{P_2}\|_{C^{k+2r-2+\delta,r-1+\delta/2}} + C'(\|u|_{P_2\cap\{t=0\}}\|_{C^{k+2r+\delta}}+\|u|_{P_2}\|_{C^{0}}).
\end{equation} 
In each case we have constants $C = C(\kappa,\Lambda)$ and $C' = C'(\kappa,\Lambda,B_1,B_2,\ep_1,\ep_2)$. 
	
\end{theorem}
\begin{proof}[Proof sketch]
	The case $r=1$ for the Sobolev estimates follows from an a-priori estimate from \cite{Kry08} Theorem 5.2.10, along with using a cutoff function to get an interior estimate. The case $r>1$ follows by induction over an intermediate sequence of balls as before. The H\"{o}lder estimates follow similarly from \cite[Exercise 9.2.5]{Kry96}.
\end{proof}

\subsubsection{Parabolic Sobolev Embedding}

We recall the following Sobolev embedding theorem for $L^2$-parabolic Sobolev spaces, which is Theorem 2.2.6 in \cite{Kry08}. 
\begin{theorem}\label{thm_sobolevembeddingparabolic}
	Define the quantity 
	\begin{equation}\label{eq_parabolicembeddingparameter}
		\kappa:= \frac{2s+1}{2r}+\frac{2(l+2s)+n}{2(k+2r)}.
	\end{equation}
	If $\kappa <1$, we have the embedding
	\begin{equation}\label{eq_parabolicholderembedding}
		L^2_{k+2r,r}(\bR^{n+1}_0)\hookrightarrow C^{l+2s,s}(\bR^{n+1}_0).
	\end{equation}
	
\end{theorem}

\subsubsection{Sobolev spaces on $\mancac$}\label{subsec_ordinarysobolev}
Let $X,Y\in \mancac$ and $f:X\rt Y$ be an m-fibration. Fix an m-metric on $X$, with corresponding m-volume form $dV_X$, and an m-connection $^\m\nabla$ on $^{\m}TX$. We can define standard Sobolev spaces via the following norms, which generalises the Euclidean parabolic Sobolev spaces:\[
\|u\|^2_{L^2}:=\int_X|u|^2dV_X,
\]
\[
\|u\|^2_{L^2_{k+2r,r}}:= \sum_{\substack{\beta\leq r \\
		\alpha\leq k+2r-2\beta}} \int_{X} |({^\m\nabla})^\beta({^\m\nabla}_{\textrm{rel}})^\alpha u|^2dV_X
\]
Note that here we are utilising the m-metric in order to dualise the inclusion \[
0\rt{^\m T}_\rel X \rt {^\m }TX
\]
to\[
0\rt {^\m}T^*_\rel X \rt {^\m}T^*X,
\]
which allows us to differentiate $({^\m}\nabla_\rel)^\alpha u$ with the connection $^\m\nabla$. Similarly we can define H\"{o}lder spaces as the completion of functions having continuous m-derivatives with norm\[
\|u\|_{C^{k+2r,r}(X)} := \sup_{x\in X,\; \beta\leq r, \; \alpha + 2\beta \leq k+2r}|({^{\m}}\nabla)^\beta({^\m\nabla}_{\textrm{rel}})^\alpha u|.
\]
Note that $\cap_{k,r\geq1} C^{k+2r,r} = \{\textrm{r-smooth functions on $X$}\}\supsetneq C^\infty(X)$, because a-smoothness requires a decay on the derivatives.
\subsubsection{Vector-valued operators}\label{subsec_vectorvalued}
Let $V$ be a finite-dimensional $\bR$-vector space, equipped with an inner product $h$. Define the trivial flat connection $\nabla^V$ on the trivial vector bundle $V\times \bR^{n+1}\downarrow\bR^{n+1}$. Consider Sobolev spaces as before,
\[
\|\sigma\|^2_{L^2_{k+2r,r}}:= \sum_{i\leq r}\sum_{j\leq k+2r-2i}\|\pa_t^i(\nabla^V)^j\sigma\|^2_{L^2},
\]
where the norm on vectors is with respect to the inner product. We consider a heat-type operator of the form\[
\pa_t-L: \Gamma^\infty(\uV)\rt\Gamma^\infty(\uV),
\]
where $L = a^{ij}\pa^2_{ij}+b^i\pa_i+c$ is a constant-coefficients strongly elliptic
 linear differential operator with coefficients $a^{ij},b^i,c\in \textrm{End}(V)$. Given this, we can prove the exact same interior estimates \eqref{eq_ffestimate}, \eqref{eq_ffestimateholder}, \eqref{eq_bfestimate}, \eqref{eq_bfestimateholder} for the vector-valued operator as before, see e.g.\ Chapter VII of the book by Ladyzhenskaya et. al. \cite{LSU68}.
\section{Motivating examples}\label{sec_applications}

We provide a sample of problems that could be approached by the theory to be developed, focusing for simplicity on the case of hypersurface mean curvature flow where the p.d.e.\ is scalar-valued. Another application to scalar-valued p.d.e.\ is Lagrangian mean curvature flow once the Hamiltonian isotopy classes are prescribed. These problems could also be considered for general codimension mean curvature flow and Ricci flow, by considering vector-valued operators (with a prerequisite gauge fixing via sections of the normal bundle for MCF and the DeTurck trick for Ricci flow). 
\subsection{Immersed mean curvature flow}
As the simplest possible application and a useful sanity check, consider the immersed mean curvature flow in Euclidean space $(X^{m},g) = (\bR^m,g_{\textrm{Eucl}})$ given by a family of smooth immersions $\iota_t: M^{m-1} \hookrightarrow X$ for time $t\in[0,T]$ of a compact manifold $M^{m-1}$. Suppose that at time $T$, there is a point $p\in X$ through which $l\geq1$ sheets of the flow pass. We can think of this as a Type I `singularity', with shrinker being a disjoint union of planes $\sqcup_l\bR^{m-1}\subset \bR^{m}$. We may ask the following question: suppose we perturb the initial data $\iota_0(M)$ slightly, via the graph of a section of the normal bundle. Can we characterise the perturbations such that $l$ sheets of the flow hit $p$ at time $T$? Among these, can we further characterise perturbations that preserve the tangent plane at time $T$? Among the previous, can we characterise perturbations that preserve the second-order Taylor expansion at time $T$? And so on.

To answer this question, we begin by considering the flow as an immersed submanifold in spacetime $\iota=\bigsqcup_{t\in[0,T]} \iota_t: M\times[0,T]\rt X\times [0,T]$. Then, we can write the MCF equation as a p.d.e.\ on sections of the normal bundle,\[
(\pa_t-\Delta-|A|^2)u+Q[u]=0, \;\; u|_{t=0}=u_0.
\]
Given that $\iota_T^{-1}(p)=\{p_1,\dots,p_l\}$, we can write the condition for the flow to pass through $p$ with $l$ sheets as $\{u_t\}$ satisfying $u_T(p_1)=u_T(p_2)=\cdots=u_T(p_l)=0$ as sections of the normal bundle. Now, we can set up a nonlinear mapping,
\begin{equation}\label{eq_pointmapping}
u_0\mapsto (u_T(p_1),\cdots,u_T(p_l)).
\end{equation}
Then, the answer to our first question is given by the elements in the kernel of this mapping. Similarly, studying the tangent plane at the point amounts to studying the mapping 
\begin{equation}\label{eq_higherpointmapping}
u_0\mapsto (u_T(p_1),\dots,u_T(p_l), \nabla u_T(p_1)\dots,\nabla u_T(p_l)),
\end{equation}
and so on for higher-order Taylor series expansions. To study what happens near the point $(p,T)$ in space-time, we perform a \textit{parabolic blow-up} of the space-time (for details, see e.g.\ \S 3 of \cite{LMPS23}, \S 4.2 of \cite{Gho25}). This construction is the key tool allowing us to provide a precise description of the geometry at the `singular' point.

\begin{figure}[!htb]
	\centering
	\def\svgwidth{0.4\columnwidth}
	%% Creator: Inkscape 1.4.3 (0d15f75042, 2025-12-25), www.inkscape.org
%% Accompanies image file 'approxsol_Layer 1_Type I.pdf' (pdf, eps, ps)
\begingroup%
  \makeatletter%
  \providecommand\color[2][]{%
    \errmessage{(Inkscape) Color is used for the text in Inkscape, but the package 'color.sty' is not loaded}%
    \renewcommand\color[2][]{}%
  }%
  \providecommand\transparent[1]{%
    \errmessage{(Inkscape) Transparency is used (non-zero) for the text in Inkscape, but the package 'transparent.sty' is not loaded}%
    \renewcommand\transparent[1]{}%
  }%
  \providecommand\rotatebox[2]{#2}%
  \newcommand*\fsize{\dimexpr\f@size pt\relax}%
  \newcommand*\lineheight[1]{\fontsize{\fsize}{#1\fsize}\selectfont}%
  \ifx\svgwidth\undefined%
    \setlength{\unitlength}{257.7775802bp}%
    \ifx\svgscale\undefined%
      \relax%
    \else%
      \setlength{\unitlength}{\unitlength * \real{\svgscale}}%
    \fi%
  \else%
    \setlength{\unitlength}{\svgwidth}%
  \fi%
  \global\let\svgwidth\undefined%
  \global\let\svgscale\undefined%
  \makeatother%
  \begin{picture}(1,0.59632825)%
    \lineheight{1}%
    \setlength\tabcolsep{0pt}%
    \put(0,0){\includegraphics[width=\unitlength,page=1]{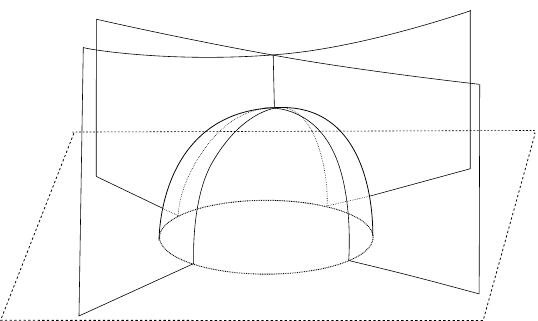}}%
    \put(0.90592775,0.5824601){\color[rgb]{0,0,0}\makebox(0,0)[lt]{\lineheight{1.25}\smash{\begin{tabular}[t]{l}$t=0$\end{tabular}}}}%
    \put(0.95435942,0.14847672){\color[rgb]{0,0,0}\makebox(0,0)[lt]{\lineheight{1.25}\smash{\begin{tabular}[t]{l}$t=T$\end{tabular}}}}%
  \end{picture}%
\endgroup%

	\caption{The compactification in the parabolic blow-up. This is a manifold with corners and a-corners with one a-boundary and two ordinary boundaries.}
	\label{fig_immblowup}
\end{figure}
Then we can consider the linearisation of \eqref{eq_pointmapping}, which will have well-behaved Fredholm theory, and its kernel will provide a finite-codimensional submanifold of initial conditions answering the original question. Similarly controlling the rate of decay at the front face for higher rates provides a Fredholm mapping away from a set of critical rates, answering the follow-up questions.

\subsection{Asymptotically conical Type I singularities}
Consider the mean curvature flow of an immersed submanifold developing a Type I singularity at an isolated point $(p,T)$, modelled on an asymptotically conical shrinker. We can study the geometry of the mean curvature flow near the singular point by the same method as before: we construct the parabolic blow-up at the singular point, and the evolving submanifold will compactify to a locally embedded manifold with corners and a-corners in the blown-up space-time. Having a well-behaved Fredholm theory for the heat operator allows us to perform precise analysis of singularity formation, and moreover determine the codimension of perturbations of the initial condition in which the singularity forms. For example, we can consider a weighted mapping between weighted Sobolev spaces with weight function $\rho^\lambda$ for $\rho = \sqrt{|x-p|^2+|T-t|}$ and $\lambda\in \bR$,
\[
(\pa_t-\Delta,|_{t=0}):\rho^\lambda L^2_{k+2r,r}\rt \rho^{\lambda-2}L^2_{k+2r-2,r-1}\oplus L^2_{k+2r-1},
\]
and the dimension of the cokernel of this operator at a rate $\lambda =2+\ep$ tells us the codimension of initial data in which the singularity forms at the same point $(p,T)$ with the same tangent flow. This would also allow us to construct singularities modelled on a given asymptotically conical shrinker by using a `parabolic gluing' method. There have been examples of such constructions in the literature, for example Stolarski \cite{Sto26} for Ricci flow and Lee--Zhao \cite{LZ24} for codimension $1$ mean curvature flow.
\subsection{Flowing out of conical singularities}
Consider the `inverse' problem to the previous one: given an initial conically singular submanifold, can we construct a mean curvature flow out of it by gluing in expanders into the tangent cone? This problem was treated using the language of manifolds with corners in the setting of network flow by Lira--Mazzeo--Pluda--Sa\'{e}z \cite{LMPS23} and for Lagrangian mean curvature flow by the author in \cite{Gho25}. Again, the main idea is to construct the parabolic blow-up in which the glued family of approximate solutions compactifies and study the heat operator on this family. We mention that this problem has also been studied in the literature for Ricci flow, e.g.\ Gianniotis--Schulze \cite{GS18}.

\subsection{Asymptotically conical Type II singularities}
Now, we come to a more difficult problem of modelling finite-time Type II singularities. The main difference from the Type I case is that one expects multiple scales at which the analysis is important, corresponding to $O(\sqrt{|T-t|})$ (Type I) and $o(\sqrt{|T-t|})$. We explain a heuristic picture of how a class of such singularities modelled on a minimal asymptotically conical submanifold may be expected to form, following the outline in the paper by Joyce on Lagrangian mean curvature flow \cite[Principle 3.9]{Joy15}: to first order, the flow quickly becomes approximately the tangent cone of the model at scale $O(\sqrt{|T-t|})$, or equivalently, if we do parabolic rescaling, we should observe a minimal submanifold at scale $o(1)$ with a fixed asymptotic cone given by the stationary tangent flow. As the rescaled time $\tau= -\log(|T-t|) \uparrow \infty$, this minimal submanifold slowly wanders in the moduli space of all minimal submanifolds with the given cone, until it shrinks to zero size at $\tau=\infty$. We treat this problem for Lawlor neck pinching in Lagrangian MCF in a follow-up paper \cite{Gho26b}.

\section{Heat operators over a mixed fibration}\label{sec_heatoperators}

Let $X^{m+1}$ be a compact, $m+1$-dimensional manifold with corners and a-corners, such that the connected components of the boundary are faces, and each boundary $\Sigma$ of $X$ admits a collar neighbourhood that is diffeomorphic to $\Sigma\times\llbracket0,\ep)$ if it is an a-boundary, or $\Sigma\times[0,\ep)$ if it is an ordinary boundary (i.e.\ the $b$-normal bundle $^bN\pa X$ is \textit{untwisted} in the sense of \cite[\S 4.4]{Joy16}). We will consider Sobolev spaces weighted by powers of boundary-defining functions at the a-boundaries, so henceforth we fix a numbering $\Sigma_1,\dots,\Sigma_N$ of the a-boundaries and corresponding boundary-defining functions $\{\rho_1,\dots,\rho_N\}$. 
\begin{figure}[!tb]
	\centering
	\def\svgwidth{0.35\columnwidth}
	%% Creator: Inkscape 1.4.3 (0d15f75042, 2025-12-25), www.inkscape.org
%% Accompanies image file 'example.pdf' (pdf, eps, ps)
\begingroup%
  \makeatletter%
  \providecommand\color[2][]{%
    \errmessage{(Inkscape) Color is used for the text in Inkscape, but the package 'color.sty' is not loaded}%
    \renewcommand\color[2][]{}%
  }%
  \providecommand\transparent[1]{%
    \errmessage{(Inkscape) Transparency is used (non-zero) for the text in Inkscape, but the package 'transparent.sty' is not loaded}%
    \renewcommand\transparent[1]{}%
  }%
  \providecommand\rotatebox[2]{#2}%
  \newcommand*\fsize{\dimexpr\f@size pt\relax}%
  \newcommand*\lineheight[1]{\fontsize{\fsize}{#1\fsize}\selectfont}%
  \ifx\svgwidth\undefined%
    \setlength{\unitlength}{205.27114988bp}%
    \ifx\svgscale\undefined%
      \relax%
    \else%
      \setlength{\unitlength}{\unitlength * \real{\svgscale}}%
    \fi%
  \else%
    \setlength{\unitlength}{\svgwidth}%
  \fi%
  \global\let\svgwidth\undefined%
  \global\let\svgscale\undefined%
  \makeatother%
  \begin{picture}(1,1.46563601)%
    \lineheight{1}%
    \setlength\tabcolsep{0pt}%
    \put(0,0){\includegraphics[width=\unitlength,page=1]{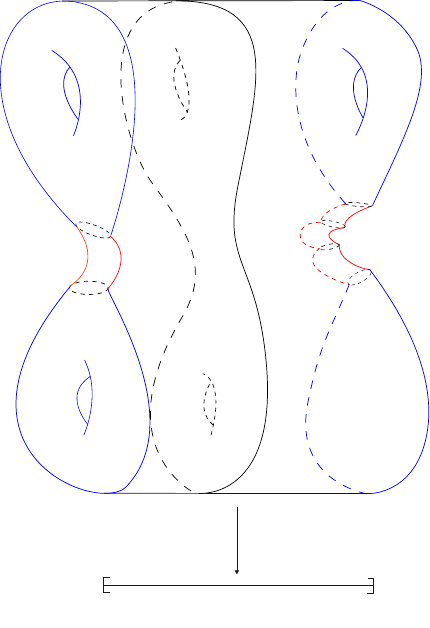}}%
    \put(0.20181283,0.00409591){\color[rgb]{0.10196078,0.10196078,0.10196078}\makebox(0,0)[lt]{\lineheight{1.25}\smash{\begin{tabular}[t]{l}$a$\end{tabular}}}}%
    \put(0.80725199,0.00437321){\color[rgb]{0.10196078,0.10196078,0.10196078}\makebox(0,0)[lt]{\lineheight{1.25}\smash{\begin{tabular}[t]{l}$b$\end{tabular}}}}%
    \put(0,0){\includegraphics[width=\unitlength,page=2]{example.pdf}}%
  \end{picture}%
\endgroup%

	\caption{An example of the manifold-with corners $X$ with an m-fibration to $[a,b]$. Blue denotes ordinary boundaries, and red denotes a-boundaries.}
	\label{fig_example}
\end{figure}
Here, $X$ represents the compactified space-time manifold corresponding to a geometric flow with $1$ dimension of time and $m$ dimensions of space (that may have developed a singularity at final time, or begins from a singularity at initial time, or flows with singularities). Thus, we have a natural global time function, $t:X\rt I$, where $I$ is a compact one-dimensional interval. We assume that $t: X\rt I$ is a proper m-fibration as in Definition \ref{defn_mfibration}.

For simplicity, we assume that the time interval $I$ is $[a,b]$ with $-\infty <a<b <+\infty$ (the other natural possibilities could be $\llbracket-\infty,a]$ or $[a,+\infty\rrbracket$ or $\llbracket-\infty,+\infty\rrbracket$). Over the interior $I^\circ=(a,b)$, the time function gives a fibration, $t^{-1}(I^\circ) \cong M\times I^\circ$ for $M$ being a compact manifold with corners and a-corners. For simplicity, we assume that $M$ is a manifold without boundary (we would let $M$ have ordinary boundaries to prescribe time-dependent boundary conditions, while we would let $M$ have a-boundaries to prescribe asymptotics near infinity or singularities). We finally assume for simplicity that $X$ has only codimension $2$ corners. Combined with the condition that $t: X\rt I$ is an m-fibration, this implies there are only five kinds of local models $V\subset X$ for this fibration:
\begin{defn}\label{defn_localmodels} 
 \begin{enumerate}
	\item (the interior) For $B_R\subset \bR^{m}$, $(x,t')\in V=B_R\times(-\delta,\delta)$, with a local time variable $t'=t-c\in(-\delta,\delta)$ for some $c\in(a,b)$.
	\item (at an ordinary boundary)For $B_R\subset \bR^{m}$, $(x,t')\in V= B_R\times[0,\delta)$, with a local time variable $t'=t-a \textrm{ (or }b-t)\in[0,\delta)$.
	\item (at an a-boundary)For $B_R\subset \bR^{m}$, $(x,t')\in V= B_R\times\llbracket0,\delta)$, we have a local time variable $t'=t-a \textrm{ (or }b-t)\in [0,\delta)$.
	\item (at the intersection of two a-boundaries) For $B_R\subset \bR^{m-1}$, $(\sigma,r,s)\in  V=B_R\times \llbracket0,\delta)\times\llbracket0,\delta)$, where $r=0, s=0$ are a-boundaries, we have a local time variable $t'=t-a \textrm{ (or } b-t)=rs \in [0,\delta)$.
	\item (at the intersection of an a-boundary and an ordinary boundary) For $B_R\subset \bR^{m-1}$, $(\sigma,r,s) \in V= B_R\times \llbracket0,\delta)\times [0,\delta)$, where $r=0$ is an a-boundary, $s=0$ is an ordinary boundary, we have a local time variable $t'= t-a \textrm{ (or } b-t)=rs \in [0,\delta)$
\end{enumerate}
\end{defn}
\begin{figure}[!htb]
	\centering
	\def\svgwidth{0.8\columnwidth}
	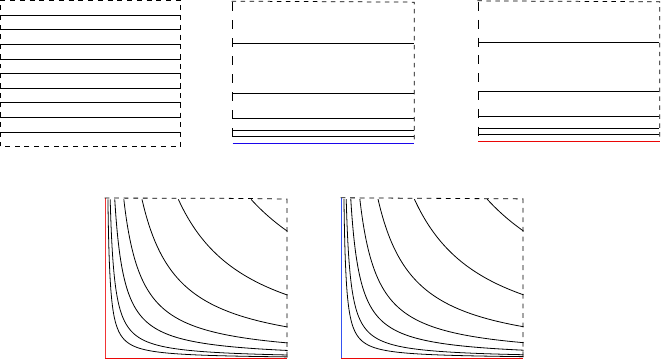
	\caption{The five kinds of local models along with their constant-time slices. Blue denotes ordinary boundaries, and red denotes a-boundaries.}
	\label{fig_localmodels}
\end{figure}

We will consider a heat-type operator on $X$, of the form \begin{equation}
	\cT=\pa_v-\cL,
\end{equation}	
 where $\cL$ is an elliptic second order relative differential operator (in the sense of \S \ref{subsubsec_reloper}) with respect to the time fibration and $\pa_v\in \Gamma^{\infty}(^\m T  X)$ is a \textit{time-like vector field} (i.e.\ $\pa_v(t)>0$ away from the a-boundaries of $  X$), which may degenerate at some of the a-boundaries, in a prescribed manner which we now  describe. Let $\be=(\eta_1,\dots,\eta_N) \in \bR^N_{\geq0}$ be a collection of weights corresponding to the a-boundary components $\{\Sigma_1,\dots,\Sigma_N\}$ of $X$. Denote the weight function $\brh^{\be}:= \rho_1^{\eta_1}\dots\rho_N^{\eta_N}$. Then we can locally write near each a-boundary $\Sigma_j$, \begin{equation}\label{eq_localformv}
 \pa_v = F\brh^{\be}\cdot \cV_t +\pa_\sigma,
 \end{equation}
  for $F$ being a non-vanishing function, $\pa_\sigma$ being a smooth relative vector field and $\cV_t$ being a local choice of vector field such that $\la\cV_t\ra\oplus {^\m T_\rel X}\cong {^\m TX}$ (recall that since $t$ is an m-fibration, $^mT_\rel X$ is well-defined by Proposition \ref{prop_mfibre}). For example, in each of the local models in Definition \ref{defn_localmodels} we can take $\cV_t = \pa_{t'}$ for (1),(2), $\cV_t = t'\pa_{t'}$ for (3), $\cV_t= r\pa_r$ for (4), $\cV_t= \pa_s$ for (5). Alternatively, we could fix an m-metric to make a global choice of $\cV_t$ as the unit-length perpendicular to $^\m T_\rel X\subset {^\m TX}$. We make the following definition of decay which is easier to check:
  \begin{defn}\label{defn_vectorv}
	A time-like vector field $\pa_v$ over the m-fibration $t:X\rt[a,b]$ \textit{degenerates to order $\be$} if $\pa_v(t):X\rt [0,\infty)$ is an a-smooth mapping as in Definition \ref{defn_asmoothmap}, and in the interior of an a-boundary $\Sigma_j$, the function $\pa_v(t)/t'$ decays to leading order $\rho_j^{\eta_j}$ for $t'$ being the local time variable at $\Sigma_j$.
\end{defn}
\begin{lem}
	The above definition is equivalent to the local form \eqref{eq_localformv}.
\end{lem}
\begin{proof}
	Near each a-boundary, we can write the local time function $t'$ on a collar neighbourhood as an a-smooth mapping from $\Sigma_j\times\llbracket 0,\ep)\rt[0,\infty)$ since $t$ is a $b$-fibration. Now for a fixed choice of vector field $\cV_t$ with $\la \cV_t\ra \oplus {^\m T_\rel X}\cong {^\m TX}$, it follows that the mapping $\cV_t(t'):\Sigma\times\llbracket0,\ep)\rt [0,\infty)(\textrm{or }(-\infty,0])$ is a-smooth with order of vanishing at $\Sigma_j$ the same as the order of vanishing of $t'$ as can be checked by computing in the local models of Definition \ref{defn_localmodels}. So, the vector field being of the form \eqref{eq_localformv} implies Definition \ref{defn_vectorv}.
	
	Conversely, we can uniquely write the vector field $\pa_v = f\cdot \cV_t+\pa_\sigma$ for $f: \Sigma_j\times\llbracket0,\infty)\rt \bR$ being a-smooth. Then again, $\cV_t(t')$ is an a-smooth mapping from $\Sigma\times\llbracket0,\ep) \rt [0,\infty)$. So $\pa_v(t)=\pm f\cV_t(t'): \Sigma\times\llbracket0,\ep) \rt [0,\infty)$ being a-smooth implies that it is uniquely a product of boundary-defining functions times a positive function, which implies that $\pm f : \Sigma_j\times \llbracket0,\ep) \rt [0,\infty)$ is itself a-smooth with decay at the a-boundaries given by $\rho_j^{\eta_j}$. Finally, $f$ is non-vanishing at the ordinary boundaries because $\pa_v(t)>0$ at the interior of ordinary boundaries. 
\end{proof}

We further assume that the decay rate $\eta_j$ is $0$ at a-boundaries that intersect an ordinary boundary, i.e.\ $\eta_j= 0$ if $\Sigma_j\cap \pa^{\ord}  X \neq \emptyset$ (this assumption is necessary in the analysis at a vanishing a-boundary). 
\begin{ass}\label{ass_cond}
We collect all our assumptions in the following list:
\begin{enumerate}
    \item\label{item_time} The function $t:X\rt I = [a,b]$ with $-\infty<a<b<\infty$ is a proper m-fibration.
    \item\label{item_xgeo} The fibration over the interior $I^\circ$ is of the form $M\times(a,b)$ with $M$ a manifold without boundary.
    \item\label{item_codim} The space-time $X$ has at most codimension $2$ corners.
    \item\label{item_operator} The operator is of the form \[
    \cT=\pa_v-\cL,\; \pa_v=\brh^{\be}\cV_{t}
    \]
    for a tuple of weights $\be=\{\eta_1,\dots,\eta_N\}$ corresponding to the a-boundaries $\Sigma_1,\dots,\Sigma_N$.
    \item\label{item_etanonvanish} The decay rate is $0$ at a-boundaries that intersect an ordinary boundary, i.e.\ $\eta_j= 0$ if $\Sigma_j\cap \pa^{\ord}  X \neq \emptyset$.
\end{enumerate}

\end{ass}

\begin{rem}
	If the weight $\be \equiv 0$, then the analysis simplifies greatly. In particular, we obtain optimal interior estimates at all rates $\bla$, and the normal operators at a-boundaries are enough to control Fredholmness (instead of the degenerate normal operators as defined in \S \ref{subsec_vanabound}). We note that in the cases where the singularity behaviour is Type I, i.e.\ either the singularity formation modelled on a shrinker or flowing out of a singularity modelled on an expander, the weight $\be$ does vanish, as in the cases of \cite{LMPS23,Gho25}. It is for applications to modelling Type II singularities that we have introduced a non-zero weight $\be$. 
\end{rem}

Henceforth, we fix an m-metric $g\in \Sym^2(^\m T^*  X)$ on $  X$, which gives an associated m-volume form $dV_{  X}$ and Levi-Civita connection $^\m\nabla$ on $^\m T  X$. This also induces a relative metric $g_\rel\in \Sym^2(^\m T^{*}_\rel  X)$, and corresponding relative Levi-Civita connection $^\m\nabla_\rel$ on $^\m T_\rel   X$. We note a lemma regarding the existence of nice cutoff functions on $  X$:
\begin{lem}\label{lem_partition}
	Given a positive bound $1\gg\epsilon>0$, numbers of derivatives $k\geq1$, and open neighbourhoods $\{U_1,\dots,U_l\}$ of some a-boundaries $\{\Sigma_1,\dots,\Sigma_l\}$, there exists an open covering of $  X$ given by $\{U_0,U_1,\dots,U_l\}$, where $U_0$ is an open set not intersecting $\{\Sigma_1,\dots,\Sigma_l\}$, and a partition of unity $\{\chi_0,\chi_1,\dots,\chi_l\}$ subordinate to this cover, such that the derivatives satisfy $\|^\m\nabla^j\chi_i\|_{C^{0}}<\epsilon$ for $i=0,1,\dots,l$ and $j=1,\dots,k$.
\end{lem}
\begin{proof}
	We will prove this by an inductive construction. If we have a single a-boundary $\Sigma_1$, it is easy to construct such a pair of functions $\chi_0,\chi_1$ as follows: take a collar neighbourhood $U=\Sigma_1\times \llbracket0,\delta)\ni(\sigma,\rho)$, and define functions on this neighbourhood, \[
	\chi_1:= 1-\chi\left(\frac{\log\rho}{R_1}+R_2\right), \; \chi_0 := \chi\left(\frac{\log\rho}{R_1}+R_2\right),
	\] 
	where $\chi:\bR\rt [0,1]$ is a fixed cutoff function satisfying $\chi(x)=0$ for $x<0$ and $\chi(x)=1$ for $x>1$. Taking $R_1\gg1$, we can ensure the m-derivatives of $\chi_0$ are small, and taking $R_2\gg0$, we can ensure that the support of $\chi_1$ is within $U_1$, and then extend these to $\{0,1\}$ outside $U$ respectively.
	    
	Now suppose we have more than one a-boundary, then we choose a $\epsilon' \ll \epsilon$ and do the construction for the a-boundaries $\Sigma_1,\dots,\Sigma_{l-1}$ with constant $\epsilon'$. As a result we have an open covering $\{U_0,U_1,\dots,U_{l-1}\}$ with cutoff functions $\chi_0,\chi_1,\dots,\chi_{l-1}$. Then, we do the construction for the single boundary $\Sigma_l$ with constant $\epsilon'$, and obtain a covering $\{U_0',U_1'\}$ with cutoff functions $\chi_0',\chi_1'$. Then $\chi_0'\chi_0,\chi_0'\chi_1,\dots,\chi_0'\chi_{l-1},\chi_1'$ is a partition of unity subordinate to the cover $\{U_0\cap U_0',U_1,\dots,U_{l-1},U_1'\}$, and moreover $
	\|^\m\nabla \chi_0'\chi_j\|_{C^0}\leq \|^\m\nabla \chi_0'\|_{C^0}\cdot\|\chi_j\|_{C^0}+\|\chi_0'\|_{C^0}\cdot\|^\m\nabla \chi_j\|_{C^0}\leq 2\epsilon'$,
	and similarly for higher derivatives. Therefore, we choose $\epsilon'$ sufficiently small and obtain the claim.
\end{proof}

\subsection{Local operators}
We can write down a covering of $  X$ by the local models in Definition \ref{defn_localmodels}, and note how the operator and corresponding Sobolev spaces look like in these local models. We need to work with pairs of local models to prove interior estimates:
\begin{defn}\label{defn_modelpairs}
	Define a \textit{model pair} to be a pair $(U,V)$ of open neighbourhoods $U:=B_R(\ze)\times(-\delta,\delta)\ssubset V:= B_R(\ze)\times (-2\delta,2\delta)\subset \bR^n$ with $B_R(\ze)\ssubset B_{2R}(\ze)\ssubset\bR^{n-1}$, or $U:=\llbracket0,\delta)^k\times[0,\delta)^l\times B_R(\ze) \ssubset V:= \llbracket0,2\delta)^k\times[0,2\delta)^l\times B_{2R}(\ze) \subset \bR^{k,n}_l$ with $B_R(\ze)\ssubset B_{2R}(\ze)\subset \bR^{n-k-l}$. We have a time fibration $t(x_1,\dots,x_n)=x_n$ for the former, and $t(x_1,\dots,x_n)=x_1\cdots x_{k+l}$ for the latter.
\end{defn}
 Now, we note the pull-back of the heat operator on the local models in Definition \ref{defn_localmodels}.
\begin{defn}\label{defn_modeloperators}
	 In each of the following, $\cL = a^{ij}\pa^2_{ij}+b^i\pa_i+c $ denotes a uniformly elliptic operator in the vector fields $\{\pa_{x^1},\dots,\pa_{x^m}\}$ spanning the relative tangent bundle, $F$ is a non-vanishing a-smooth function on the model space, $t'$ is the local time variable, and $\eta\in \bR_{\geq0}$ denotes a decay parameter. \begin{enumerate}
		\item\label{item_interior} (the interior) On $(x,t')\in V=B_R\times(-\delta,\delta)$, we have an operator $F\pa_{t'}-\cL$.
		\item\label{item_ordbound} (at an ordinary boundary) On $(x,t')\in V=B_R\times[0,\delta)$, we have an operator $F\pa_{t'}-\cL$.
		\item\label{item_abound} (at an a-boundary) On $(x,t')\in V=B_R\times\llbracket0,\delta)$, we have an operator $F{t'^\eta}\cdot t'\pa_{t'}-\cL$.
		\item\label{item_aboundabound} (at the intersection of two a-boundaries) On $(\sigma,r,s)\in V=B_R\times \llbracket0,\delta)\times\llbracket0,\delta)$, where $\{r=0\}, \{s=0\}$ are a-boundaries, we have a local time variable $t'=rs$, the relative tangent bundle is spanned by sections $\la \pa_{x^1}:=\pa_{\sigma^1},\dots,\pa_{x^{m-1}}:=\pa_{\sigma^{m-1}},\pa_{x^m}:=r\pa_r-s\pa_s\ra$, a non-degenerate time-like vector field is $s\pa_s+r\pa_r$, and we have operator $F r^{\eta_r} s^{\eta_s}\cdot(r\pa_r+s\pa_s)-\cL$.
		\item\label{item_aboundordbound} (at the intersection of an a-boundary and an ordinary boundary) On $(\sigma,r,s)\in V=B_R\times \llbracket0,\delta)\times [0,\delta)$, where $\{r=0\}$ is an a-boundary, $\{s=0\}$ is an ordinary boundary, we have a local time variable $t'=rs$, the relative tangent bundle is spanned by sections $\la \pa_{x^1},\dots,\pa_{x^{m-1}},\pa_{x^m}:=r\pa_r-s\pa_s\ra$, a non-degenerate time-like vector field is $\pa_s+sr\pa_r$, and we have operator $F\cdot (\pa_s+sr\pa_r)-\cL$. Note that the vector field is non-degenerate in this model by Assumption \ref{item_etanonvanish}.
	\end{enumerate}
\end{defn}

On each model, we can write down a flat m-metric $g_0$, i.e.\ $g_0=\sum_{1\leq i,j \leq n}g_{ij}de^ide^j$ for $de^i=x_i^{-1}dx_i$ for $1\leq i\leq k$ or $de^i=dx^i$ for $k+1\leq i \leq n$, and $[g_{ij}]$ being a positive-definite matrix. Associated to this we have a relative metric $g_\rel$ on the time slices, the Levi-Civita connection $^\m\nabla$ on $^\m TV$ and relative Levi-Civita connection $^\m\nabla_\rel$ on $^\m T_\rel V$. On each model space, we define parabolic Sobolev spaces $L^2_{\be;k+2r,r}(V)$ as the completion of $C^\infty_0(V)$-functions under the following norm, by taking up to $k+2r$ relative derivatives and $r$ derivatives in the time-like vector field: \begin{equation}
	\|u\|^2_{L^2_{\be;k+2r,r}(V)}:=\sum_{\substack{\beta\leq r \\
			\alpha+2\beta\leq k+2r}} \int_V |({\pa_v})^\beta{^\m\nabla}_{\textrm{rel}}^\alpha u|^2 d\vol.
\end{equation}
We can write it more explicitly in terms of derivatives of the local relative vector fields $\{\pa_{x^1},\dots,\pa_{x^m}\}$ and time-like vector field $\pa_v$ as \begin{equation}\label{eq_explicitnorm}
\|u\|^2_{L^2_{\be;k+2r,r}(V)}\equiv \sum_{\substack{\beta\leq r \\
		(\sum_{i=1}^m\alpha_i)+2\beta\leq k+2r}} \int_V |({\pa_v})^\beta\pa_{x^1}^{\alpha_1}\dots\pa_{x^m}^{\alpha_m} u|^2 d\vol.
\end{equation}

Note that for $r=0$, we have $L^2_{\be;k,0}\equiv L^2_{k,0}$ are the standard Sobolev spaces as defined in \S \ref{subsec_ordinarysobolev}, but for $r\geq 1$, the $L^2_{\be}$-norm is weaker than the standard norm if the vector field $\pa_v$ degenerates at some a-boundaries, i.e.\ possibly in cases \eqref{item_abound}, \eqref{item_aboundabound} of Definition \ref{defn_modeloperators} (for an example see Remark \ref{rem_lowreg}). Now, we can define a mapping on each model space, \begin{equation}\label{eq_sobolevmappingmodel}
	(\cT:=\pa_v-\cL, |_{\pa^{\ord}  V^-}) : L^2_{\be;k+2r,r}( V) \rt L^2_{\be;k+2r-2,r-1}(  V) \oplus L^2_{k+2r-1}(\pa^\ord  V^-),
\end{equation}
where $\pa^\ord V^-$ denotes the ordinary boundaries of V where the time-like vector field is inward pointing. That the above mapping is well-defined follows from the local expression of the operator in Definition \ref{defn_modeloperators} and the explicit norm in \eqref{eq_explicitnorm}, noting from calculation in local models that \[
\pa_{x^j}(\pa_vu) = \pa_v(\pa_{x^j}u)+O(\pa_vu,\pa_{x^j}u),
\]
due to decay properties of $\pa_v$.

We will define the \textit{H\"{o}lder distance} between points $x,y\in   V$ for $d_  V(x,y)$ sufficiently small as follows:
\begin{equation}\label{eq_holderdistance}
d(x,y) := \sqrt{d_{  V}(x,y)^2+\int_0^l\frac{|\la\gamma'(s),\pa_v^\perp(\gamma(s))\ra|}{|\pa_v^\perp(\gamma(s))|^2}ds },
\end{equation}
where $d_{V}$ is the distance function for the m-metric, $\gamma: [0,l]\rt   V$ is the unit-speed geodesic in the m-metric joining the points $x,y$, and $\pa_v^\perp(\gamma(s))$ denotes the projection of $\pa_v$ onto the normal to the relative tangent bundle at the point $\gamma(s)$. We will see in \S \ref{subsubsec_goodcover} that this localises correctly to the usual notion of H\"{o}lder distances as defined in \S\ref{subsec_holdereuclidean}. Then, we can define H\"{o}lder spaces with norm\[
\|u\|_{C^{\be;k+2r+\delta,r+\delta/2}(V)} := \sup_{x\in V,\;  \beta \leq r,\; \alpha + 2\beta \leq k+2r}|\pa_v^\beta({^\m\nabla}_{\textrm{rel}}^\alpha u)| + \sum_{\beta\leq r,\; \alpha+2\beta=k+2r}[\pa_v^\beta({^\m\nabla}_{\textrm{rel}}^\alpha u)]_{\delta/2;\delta},
\]
where $[\cdot]_{\delta/2;\delta}$ is the norm given by \[
[v]_{\delta/2,\delta}:= \sup_{d_V(x,y)<\textrm{inj}_V, \; x \neq y}\frac{|v(x)-v(y)|}{d(x,y)^\delta}.
\]
As before \eqref{eq_sobolevmappingmodel}, we have a well-defined continuous mapping on H\"{o}lder spaces,
\begin{equation}\label{eq_holdermappingmodel}
	(\cT, |_{\pa^{\ord}  V^-}) : C^{\be;k+2r+\delta,r+\delta/2}(V) \rt C^{\be;k+2r-2+\delta,r-1+\delta/2}(V) \oplus C^{k+2r+\delta}(\pa^\ord  V^-).
\end{equation}

 Now, we claim the following interior estimates on each pair $(U,V)$, analogous to the interior estimate on Euclidean spaces in \S \ref{subsec_interior}:
\begin{prop}\label{prop_localapriori}
	As in Definition \ref{defn_modeloperators}, suppose we have a parabolic operator of the form $\pa_v-a^{ij}\pa^2_{ij}-b^i\pa_i-c$ with the coefficients $a^{ij},b^i,c \in C^{\be;k+2r-2+\delta,r-1+\delta/2}$ and $a^{ij}$ being uniformly elliptic. Given a function $u\in L^2_{{k+2r,r,\loc}}(V)$, we have that for $h:= (\pa_v-\cL)u$, \begin{equation}\label{eq_modelestimatesobolev}
	\|u|_U\|_{L^2_{\be;k+2r,r}}\leq C(\|h|_V\|_{L^2_{\be;k+2r-2,r-1}}+\|u|_V\|_{L^2}+\|u|_{ \pa^\ord V^-}\|_{L^2_{k+2r-1}}).
	\end{equation}
	If instead $u\in C^{k+2r+\delta, r+\delta/2}_{\loc}(V)$, we have that\begin{equation}\label{eq_modelestimateholder}
		\|u|_{U}\|_{C^{\be;k+2r+\delta,r+\delta/2} } \leq C(\|h|_{V}\|_{C^{\be;k+2r-2+\delta,r-1+\delta/2}}+\|u|_{\pa^\ord V^-}\|_{C^{k+2r+\delta}} + \|u|_{V}\|_{C^{0}}).
	\end{equation} 
\end{prop}
We will also prove the following Sobolev embedding theorem. We continue to denote by $\brh^{\be}$ the function given by $t^\eta$ on local model \eqref{item_abound} or $r^{\eta_r}s^{\eta_s}$ on local model \eqref{item_aboundabound} (the reason for this notation is that the function $\brh^{\be}$ on $X$ localises to a bounded multiple of these functions on each local model).  

\begin{prop}\label{prop_localembedding}
	Given $l,s\geq1$, for all $k,r\geq1$ such that for $\xi<1$ (as defined in \eqref{eq_parabolicembeddingparameter}), we have a continuous embedding of Sobolev spaces, $L^2_{\be;k+2r,r}(V) \hookrightarrow \brh^{-\frac{\be}{2}}C^{\be;l+2s+\delta,s+\delta/2}(  V)$. If $\chi:   V \rt \bR$ is a smooth function vanishing at the a-boundaries, then the mapping $\chi\cdot:L^2_{\be;k+2r,r} \rt L^2_{\be;k+2r-2,r-1}$ given by $u\mapsto \chi u$ is compact.
\end{prop}
To prove these, we will further subdivide the model pair $(U,V)$ into a \textit{good cover}.

\subsubsection{Good covers}\label{subsubsec_goodcover}
\begin{defn}
Fix positive constants $R,\ep>0$ and numbers of derivatives $k,r\geq1 $. A countable collection $\{(S_i,T_i,\Psi_i)\}_{i\in I}$ of open sets $(S_i,T_i)$ of the form $S\subset T$ with $S=B_{R}\times [0,\ep)$, $T=B_{2R}\times [0,2\ep)$ or $S=B_{R}\times (-\ep,0]$, $T=B_{2R}\times (-2\ep,0]$ or $S=B_R\times(-\ep,\ep), T = B_{2R}\times(-2\ep,2\ep)$, along with a-diffeomorphisms to $V$ of the form $\Psi_i: T_i\rt \Psi_i(T_i)$, is defined to be a \textit{good cover}, if the following are satisfied:
\begin{enumerate}
	\item\label{item_timeslices} The time-slices are preserved.
	\item\label{item_cover} The images $\Psi_i(S_i)$ cover $U\backslash \pa^{\ma}U$.
	\item\label{item_fincover} Each point $x\in V$ lies in the image of at most $N<\infty$ mappings $\Psi_i(T_i)$.
	\item\label{item_normequiv} We have an equivalence of norms,\begin{equation}\label{eq_psiequiv}
	\|u|_{\Psi_i(T_i)}\|_{C^{\be;k+2r+\delta,r+\delta/2}}\equiv \|\Psi_i^*(u)\|_{C^{k+2r+\delta,r+\delta/2}(T_i)},
	\end{equation}
	with the equivalence with respect to constants uniform in $i\in I$.
	\item\label{item_operatorpull} The pull-back of $\pa_v-\cL$ on each $T_i$ is a parabolic operator of the form $F\pa_t-a^{ij}\pa^2_{ij}-b^i\pa_i-c$ with coefficients having uniform bounded $C^{k+2r-2+\delta,r-1+\delta/2}(T_i)$-norms.
	\item The pull-back of the volume forms is uniformly equivalent up to a constant weight, i.e. there is a $\kappa>1$ independent of $i\in I$ such that 
	\begin{equation}\label{eq_volequiv}
\kappa^{-1}\|\Psi^*u\|^2_{L^2(T_i)}\leq C_i\cdot \|u|_{\Psi_i(T_i)}\|^2_{L^2}\leq	\kappa\|\Psi^*u\|^2_{L^2(T_i)},
	\end{equation}
	for a constant $C_i = \brh^{-\be}(\Psi_i(0,0))\in(0,\infty)$ which may depend on $i\in I$. 
\end{enumerate}
\end{defn}

Now we prove the existence of a good cover on each local model:
	\begin{enumerate}
		\item On the local models \eqref{item_interior}, \eqref{item_ordbound}, it is clear that we can cover with finitely many mappings satisfying conditions \eqref{item_timeslices} and \eqref{item_cover}, thus they give a good cover.
		\item On the local model \eqref{item_abound}, we have two cases: \begin{itemize}
			\item Suppose $\eta=0$. Then, we can first write down a good covering of $B_R\times(\delta/2,\delta)$ with finitely many charts $N_0$. Now, we observe that the metric and the norms are invariant under the action of $(0,1)\subset \bR_+$ given by $(x,t)\mapsto (x,\lambda t)$. Composing the charts with this action, we get a countable covering of $V\backslash \pa^{\ma}V$, and due to invariance of the metric and the operator, we obtain the good cover. Note that the constant $N=2N_0$ for this covering. 
			\item Suppose $\eta>0$. Then the same cover as before satisfies the requirements \eqref{item_timeslices}, \eqref{item_cover}, \eqref{item_fincover} for the good cover. However, the pull-back norms on the charts are not equivalent, because of the decaying vector field. The pull back operator can be written as $F \cdot t_i^\eta \exp(t\eta)\pa_t-\cL$, where $t_i\downarrow 0$. We modify this cover by `stretching' the time variable of each $T_i$, so that the pull-back operator becomes $F\cdot \pa_t-\cL$. It is clear that the norms of $F$ are bounded in these new charts. Then, we can divide each chart into about $O(t_i^{-\eta})$ uniformly sized balls and composing this sub-division with our initial mappings gives us the good cover on $V$. Since each point is covered by $\leq2N_0$ images of the initial cover, each point is covered by at most $N=4N_0$ images of the refinement.
		\end{itemize}
		\item On the local model \eqref{item_aboundabound}, we again split into two cases:
		\begin{itemize}
			\item Suppose $(\eta_r,\eta_s)=(0,0)$. Then, we first write down a good cover of $B_R\times(\delta/2,\delta)\times(\delta/2,\delta)$ by finitely many charts $N_0$. Now, we have an isometric action of $(0,1)^2\subset\bR_+\times \bR_+$ on the local model, by $(\sigma,r,s)\mapsto (\sigma,\lambda r, \mu s)$ preserving the metric and the pull-back norms. Composing the charts with this action, we get a covering of $V\backslash \pa^{\ma}V$ which is a good cover, with constant $N=4N_0$.
			\item Suppose $(\eta_r,\eta_s) \neq (0,0)$. As before, we first write down a covering as above, satisfying requirements \eqref{item_timeslices}, \eqref{item_cover}, \eqref{item_fincover}. The pull-back operator can be written as $F\cdot r_i^{\eta_r}s_i^{\eta_s}\exp(t\eta'+x\mu')\pa_t-\cL$. Again, we can `stretch' in the time variable of each $T_i$ by a factor of $r_i^{-\eta_r}s_i^{-\eta_s}$ and obtain a sub-division giving a good cover, with each point covered by at most $N=8N_0$ images.
		\end{itemize}
		\item On the local model \eqref{item_aboundordbound}, we first choose a good covering of $B_R\times (\delta/2,\delta)\times[0,\delta)$ by finitely many charts $N_0$. Then we again have an isometric action by $(0,1)\subset\bR_+$ on the $\llbracket0,\delta)$-factor preserving the metric and the norms such that translates of the charts under this action cover $V\backslash \pa^{\ma}V$, and similarly as above this gives a good cover with each point covered by at most $N=2N_0$ images.
	\end{enumerate}
 This proves the existence of a good cover. 

\begin{proof}[Proof of Proposition \eqref{prop_localapriori} and Proposition \eqref{prop_localembedding}:]Let $u\in L^2_{\be;k+2r,r,\loc}$. Then, \[
	\|u|_U\|^2_{L^2_{\be;k+2r,r}}\leq \sum_{i\in I}\|u|_{\Psi_i(S_i)}\|^2_{L^2_{\be;k+2r,r}}\leq \sum_{i\in I}C'C_i\kappa(\|u|_{S_i}\|^2_{L^2_{\be;k+2r,r}})\]
	\[\leq \sum_{i\in I}C'C_i\kappa(\|h|_{T_i}\|^2_{L^2_{\be;k+2r-2,r-1}}+\|u|_{\pa^\ord{T_i}^-}\|^2_{L^2_{k+2r-1}}+\|u|_{T_i}\|^2_{L^2})\]
	\[	\leq \sum_{i\in I}C'^2\kappa^2(\|\cT u|_{\Psi_i(T_i)}\|^2_{L^2_{\be;k+2r-2,r-1}}+\|u|_{\pa^\ord\Psi_i(T_i)^-}\|^2_{L^2_{k+2r-1}}+\|u|_{\Psi_i(T_i)}\|^2_{L^2})\]
    \[\leq NC'^2\kappa^2(\|h|_{V}\|^2_{L^2_{\be;k+2r-2,r-1}}+\|u|_{\pa^{\ord}V^-}\|^2_{L^2_{k+2r-1}}+\|u|_V\|^2_{L^2}),
	\]
	for $C'$ being the constant in \eqref{eq_psiequiv}. Here the first estimate is the triangle inequality, the second estimate follows from \eqref{eq_psiequiv} and \eqref{eq_volequiv} on $L^2_\eta$-norms, the third estimate follows from \eqref{eq_ffestimate} and \eqref{eq_bfestimate}, and the last estimate follows again from \eqref{eq_psiequiv} and \eqref{eq_psiequiv}. This proves the estimate \eqref{eq_modelestimatesobolev}, and \eqref{eq_modelestimateholder} follows similarly. The first Sobolev embedding in Proposition \ref{prop_localembedding} follows from standard Euclidean Sobolev embedding as in Theorem \eqref{thm_sobolevembeddingparabolic} and noting $C_i \approx\rho^{\be}$. The second embedding follows from an Arzela-Ascoli type argument along with standard Sobolev embedding.
\end{proof}

\subsection{A-priori estimates}
Let recall the decomposition of the a-boundary of $  X$ as the disjoint union of connected components $\pa^{\ma}  X=\Sigma_1\sqcup\dots\sqcup \Sigma_N$, each having a globally defined boundary-defining function $\rho_j$ and further a collar neighbourhood diffeomorphic to $\Sigma_j\times\llbracket0,\delta)\ni(\sigma,\rho_j)$, with $\rho_j$ smoothly interpolated to be $1$ outside a slightly bigger collar neighbourhood. Similar to the definition on model spaces, we define new parabolic Sobolev spaces $L^2_{\be;k+2r,r}(X)$ on $X$ as the completion of $C^\infty_0(X)$-functions under the following norm, by taking up to $k+2r$ relative derivatives and $r$ derivatives in the time-like vector field: 
\begin{equation}\label{eq_sobolevnormeta}
	\|u\|^2_{L^2_{\be;k+2r,r}(X)}:=\sum_{\substack{\beta\leq r \\
			\alpha\leq k+2r-2\beta}} \int_{X} |({\pa_v})^\beta{^\m\nabla}_{\textrm{rel}}^\alpha u|^2 dV_{  X}\equiv\sum_{\substack{\beta\leq r \\
			\alpha\leq k+2r-2\beta}} \int_{X} |{\brh^{\beta\cdot \be}}({^\m\nabla})^\beta({^\m\nabla}_{\textrm{rel}})^\alpha u|^2dV_{  X},
\end{equation}
where the weight $\brh^{\beta\cdot \be}:=\rh^{\beta\cdot\eta_1}\rg^{\beta\cdot\eta_2}\dots\rho_N^{\beta\cdot\eta_N}$. The above equivalence of norms follows because we have an equivalence of norms $|\pa_v u|+|^\m\nabla_\rel u| \equiv |\brh^{\be} \cdot {^\m}\nabla u|+|^\m\nabla_\rel u|$ and so on for higher derivatives, due to the local form of $\pa_v$ as in \eqref{eq_localformv}.
\begin{rem}
	We are taking derivatives in $\pa_v$ in the following way: we have $^\m \nabla_\rel^\alpha u \in \Gamma^\infty(\otimes^{\alpha} {^\m} T^*_\rel)$, then we include $^\m T_\rel ^*X\hookrightarrow {^\m}T^*X$ via dualising $^\m T_\rel X\hookrightarrow {^\m} TX$ with respect to an m-metric, and finally we take $^\m \nabla_{\pa_v}(^\m \nabla^\alpha_\rel u)$. 
\end{rem}

\subsubsection{Weighted Sobolev and H\"{o}lder spaces}\label{subsubsec_sobolevspaces}
 For a tuple of real-valued weights $(\lah,\lag,\dots,\lambda_N) \in \bR^N$, we consider Sobolev spaces weighted by the function $\brh^{\bla}:=\rh^{\lah}\rg^{\lag}\dots\rho_N^{\lambda_N}$:\[
\|u\|^2_{\brh^{\bla}L^2_{\be;k+2r,r}(X)}:= \sum_{\substack{\beta\leq r \\
		\alpha\leq k+2r-2\beta}} \int_{\mathbb{L}} |\rh^{-\lah}\rg^{-\lag}\dots\rho_N^{-\lambda_N}(\pa_v)^\beta{^\m\nabla}_{\textrm{rel}}^\alpha u|^2 dV_{  X}\equiv \|\brh^{-\bla} u\|^2_{L^2_{\be;k+2r,r}(X)}.
\]

The above equivalence of norms follows because the m-derivatives of the weights $\brh^{-\bla}$ are bounded by their corresponding powers, $|^\m\nabla^\beta(\brh^{-\bla})| = O(\brh^{-\bla})$. Taking a finite cover of $X$ by model pairs as in Definition \ref{defn_modeloperators} and a partition of unity subordinate to it, we can easily transfer all the preceding results for the local models to $X$. So, we get a well-defined bounded linear map,
	\begin{equation}\label{eq_sobolevmappingL}
		(\cT, |_{ \pa^{\ord}  X^-}) : \brh^{\bla}L^2_{\be;k+2r,r}(  X) \rt \brh^{\bla}L^2_{\be;k+2r-2,r-1}(  X) \oplus \brh^{\bla}L^2_{k+2r-1}(\pa^\ord  X^-),
	\end{equation}
	where $\pa^\ord  X^-$ denotes the ordinary boundary of $  X$ where the vector field $\pa_v$ points inwards, or equivalently the ordinary boundary at initial time $t=a$.
 We define H\"{o}lder spaces on $X$, given by $C^{\be;k+2r+\delta,r+\delta/2}(X)$, with norm\begin{equation}\label{eq_holdernormeta}
\|u\|_{C^{\be;k+2r+\delta,r+\delta/2}(X)} := \sup_{x\in X,\;\beta\leq r,\; \alpha + 2\beta \leq k+2r}|\brh^{\beta\cdot\be}{^{\m}}\nabla^\beta{^\m\nabla}_{\textrm{rel}}^\alpha u| + \sum_{\beta\leq r,\;\alpha+2\beta=k+2r}[\brh^{\beta\cdot\be}{^{\m}\nabla}^\beta{^\m\nabla}_{\textrm{rel}}^\alpha u]_{\delta/2;\delta},
\end{equation}
where $[\cdot]_{\delta/2;\delta}$ is the norm given by \[
[v]_{\delta/2,\delta}:= \sup_{d_  X(x,y)<\textrm{inj}_  X, \; x \neq y}\frac{|v(x)-v(y)|}{d(x,y)^\delta},
\]
with the distance defined exactly as \eqref{eq_holderdistance} with $V$ replaced by $X$. As before, we have a well-defined continuous mapping on H\"{o}lder spaces,
	\begin{equation}\label{eq_holdermappingL}
	(\cT, |_{\pa^{\ord}  X^-}) : \brh^{\bla}C^{\be;k+2r+\delta,r+\delta/2}(X) \rt \brh^{\bla}C^{\be;k+2r-2+\delta,r-1+\delta/2}(  X) \oplus \brh^{\bla}C^{k+2r+\delta}(\pa^\ord  X^-),
\end{equation}

\begin{rem}\label{rem_lowreg}
	These $\be$-dependent Sobolev and H\"{o}lder spaces are of far weaker regularity than being r-smooth, i.e.\ $\bigcap_{k,r\geq1}C^{\be;k+2r,r}\supsetneq \{\textrm{r-smooth functions}\}$. For example, consider on $\llbracket0,1/2]$, the function $\log(x)^{-1}e^{i\log(x)^3}$. Then, we have for weight $\be=\{1\}$ associated to the a-boundary $\{0\}$, the function lies in all $\be$-spaces. However, we have that $|(x\pa_x)u|=\Theta(|\log(x)|)$, so it is not in any ordinary Sobolev space with $r\geq1$ derivatives.
\end{rem} 
Note that working on weighted spaces is equivalent to working with a conjugated operator $\brh^{-\bla}\cT(\brh^{\bla}\cdot)$, and we can write this as $\brh^{-\bla}\cT(\brh^{\bla}\cdot) = \cT + \textrm{ lower order relative terms}$. As a result, the conjugated operator is also a heat-type operator as discussed, and we immediately obtain weighted a-priori estimates and further obtain Sobolev embeddings on $X$.
\begin{theorem}\label{thm_degenerateapriori}
	Let us assume $\cT$ is a heat operator with coefficients in $C^{\be;k+2r+\delta,r+\delta/2}$. We have the following estimate for $u\in L^2_{\be;2,1,\loc}$: 
	\begin{equation}\label{eq_apriori}
		\|\pa_vu\|_{L^2} +\|^\m\nabla_\rel^2u\|_{L^2}+\|^\m\nabla_\rel u\|_{L^2}\leq C(\|\cT u\|_{L^2}+\|u|_{\pa^{\ord}  X^-}\|_{L^2_{1}}+\|u\|_{L^2}).
	\end{equation}
	We have weighted estimates for all higher orders for $u\in L^2_{\be;k+2r,r,\loc}$,
	\begin{equation}\label{eq_apriorihigher}
		\begin{split}
			\|u\|_{\brh^{\bla}L^2_{\be;k+2r,r}} \leq C(\|\cT u\|_{\brh^{\bla}L^2_{\be;k+2r-2,r-1}}+\|u|_{\pa^{\ord}  X^-}\|_{\brh^{\bla}L^2_{k+2r-1}}+\|u\|_{\brh^{\bla}L^2}).
		\end{split}
	\end{equation}
	We also have Schauder estimates for $u\in C^{\be;k+2r+\delta,r+\delta/2,\loc}$, \begin{equation}\label{eq_schauder}
		\|u\|_{\brh^{\bla}C^{\be;k+2r+\delta,r+\delta/2} } \leq C(\|\cT u\|_{\brh^{\bla}C^{\be;k+2r-2+\delta,r-1+\delta/2}}+\|u|_{\pa^{\ord}  X^-}\|_{\brh^{\bla}C^{k+2r+\delta}} + \|u\|_{\brh^{\bla}C^{0}}).
	\end{equation} 
\end{theorem}
\begin{proof}[Proof sketch]
	Since $X$ is compact, we can cover it with finitely many model pairs as in Definition \ref{defn_modelpairs}. Then, we can add up the interior estimates \eqref{eq_modelestimatesobolev} and \eqref{eq_modelestimateholder} of Proposition \ref{prop_localapriori} over this cover to obtain the claimed estimates.
\end{proof}

\begin{theorem}\label{thm_sobolevembedding}
	Given $l,s\geq1$, for all $k,r\geq1$ such that for $\xi<1$ (as defined in \eqref{eq_parabolicembeddingparameter}), we have a continuous embedding of Sobolev spaces, $\brh^{\bla}L^2_{\be;k+2r,r}(  X) \hookrightarrow \brh^{\bla-\frac{\be}{2}}C^{\be;l+2s,s}(X)$. If $\chi:   X \rt \bR$ is a smooth function vanishing at the a-boundaries, then the mapping $\chi\cdot:\brh^{\bla}L^2_{\be;k+2r,r} \rt \brh^{\bla}L^2_{\be;k+2r-2,r-1}$ given by $u\mapsto \chi u$ is compact.
\end{theorem}
\begin{proof}[Proof sketch]
	As before, we can cover $X$ with finitely many model pairs as in Definition \ref{defn_modelpairs}, then apply Proposition \ref{prop_localembedding} to obtain the claim.
\end{proof}

\subsection{Global mapping properties}
We recall the following functional analytic result, see e.g.\ Buhler--Salamon \cite[Lemma 4.3.9]{BH18}: 
\begin{lem}\label{lem_funcana}
	Let $X,Y,Z$ be Banach spaces and $\cT: X\rightarrow Y$ be a bounded, linear map and $\cK: X \rightarrow Z$ be a bounded, linear, compact map. Suppose we have the estimate \[
	\|u\|_X \leq C(\|\cT u\|_Y + \|\cK u\|_Z)
	\]
	for a constant $C>0$. Then, we have that $\cT$ has finite dimensional kernel and closed image, i.e.\ it is semi-Fredholm.
\end{lem}

\begin{theorem}\label{thm_fredholmlocal}
	The mapping \eqref{eq_sobolevmappingL} between weighted Sobolev spaces is semi-Fredholm if each a-boundary $\Sigma_j$ has an open neighbourhood $U_j$ such that for a function $u\in L^2_{2,1}$ supported in $U_j$ compactly away from the a-boundaries, we have the estimate
	\begin{equation}\label{eq_L^2}
		\|u\|_{\brh^{\bla}L^2}\leq C_j(\|\cT u\|_{\brh^{\bla}L^2}+\|u|_{\pa^{\ord}  X^-}\|_{\brh^{\bla}L^2_{1}}),
	\end{equation}
	for a positive constant $C_j>0$. Conversely, if the mapping is semi-Fredholm for $k=0,r=1$, then each a-boundary has an open neighbourhood on which this estimate holds.
\end{theorem}

\begin{proof}
	Let us first assume that each a-boundary $\Sigma_j$ has an open neighbourhood $U_j$ such that the estimate \eqref{eq_L^2} holds for $u$ compactly supported in $U_j$. By Lemma \ref{lem_partition} applied to the set of all a-boundaries, we can obtain an open set $U_0$ and a partition of unity $\{\chi_0,\dots,\chi_N\}$ subordinate to the covering $\{U_0,U_1,\dots, U_N\}$ with constant $\epsilon$. Then, for a general $u\in C^\infty_0$, we can write, $u = u_0+\sum_{j=1}^Nu_j$, for $u_j=\chi_j\cdot u$. Now, we have that \[
	\|u\|_{\brh^{\bla}L^2}\leq \sum_{j=0}^N\|u_j\|_{\brh^{\bla}L^2}\leq \|\chi_0u\|_{\brh^{\bla}L^2}+\sum_{j=1}^NC_j(\|\cT u_j\|_{\brh^{\bla}L^2}+\|u_j|_{\pa^{\ord}X^-}\|_{\brh^{\bla}L^2_1}).
	\] 
	Then, we can estimate by the triangle inequality and noting $\|^m\nabla \chi_j\|_{C^{2,1}}<\epsilon$,\[
	\|\cT u_j\|_{\brh^{\bla}L^2}\leq \|\cT u\|_{\brh^{\bla}L^2}+K\epsilon\|u\|_{\brh^{\bla}L^2}+K\epsilon\|^\m\nabla_\rel u\|_{\brh^{\bla}L^2}. 
	\]
	Therefore, we get \[
	\|u\|_{\brh^{\bla}L^2}\leq \|\chi_0u\|_{\brh^{\bla}L^2}+C(\|\cT u\|_{\brh^{\bla}L^2}+\|u|_{\pa^\ord X^-}\|_{\brh^{\bla}L^2_1})+C'K\epsilon\|u\|_{\brh^{\bla}L^2_{1,0}}.
	\]
	Combining this with the a-priori estimate \eqref{eq_apriorihigher}, we get \[
	\|u\|_{\brh^{\bla}L^2_{\be;k+2r,r}}\leq C'(\|\cT u\|_{\brh^{\bla}L^2_{\be;k+2r-2,r-1}}+\|u|_{\pa^{\ord}  X^-}\|_{\brh^{\bla}L^2_{k+2r-1}})+\|\chi_0u\|_{\brh^{\bla}L^2}+C''K\epsilon\|u\|_{\brh^{\bla}L^2_{1,0}}.
	\]
	Therefore, by choosing $\epsilon$ sufficiently small, we can get rid of the last term on the right, to get \[
	\|u\|_{\brh^{\bla}L^2_{\be;k+2r,r}}\leq 2C'(\|\cT u\|_{\brh^{\bla}L^2_{\be;k+2r-2,r-1}}+\|u|_{\pa^{\ord}  X^-}\|_{\brh^{\bla}L^2_{k+2r-1}})+2\|\chi_0u\|_{\brh^{\bla}L^2}.
	\]
	Finally, by Lemma \ref{lem_funcana} and the compactness of $u\mapsto \chi_0u$ by Theorem \ref{thm_sobolevembedding}, we obtain that \eqref{eq_sobolevmappingL} is semi-Fredholm.
	
	Conversely, if $\Sigma_j$ is an a-boundary such that for any open neighbourhood $U_j$ the estimate \eqref{eq_L^2} does not hold, we have a sequence of compactly supported smooth functions $\{u_i\}_{i\geq 1}$ such that $\|u_i\|_{\brh^{\bla}L^2}=1$ and $\|\cT u_i\|_{\brh^{\bla}L^2}+\|u_i|_{\pa^\ord X^-}\|_{\brh^{\bla}L^2_1} \downarrow 0$. Moreover, we can choose $u_i$ such that their supports do not intersect, by choosing shrinking open neighbourhoods. Therefore, the $\{u_i\}$ cannot subsequentially converge to an element, which implies that \eqref{eq_sobolevmappingL} has an infinite-dimensional kernel or non-closed image, contradicting semi-Fredholmness.
\end{proof}

Given the above, we make the following definition:

\begin{defn}\label{defn_injective}
	Let $M\in M^{\ma}_k(X)$ be a codimension $k\geq1$ a-corner as in Definition \ref{defn_acorners}. The operator $\cT$ between weighted Sobolev spaces at rate $\bla$, or equivalently the conjugated operator $\cT^{\bla}$, is said to be \textit{injective near $M$} if there exists an open neighbourhood $U$ of $M$ such that for $u$ being a-smooth and compactly supported in $U$ away from the a-boundaries, we have an estimate    
	\begin{equation}
		\|u\|_{L^2}\leq C(\|\cT^{\bla} u\|_{L^2}+\|u|_{\pa^{\ord}  X^-}\|_{L^2_{1}}),
	\end{equation}
	for a positive constant $C>0$.
\end{defn}

\subsubsection{Normal operator}

Now, we will define the important notion of a \textit{normal operator} at a codimension $k$ a-corner, following e.g. Loya--Melrose \cite[\S A.3.]{LM03}.
\begin{defn}\label{defn_normaloperator}
	Let $M$ be a codimension $k\geq 1$ a-corner component of $  X$. Let the boundary-defining functions corresponding to $M$ be $\rho_{i_1},\dots, \rho_{i_k}$, so that $M\subseteq \Sigma_{i_1}\cap\cdots\cap \Sigma_{i_k}$. Define the \textit{family of normal operators} associated to $\cT$ as follows: for a parameter $\xi \in \bR^k$, given a function $u\in C^{\infty}(M)$, define
	\begin{equation}\label{eq_normaloperatordefn}
	N_M(\cT)(\xi)(u):= \brh^{-i\xi}\cT(\brh^{i\xi}\tu)|_{M},
	\end{equation}
	for $\tu$ being a smooth extension of $u$ on a small neighbourhood $U$ of $M$. As a result, we obtain a family of operators dependent on the parameter $\xi\in \bR^k$.
\end{defn}
It was shown by Loya--Melrose \cite{LM03} that a uniformly $b$-elliptic operator of order $l$ on a manifold with corners as a mapping between weighted $b$-Sobolev spaces $\brh^{\bla}{L}^p_{k+l}\rt \brh^{\bla}L^p_{k}$ is Fredholm if and only if the normal operators corresponding to the conjugated operator are isomorphisms at each boundary over all parameters. In our setting, there are added complications due to the possible degeneracy of the time-like vector field, so the following definition will be useful:
\begin{defn}
	A codimension $k$ a-corner $M\in M^{\ma}_k(X)$ is called \textit{vanishing} if $\brh^{\be}|_M = 0$, and called \textit{non-vanishing} otherwise. In particular, the a-boundaries where $\eta_j>0$ are \textit{vanishing a-boundaries}, and the boundaries where $\eta_j=0$ are \textit{non-vanishing a-boundaries}.
\end{defn}
\begin{rem}
	In general, the family of normal operators fail to control the behaviour at a vanishing a-corner $M$, since\[
	N_M(\cT)(\xi)\equiv -N_M(\cL)(\xi),
	\]
	so the normal operator only sees the relatively elliptic part rather than the whole operator. To get around this, we will define analogous notions of \textit{degenerate normal operators} such that the mapping properties near the corner is controlled by this new family of operators. 
\end{rem}

\subsection{Model problems}\label{subsubsec_modelproblems}
We outline the method to prove estimates of the form \eqref{eq_L^2} near each a-boundary (and also codimension $k$ a-corners). First, we reduce to showing estimates for a p.d.e.\ on a \textbf{model space}. Then we can do a reduction of dimension in the time variable and obtain elliptic p.d.e.s. (which turn out to be the normal operators).
\subsubsection{Fourier transform}
We recall the definition of Fourier transform we will be using here. Note that we work with complex-valued $L^2$-spaces throughout. 

\begin{defn}
	Let $M$ be a compact manifold with corners and a-corners, and let $f : M\times \bR^k \rt \bR$ be a function in $C^\infty_0(M\times \bR^k)$. The \textit{Fourier transform} of $f$ in the $\bR^k$-directions is the function $\hat{f} \in L^2(M\times \bR^k)$, defined as:
	\[
	\hat{f}(p,\xi):= \int_{\bR^k}e^{-i\xi\cdot x}f(p,x)dx.
	\] 
\end{defn}

A basic property of the Fourier transform is:
\begin{equation}\label{eq_fourierder}
\cF(\pa_{x^j} f)(p,\xi) = i\xi^j\hat{f}(p,\xi).
\end{equation}
One of the key results about the Fourier transform is that:
\begin{theorem}[Plancherel's theorem]
	The Fourier transform is an isometry on $L^2(M\times \bR^k)$ up to the factor $(2\pi)^{k/2}$, i.e.\[
	\|\hat{f}\|_{L^2}=(2\pi)^{k/2}\|f\|_{L^2}.
	\]
\end{theorem}

\subsubsection{Model problems at non-vanishing a-corners}
Let us begin with defining the model space for a non-vanishing a-corner $\Xi \in M_k^{\ma}(X)$. Define the model space to be $\Xi \times \llbracket0,\infty\rrbracket^k$ with a diffeomorphism of $\Xi\times \llbracket 0,\delta)^k$ onto a small neighbourhood $U$ of $\Xi$, given by $\Psi: \Xi\times \llbracket0,\delta)^k \rt U$ mapping $(\sigma,\rho_1,\dots, \rho_k) \mapsto (\sigma,\rho_1,\dots, \rho_k)$. As $\Xi$ is non-vanishing, we can write on $U$,
\begin{equation}\label{eq_decompositionmodel}
	\cT = \cT_\Xi+\textrm{vanishing operators},
\end{equation}
for $\cT_\Xi$ defined on $\Xi\times\llbracket0,\delta)^k$ by formally `freezing coefficients' at $\Xi$. Therefore, by taking $\delta$ sufficiently small, the second term is dominated by a small factor times $C^{\be;2,1}$-derivatives. Thus for $u$ supported in $U$, we can transfer the problem to $\Xi\times \llbracket0,\infty\rrbracket^k$ and reduce to proving estimates for the translation-invariant operator $\cT_\Xi$, which gives us the \textbf{model problem}. We also recall the following identity relating the Fourier transform of $\cT_\Xi$ in the variables $(\log(\rho_1),\dots,\log(\rho_k))$ to the normal operator in Definition \ref{defn_normaloperator}: \[
\cF(\cT_{\Xi,\xi} u)(\xi) = N_\Xi(\cT_\Xi)(\xi)\hu(\xi),
\]
which follows from \eqref{eq_fourierder}. This gives a family of second-order elliptic operators on $\Xi$.

\subsubsection{Model problems at vanishing a-corners}
On the other hand, given a vanishing a-corner $\Xi\in M^{\ma}_k(X)$, we cannot proceed directly as above. Instead, we can write,\[
\cT = \brh^{\be_\Xi}F(\sigma)\cdot\cV_t-\cL_\Xi+\textrm{vanishing operators},
\]
on a small open neighbourhood $U$ which replaces \eqref{eq_decompositionmodel}, and we can work with $\brh^{\be_\Xi}F(\sigma)\cdot\cV_t-\cL_\Xi$ for $\cL_\Xi$ being the operator obtained by freezing the coefficients of the relatively elliptic part $\cL$ at $\Xi$. However, this is not translation-invariant. Therefore, we will perform a change of coordinate over the interior (by Assumption \ref{item_etanonvanish}, $\Xi$ being vanishing implies all boundaries of $U$ are a-boundaries, and since it suffices to prove estimates working with $C^\infty_0$-functions, such a change-of-coordinate is valid). Now, we split into the two possible cases (since we assumed $X$ has at most codimension $2$-corners):

\begin{enumerate}
	\item For $k=1$ (vanishing a-boundary), we do a change of coordinate for the interior of a collar neighbourhood: $(\sigma,\rho_1)\in \Sigma\times\llbracket0,\delta) \mapsto (\sigma,t') \in \Sigma\times\llbracket0,\delta')$, where $t'=\rho_1f(\sigma)$ is the local time variable. By Assumption \ref{item_etanonvanish}, all the neighbouring boundaries of $\Sigma$ are a-boundaries, so we can take $\cV_t = \rho_1\pa_{\rho_1}$. As a result, the new operator becomes \[
	F(\sigma)f(\sigma)^{-\eta_1}\cdot t'^{\eta_1}t'\pa_{t'}-\cL_\Xi,
	\]
	where $\cL_\Xi$ transforms to a relatively elliptic operator with respect to the new time-fibration, i.e.\ $\cL_\Xi= a^{ij}\pa^2_{ij}+b^i\pa_i+c$ for $\la\pa_i\ra_{i=1}^m = \la \pa_{\sigma_1},\dots,\pa_{\sigma_{m}}\ra$. Now, we further perform a change of variable to $\ttil=\exp(-\eta^{-1}_1t'^{-\eta_1})$ so that $\ttil\pa_{\ttil}=t'^{\eta_1}t'\pa_{t'}$, and thus we obtain the following \textbf{model problem}: \[
	\cT_{\Xi,\deg} = \mu \ttt-\cL_\Xi,
	\]
	for $\mu=F(\sigma) f(\sigma)^{-\eta_1}$ on the {model space} $(\sigma,t')\in\Xi\times \llbracket0,\infty\rrbracket$, which is translation-invariant in $\log(\ttil)$.
	\begin{figure}[!htb]
		\centering
		\def\svgwidth{0.35\columnwidth}
		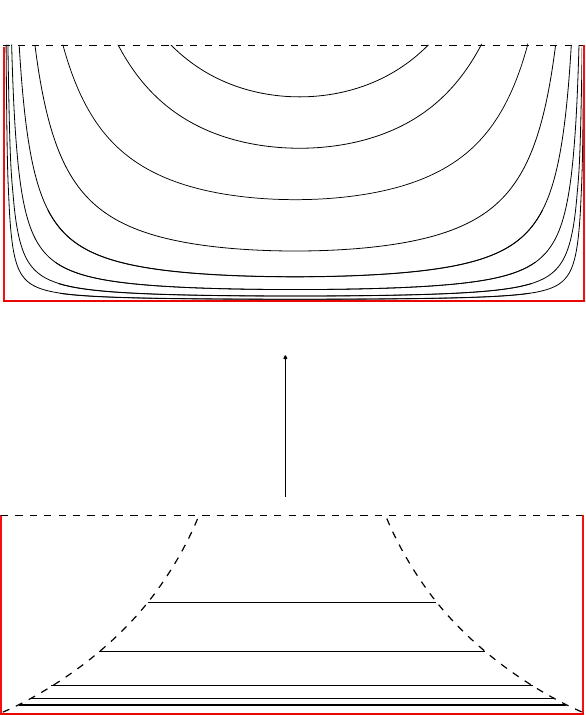
		\caption{Reparametrisation near a vanishing a-boundary.}
		\label{fig_reparamabound}
	\end{figure}

	\item For $k=2$ (vanishing a-corner), we do a change of coordinate for the interior: $(\sigma,\rho_1,\rho_2)\mapsto (\sigma, s=\rho_1\rho_2^{-1},t'=\rho_1\rho_2f(\sigma))$. This transforms the operator to \[
	s^{(\eta_2-\eta_1)/2}F(\sigma)f(\sigma)^{(\eta_1+\eta_2)/2}t'^{(\eta_1+\eta_2)/2}\ttp-\cL_\Xi
	\] and thus we obtain the following \textbf{model problem}: \[
	\cT_{\Xi,\deg} = \mu \ttt-\cL_\Xi,
	\]
	for $\ttil=\exp(-((\eta_1+\eta_2)/2)^{-1}t'^{(\eta_1+\eta_2)/2})$, and $\mu=s^{(\eta_2-\eta_1)/2}F(\sigma)f(\sigma)^{(\eta_1+\eta_2)/2}$ on the model space $(\sigma,s,\ttil)\in\Xi\times \llbracket0,\infty\rrbracket\times \llbracket0,\infty\rrbracket$, which is translation-invariant in $\log(\ttil)$.
	\begin{figure}[htb]
		\centering
		\def\svgwidth{0.37\columnwidth}
		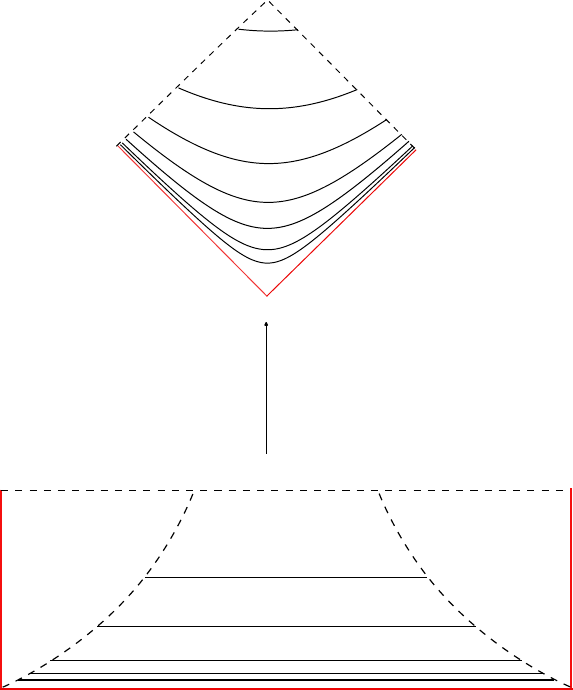
		\caption{Reparametrisation near a vanishing codimension $2$ a-corner.}
		\label{fig_reparamacorner}
	\end{figure}
\end{enumerate}

\subsubsection{Positive operators}

In practice, it may be easy to obtain injectivity near an a-corner if the operator is \textbf{positive}, as in the following definition:

\begin{defn}\label{defn_positive}
	The family of normal operators for $\cT$ at $\Xi$ is said to be \textit{positive}, if we have that with respect to choice of an m-volume form on $\Xi$, \[
	\|u|_{\pa^{\ord}\Xi^-}\|^2_{L^2(\Xi)}+\Re(\la N_\Xi(\cT)(\xi)u,u\ra) \geq C\|u\|_{L^2(\Xi)}^2,
	\]
	for a constant $C$ independent of $\xi\in \bR^k$, where $\pa^{\ord}\Xi^-$ is the subset of $\pa^{\ord}\Xi$ such that $\pa^{\ord}\Xi^- \times \llbracket 0,\delta)^k = \pa^{\ord} X^-\cap U$.
\end{defn}
Then, we have that:
\begin{prop}\label{prop_cornervan}
	Given an a-corner $\Xi$, the operator $\cT$ is injective near $\Xi$ if the family of normal operators $N_\Xi(\cT_\Sigma)(\xi)$ is positive.
\end{prop}
\begin{proof}
	Consider a small open neighbourhood $U\cong \Xi \times \llbracket0,\delta)^k$. As before, it suffices to prove the estimate for an operator\[
	\cT = \brh^{\be_{\Xi}}(\cV_t)_\Xi-\cL_\Xi,
	\]
	where $\cL_\Xi$ and $(\cV_t)_\Xi$ are obtained by freezing coefficients at $\Xi$. Assuming that the normal operators are positive, we can write for $\hu$ being the Fourier transform in translation-invariant directions,\[
	\|\hu|_{\pa^{\ord}\Xi^-}\|^2_{L^2(\Xi)}+\Re(\la N_\Xi(\cT)(\xi)\hu,\hu \ra) \geq C\|\hu\|_{L^2(\Xi)}^2
	\]
	which by Plancherel's theorem gives us either in the case $\Xi$ is non-vanishing,\[
	\|u|_{\pa^{\circ}X^-}\|^2_{L^2}+\la \cT u,u\ra\geq C\|u\|^2_{L^2},
	\]
	which combined with H\"{o}lder's estimate proves injectivity near $\Xi$. 
	
	On the other hand, if $\Xi$ is vanishing, the normal operator only sees the relatively elliptic part and we get
	\begin{equation}\label{eq_vanpos}
	\|u|_{\pa^{\circ}X^-}\|^2_{L^2}-\la\cL_\Xi u,u\ra\geq C\|u\|^2_{L^2}.
	\end{equation}
	Now, we can write the operator $\cT=\pa_v-\cL_\Xi$, and we have in this case that\begin{equation}\label{eq_ibp}
	|\la \pa_vu,u\ra| \leq  \epsilon \cdot\|u\|^2_{L^2},
	\end{equation}
	with there being no boundary terms, and the constant $\epsilon\downarrow 0$ by taking a sufficiently small neighbourhood of $\Xi$ as the divergence of $\pa_v$ decays at the vanishing a-corner. Therefore, we have\[
	\|u|_{\pa^{\circ}X^-}\|^2_{L^2}+\la \cT u,u\ra = \|u|_{\pa^{\circ}X^-}\|^2_{L^2}-\la \cL_\Xi u,u\ra+\la \pa_vu,u\ra \geq C'\|u\|^2_{L^2},
	\]
	for $u$ supported on a sufficiently small neighbourhood of $\Xi$. Combined with H\"{o}lder's inequality, the injectivity follows.
\end{proof}
However, this is a rather strong condition unsuited for most applications. Instead, we will obtain sharp conditions for injectivity in the following sections.

\section{Analysis of model problems}\label{sec_models}
\subsection{Analysis near codimension 2 a-corners}
In this subsection, we will work with $\Xi:= \Sigma_{j}\cap\Sigma_{k}$ being the intersection of two a-boundaries, which is a compact manifold. Let us write the operator in a small neighbourhood $U$ as\[
\cT = \rho_j^{\eta_j}\rho_k^{\eta_k}F\cdot(\rho_j\pa_{\rho_j}+\rho_k\pa_{\rho_k})-a^{ij}\pa^2_{ij}-b^i\pa_i-c,
\]
with coefficients in $C^{\delta,\delta/2}$ and $F>0$ or $F<0$, and $a^{ij}$ uniformly elliptic. Our goal will be to obtain criteria for the operator to be injective near $\Xi$ in the sense of Definition \ref{defn_injective}.

\subsubsection{Non-vanishing case}
In this case, $(\eta_j,\eta_k)=(0,0)$ and the model space is $\bM:=\Xi\times \llbracket0,\infty\rrbracket \times \llbracket0,\infty\rrbracket$ with the model problem obtained by formally freezing the coefficients of $\cT$ at $\Xi$, i.e.\
\begin{equation}\label{eq_modelnonvanishcorner}
\cT_\Xi:= F(\sigma)(\rho_j\pa_{\rho_j}+\rho_k\pa_{\rho_k})-a^{ij}(\sigma)\pa^2_{ij}-b^i(\sigma)\pa_i-c(\sigma),
\end{equation}
for $\la\pa_i\ra_{i=1}^m = \la \pa_{\sigma_1},\dots,\pa_{\sigma_{m-1}},\rho_j\pa_{\rho_j}-\rho_k\pa_{\rho_k}\ra$ being the vector fields spanning the relative tangent bundle. Thus we are reduced to proving an $L^2$-estimate on $\bM$ of the form \begin{equation}\label{eq_modelcornerl2}
\|u\|_{L^2(\bM)} \leq C\|\cT_\Xi u\|_{L^2(\bM)},
\end{equation}
which is equivalent, by Fourier transform in variables $(\log(\rho_j),\log(\rho_k))$ (or equivalently, a Mellin transform in variables $\rho_j,\rho_k$) and Plancherel's theorem, to proving an estimate\begin{equation}\label{eq_uniformindicialestimatecorner}
	\|\hat{u}(\xi_j,\xi_k)\|_{\chL^2(\Xi)}\leq C\|N_\Xi(\cT_\Xi){(\xi_j,\xi_k)} \hat{u}(\xi_j,\xi_k)\|_{\chL^2(\Xi)},
\end{equation}
for all $(\xi_j,\xi_k)\in \bR^2$ and $C$ being a constant independent of $(\xi_j,\xi_k)$. This gives a well-defined family of mappings\begin{equation}\label{eq_normaloperatorcorner}
N_\Xi(\cT_\Xi){(\xi_j,\xi_k)}: \chL^2_{2}(\Xi)\rt\chL^2(\Xi),
\end{equation}
for $\chL^2_2,\chL^2$ being Sobolev spaces on $\Xi$ which are the completion of $C^\infty_0(\Xi)$-functions under norms\[
\|u\|^2_{\chL^2(\Xi)}:=\int_\Xi|u|^2dV,
\]
\[
\|u\|^2_{\chL^2_2(\Xi)}:= \|\nabla^2u\|^2_{\chL^2}+\|\nabla u\|^2_{\chL^2}+\|u\|^2_{\chL^2}.
\]

Note that $\Xi$ is an ordinary compact manifold by Assumption \ref{item_codim}, so $\nabla$ is an ordinary connection on $\Xi$. These norms are restrictions of the norms on $\bM$, in the following sense: let $\psi: (0,\infty)^2\rt [0,1]$ be a compactly supported cut-off function. Then, for any $u\in C^\infty_0(\Xi)$, we have that\[
\|u(\sigma)\|_{\chL^2(\Xi)}\equiv \|\psi(\rho_j,\rho_k)u(\sigma)\|_{L^2(\bM)}, \|u(\sigma)\|_{\chL^2_{2,1}(\Xi)}\equiv \|\psi(\rho_j,\rho_k)u(\sigma)\|_{L^2_{\be;2,1}(\bM)}.
\]

\begin{prop}\label{prop_cornernonvan}
	Given a non-vanishing codimension $2$ a-corner $\Xi$, the operator is injective near the corner $\Xi$ if and only if the family of normal operators \eqref{eq_normaloperatorcorner} has no kernel for $(\xi_j,\xi_k)\in \bR^2$.
\end{prop}
\begin{proof}
Let us assume that the family of normal operators has no kernel. As $\Xi$ is a compact manifold without boundary, standard elliptic theory implies an estimate for each $(\xi_j,\xi_k)\in \bR^2$,\begin{equation}\label{eq_injectiveestimatexi}
\|u\|^2_{\chL^2(\Xi)}\leq C_{\xi_j,\xi_k}\|N_\Xi(\cT_\Xi){(\xi_j,\xi_k)}u\|^2_{\chL^2(\Xi)},
\end{equation}
with the constant $C_{\xi_j,\xi_k}$ being continuous in $(\xi_j,\xi_k)$. Now, we have the following a-priori estimate for the model operator on $\bM$,\begin{equation}\label{eq_nondegcornerapriori}
	\|\rho_j\pa_{\rho_j}u\|^2_{L^2(\bM)}+\|\rho_k\pa_{\rho_k}u\|^2_{L^2(\bM)}\leq C_\Xi(\|\cT_\Xi u\|^2_{L^2(\bM)}+\|u\|^2_{L^2(\bM)}),
	\end{equation}
	which can be proved in the same way as the proof of a-priori estimates on $X$, via a good cover. Therefore, we have by Fourier transform of \eqref{eq_nondegcornerapriori} in $(\log(\rho_j),\log(\rho_k))$ and Plancherel,\begin{equation}\label{eq_fourierestimatexi}
	\int(|\xi_j|^2+|\xi_k|^2)\|\hat{u}\|^2_{\chL^2}d\xi_jd\xi_k\leq C_\Xi\int(\|N_\Xi(\cT_\Xi){(\xi_j,\xi_k)}\hu\|^2_{\chL^2}+\|\hu\|^2_{\chL^2})d\xi_jd\xi_k.
	\end{equation}
	Since this holds for all $\hu$, we can localise in $(\xi_j,\xi_k)$ and obtain the estimate\[
	(|\xi_j|^2+|\xi_k|^2)\|\hat{u}\|^2_{\chL^2} \leq C_\Xi(\|N_\Xi(\cT_\Xi){(\xi_j,\xi_k)}\hu\|^2_{\chL^2}+\|\hu\|^2_{\chL^2}).
	\]
	So, for $|\xi_j|^2+|\xi_k|^2>2C_\Xi$, we can get rid of the $\|u\|_{L^2}$-term on the RHS of \eqref{eq_fourierestimatexi}, while for $|\xi_j|^2+|\xi_k|^2\leq 2C_\Xi$ we can uniformly bound $\|\hu\|_{L^2}\leq C'\|N_\Xi(\cT_\Xi){(\xi_j,\xi_k)}\hu\|_{L^2}$ by \eqref{eq_injectiveestimatexi}. As a result, we obtain for all $\xi$,\[
	\|\hu\|^2_{L^2}\leq C''\|N_\Xi(\cT_\Xi){(\xi_j,\xi_k)}\hu\|^2_{L^2},
	\]
	proving estimate \eqref{eq_uniformindicialestimatecorner}.

	Conversely, suppose $\cT_{\Xi,(\xi_j^0,\xi_k^0)}$ has a non-trivial kernel $u\in C^\infty(\Xi)$ for some $(\xi_1^0,\xi_2^0)$. Then, we have that for $\tu:=\rho_j^{i\xi_j^0}\rho_k^{i\xi_k^0}u$, $\cT_\Xi(\tu)=0$ on $\bM$. Now take a sequence of $C^\infty_0(\bM)$ cutoff functions $\varphi$ such that $\|^\m\nabla^2 \varphi\|_{L^2}+\|^\m \nabla \varphi\|_{L^2} = O(1)$, while $\|\varphi\|_{L^2}\uparrow \infty$. Then, we consider functions $\varphi(\rho_j,\rho_k)\tu$ and it is easy to see that they contradict \eqref{eq_modelcornerl2}.
\end{proof}

\subsubsection{Vanishing case}\label{subsubsec_vancorner}
In this case $(\eta_j,\eta_k)\neq(0,0)$ and the model space is $\bM:= \Xi\times\llbracket0,\infty\rrbracket \times \llbracket0,\infty\rrbracket\ni (\sigma,s,\ttil)$ obtained by a reparametrisation over the interior as described in \S \ref{subsubsec_modelproblems}. The model problem given by \[
\cT_{\Xi,\deg} := \mu(\sigma)s^\alpha\ttt-a^{ij}(\sigma)\pa^2_{ij}-b^i(\sigma)\pa_i-c(\sigma)= \mu(\sigma)s^\alpha\ttt - \cL_\Xi,
\] 
for $\la\pa_i\ra$ being the vector fields spanning the relative tangent bundle given by $\pa_{\sigma^1},\dots,\pa_{\sigma^{m-1}},s\pa_{s}$, and $\alpha =(\eta_j-\eta_k)/2$, so $\mu(\sigma)s^\alpha$ is a potential that may blow up or decay at $s\in\{0,\infty\}$. Now, we claim that it suffices to prove an $L^2$-estimate on $\bM$ for $u\in C^\infty_0(\bM)$ of the form \begin{equation}\label{eq_modifiedoperator}
\|u\|_{L^2(\bM)} \leq C\|(\mu(\sigma)s^\alpha \ttt-\cL_\Xi) u\|_{L^2(\bM)},
\end{equation}
with $L^2(\bM)$ denoting the natural $L^2$-space with an m-volume form.
\begin{lem}\label{lem_equivmodel}
	The operator is injective near the corner if and only if the estimate \eqref{eq_modifiedoperator} holds.
\end{lem}
\begin{proof}
	The difference between the estimate above from the previous case of a non-vanishing corner is that the pull-back of the volume form to $\bM$ is given by $\Psi^*(dV_X)=dV_\Xi s^{-1}ds|\log(t')|^{-1} t'^{-1}dt'$ up to a non-vanishing factor. So, it differs from the standard volume form on $\bM$ by a $\log$-factor. However, this does not affect the estimates, because of the following: we can absorb the $\log$-factor into $u$, to obtain a slightly deformed operator, 
	\begin{equation}\label{eq_deformoperator}
	\mu(\sigma)s^{\alpha}\ttt-\cL_\Xi-\mu(\sigma)s^\alpha\log(\ttil)^{-1}.
	\end{equation}
	By choosing the neighbourhood $U$ of $\Xi$ in $X$ sufficiently small, we can ensure $|\mu_0(\sigma)s^\alpha\log(\ttil)^{-1}|\ll 1$, so proving estimate \eqref{eq_modifiedoperator} implies estimates for this deformed operator in the pre-image of $U$, i.e.\ on $U$.
	
	Conversely, suppose the estimate fails for a sequence of compactly supported functions $u_i$, so that these satisfy \[
	\|u_i\|_{L^2(\bM)}=1, \; \|(\mu(\sigma)s^{\alpha}\ttt-\cL_\Xi)u_i\|_{L^2(\bM)}\downarrow 0.
	\] Then, we can translate these functions in the $\log(\ttil)$-direction to functions $u_i'(\sigma,s,\ttil):= u_i(\sigma,s,\ttil+T_i)$ for $T_i\gg1$ so that i) they lie in the pre-image of $U$ and ii) under the deformed operator, their norms are still vanishing, i.e.\[
	\|(\mu(\sigma)s^{\alpha}\ttt-\cL_\Xi-\mu(\sigma)s^\alpha\log(\ttil)^{-1})u'_i\|_{L^2}\downarrow 0,
	\]
	which is possible because $u_i$ are compactly supported and $\log(\ttil)^{-1} \downarrow 0$ as $\ttil \downarrow 0$. Then, their push-forwards will contradict injectivity near $\Xi$ on $X$, proving the claimed equivalence.
\end{proof}
The estimate \eqref{eq_modifiedoperator} is equivalent by Fourier transform in $\log(\ttil)$ and Plancherel, to proving estimates \begin{equation}\label{eq_modelvanishingcorner}
	\|(\mu(\sigma)s^\alpha i\xi-\cL_\Xi)\hu(\xi)\|_{L^2(\Xi\times\llbracket0,\infty\rrbracket)}\geq C\|\hu(\xi)\|_{L^2(\Xi\times\llbracket0,\infty\rrbracket)},
\end{equation}
for a constant $C$ independent of $\xi\in \bR$. Therefore, we define a family of mappings,\begin{equation}\label{eq_mapP}
(\mu(\sigma)s^\alpha i\xi-\cL_\Xi):\chL^2_{\xi;2,1}(\Xi\times\llbracket0,\infty\rrbracket) \rt \chL^2(\Xi\times\llbracket0,\infty\rrbracket),
\end{equation}
with the Sobolev spaces being the completions of $C^\infty_0$-functions under norms:
\[
\|u\|^2_{\chL^2(\Xi\times\llbracket0,\infty\rrbracket)}:=\int_{\Xi\times\llbracket0,\infty\rrbracket}|u|^2dV
\]
\begin{equation}\label{eq_xicylnorm}
\|u\|^2_{\chL^2_{\xi;2,1}(\Xi\times\llbracket0,\infty\rrbracket)}:=\|^\m\nabla^2u\|^2_{\chL^2}+\|^\m\nabla u\|^2_{\chL^2}+\|u\|^2_{\chL^2}+\|\xi\mu(\sigma)s^\alpha u\|^2_{\chL^2}
\end{equation}
\begin{rem}
	Even though the weight $\mu(\sigma)s^\alpha$ may blow-up, these spaces are equivalent to functions in $\chL^2_{2,\loc}(\Xi\times\llbracket0,\infty\rrbracket)$ having finite \eqref{eq_xicylnorm}-norm.
\end{rem}

 Now the above is a family of elliptic operators in $\xi$ on the cylinder $\Xi\times \llbracket0,\infty\rrbracket$, with very different behaviour depending on $\alpha=0$ (in which case it is translation-invariant in $\log(s)$) or $\alpha\neq 0$ (in which case the constant term blows up at one end and vanishes at the other end). If $\alpha=0$, then the analysis is standard. So, let's assume $\alpha \neq 0$. Then the operator $\cL_\Xi$ is translation-invariant in $\log(s)$, so under the action of $\lambda\in \bR^+$, \[
(\mu(\sigma)s^{\alpha}i\xi -\cL_\Xi)u(\sigma,s) = (\mu(\sigma)s^{\alpha}i\xi \lambda^{-\alpha}-\cL_\Xi)(u(\sigma,\lambda s))
\] 
Therefore, for $\alpha\neq 0$, proving the estimate \eqref{eq_modelvanishingcorner} for all $\xi\in \bR$ is equivalent to proving it for $\xi\in\{-1,0,+1\}$. Clearly, for the estimate to hold, we require that there are no $L^2$-solutions to the p.d.e.
\begin{equation}\label{eq_preciseoperator}
	(\mu(\sigma)s^{\alpha}i\xi-\cL_\Xi)u=0 
	\end{equation}
	for $\xi \in\{-1,0,+1\}$. We show that combined with the standard elliptic estimate at $\xi=0$, it is also sufficient:
\begin{lem}
	Given $\alpha\neq 0$, the estimate \eqref{eq_modelvanishingcorner} holds for $\xi \in \bR$ if and only if the estimate holds for $\xi=0$ and the p.d.e.\ \eqref{eq_preciseoperator} has no $L^2$-kernel for $\xi \in \{-1,+1\}$.
\end{lem}
\begin{proof}
	We can assume without loss of generality that $\alpha >0$. Suppose to the contrary, we have a sequence of functions $u_i\in C^\infty_0$ such that $\cT u_i\downarrow 0$ yet $\|u_i\|_{L^2}=1$. Now by the standard a-priori estimate on $\bM$ and localisation, we get \[
	\|u\|_{L^2}+\|(\mu(\sigma)s^\alpha i\xi -\cL_\Xi)u\|_{L^2}\geq \|\xi\mu(\sigma)s^\alpha u\|_{L^2}.
	\]
	Further, injectivity at $\xi=0$ implies we can write using a cutoff function,\[
	\|\chi_0 u\|_{L^2}+\|(\mu(\sigma)s^\alpha i\xi -\cL_\Xi)u\|_{L^2}\geq \|\xi\mu(\sigma)s^\alpha u\|_{L^2},
	\]
	for $\chi_0$ a cutoff vanishing at $s=0$. Now we note that the mapping $\chL^2_{\pm1;2,1}\rt \chL^2$ given by $u\mapsto \chi_0 u$ for $\chi_0$ being a cutoff function vanishing at $s=0$ is compact. This is because the restricted norm $\|u_i|_{[a,\infty\rrbracket}\|_{L^2} < C\cdot a^{-\alpha}$, so the norm of $u_i$ cannot `escape to infinity'. As a result, we obtain by the above estimate and Lemma \ref{lem_funcana} that the mapping \eqref{eq_mapP} is semi-Fredholm. Then, injectivity gives us the estimate \eqref{eq_modelvanishingcorner}.
\end{proof}
So we define the degenerate normal operator at vanishing corners as follows, where $\beta = (\eta_j+\eta_k)/2$ and $s=\rho_j\rho_k^{-1}$: 
\begin{defn}\label{defn_cornerdegenop}
At a vanishing a-corner $\Xi$ with $\alpha \neq 0$, the family of \textit{degenerate normal operators} is defined to be
	\begin{equation}\label{eq_indicialcorn}
		N_{\Xi,\deg}(\cT_{\Xi,\deg})(\xi):=\mu(\sigma)s^\alpha i\xi+\cT|_{\Xi}: \chL^2_{\xi;2,1}(\Xi\times \bR) \rt\chL^2(\Xi\times \bR), 
	\end{equation}
	where $\mu(\sigma,s)=-\frac{1}{\beta}\pa_v(t^{-\beta})|_\Xi = \mu(\sigma)s^\alpha$ with restriction referring to the leading order term. On the other hand if $\alpha = 0$, we instead define the \textit{degenerate normal operators} as a two-parameter family,
	\begin{equation}\label{eq_indicialcornalphaneqzero}
		N_{\Xi,\deg}(\cT_{\Xi,\deg})(\xi_1,\xi_2):= \mu(\sigma)i\xi_1+N_\Xi(\cT|_{\Xi})(\xi_2).
	\end{equation}
\end{defn}

As a result, we have:
\begin{prop}\label{prop_degencorner}
	Given a vanishing codimension $2$ a-corner $\Xi$, the operator is injective near $\Xi$ if and only if i) the family of normal operators (as in Definition \ref{defn_normaloperator} is injective over all parameters $(\xi_1,\xi_2)\in \bR^2$ ii) for $\alpha=0$, the family of degenerate normal operators \eqref{eq_indicialcornalphaneqzero} is injective for $(\xi_1,\xi_2)\in \bR^2$ 
	iii) for $\alpha \neq 0$, the family of operators \eqref{eq_indicialcorn} has no kernel for $\xi \in \{-1,+1\}$.
\end{prop}

\subsection{Analysis near a non-vanishing a-boundary}\label{subsec_nonvanabound}

In this subsection, we will work with $\Sigma$ being a non-vanishing a-boundary and model space $\bM:= \Sigma \times \llbracket0,\infty\rrbracket\ni(\sigma,\rho)$with time-fibration $t'=\rho f(\sigma)$ and model problem\[
\cT_\Sigma:= \pa_w-\cL_\Sigma,
\]
where $\pa_w = (\pa_v)_\Sigma$ is a non-vanishing $\rho$-invariant time-like vector field and $\cL_\Sigma$ is a $\rho$-invariant relatively elliptic operator along the time-fibration, both obtained by freezing coefficients at $\Sigma$. Note that $\Sigma$ may have both ordinary boundaries and a-boundaries.

	\begin{figure}[htb]
	\centering
	\def\svgwidth{0.5\columnwidth}
	%% Creator: Inkscape 1.4.3 (0d15f75042, 2025-12-25), www.inkscape.org
%% Accompanies image file 'nonvanbound.pdf' (pdf, eps, ps)
\begingroup%
  \makeatletter%
  \providecommand\color[2][]{%
    \errmessage{(Inkscape) Color is used for the text in Inkscape, but the package 'color.sty' is not loaded}%
    \renewcommand\color[2][]{}%
  }%
  \providecommand\transparent[1]{%
    \errmessage{(Inkscape) Transparency is used (non-zero) for the text in Inkscape, but the package 'transparent.sty' is not loaded}%
    \renewcommand\transparent[1]{}%
  }%
  \providecommand\rotatebox[2]{#2}%
  \newcommand*\fsize{\dimexpr\f@size pt\relax}%
  \newcommand*\lineheight[1]{\fontsize{\fsize}{#1\fsize}\selectfont}%
  \ifx\svgwidth\undefined%
    \setlength{\unitlength}{329.56815224bp}%
    \ifx\svgscale\undefined%
      \relax%
    \else%
      \setlength{\unitlength}{\unitlength * \real{\svgscale}}%
    \fi%
  \else%
    \setlength{\unitlength}{\svgwidth}%
  \fi%
  \global\let\svgwidth\undefined%
  \global\let\svgscale\undefined%
  \makeatother%
  \begin{picture}(1,0.52037746)%
    \lineheight{1}%
    \setlength\tabcolsep{0pt}%
    \put(0,0){\includegraphics[width=\unitlength,page=1]{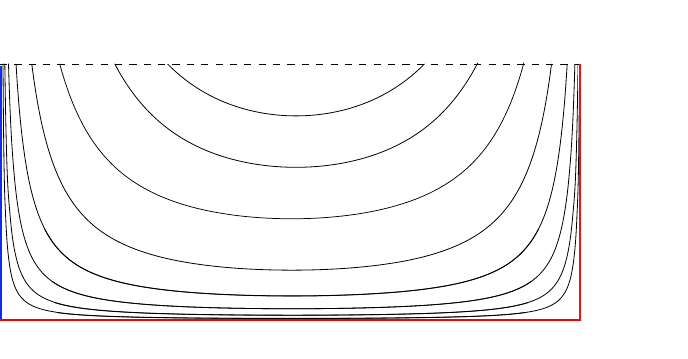}}%
    \put(0.36440189,0.00697549){\color[rgb]{0,0,0}\makebox(0,0)[lt]{\lineheight{1.25}\smash{\begin{tabular}[t]{l}$\Sigma$\end{tabular}}}}%
    \put(0,0){\includegraphics[width=\unitlength,page=2]{nonvanbound.pdf}}%
    \put(0.89978628,0.04798961){\color[rgb]{0,0,0}\makebox(0,0)[lt]{\lineheight{1.25}\smash{\begin{tabular}[t]{l}$0$\end{tabular}}}}%
    \put(0.89528749,0.41352781){\color[rgb]{0,0,0}\makebox(0,0)[lt]{\lineheight{1.25}\smash{\begin{tabular}[t]{l}$\delta$\end{tabular}}}}%
  \end{picture}%
\endgroup%

	\caption{A prototypical neighbourhood of a non-vanishing a-boundary, with the time slices. The blue boundaries are ordinary, while the red boundaries are a-boundaries.}
	\label{fig_nonvanbound}
\end{figure}
 As before, we are reduced to proving an $L^2$-estimate on $\bM$ of the form \begin{equation}\label{eq_estimateaboundary}
	\|u\|_{L^2(\bM)} \leq C(\|\cT_\Sigma u\|_{L^2(\bM)}+\|u|_{\pa^{\ord}\bM^-}\|_{L^2_1})
\end{equation}
which is equivalent by Fourier transform in $\log(\rho)$ and Plancherel's theorem, to proving an estimate\begin{equation}\label{eq_uniformindicialestimate}
	\|\hat{u}\|_{\chL^2(\Sigma)}\leq C(\|N_\Sigma(\cT_{\Sigma})(\xi) \hat{u}\|_{\chL^2(\Sigma)}+\|\hu|_{\pa^{\ord}\Sigma^-}\|_{\chL^2_1}+|\xi|\|\hu|_{\pa^\ord \Sigma^-}\|_{\chL^2}),
\end{equation}
for all $\xi \in \bR$ and $C$ being a constant independent of $\xi$. Note that we have \[
N_\Sigma(\cT_{\Sigma})(\xi) =\pa_{w,\xi}-\cL_{\Sigma,\xi},
\]
and $\cL_{\Sigma,\xi}$ is a family of uniformly $b$-elliptic operators on $\Sigma$ and $\pa_{w,\xi}$ is a family of vector fields, as can be easily checked from calculation in the local models in Definition \ref{defn_modeloperators}. We have a family of mappings, \begin{equation}\label{eq_indicialfamily}
	(N_\Sigma(\cT_{\Sigma})(\xi),|_{\pa^{\ord}\Sigma^-}): \chL^2_{2,1}(\Sigma)\rt \chL^2(\Sigma)\oplus L^2_1(\pa^{\ord}\Sigma^-),
\end{equation}
sending\[
u\mapsto (\rho^{-i\xi}(\cT_\Sigma(\rho^{i\xi}u(\sigma))|_{\Sigma}, u|_{\pa^{\ord}\Sigma^-}),
\]
 where we define the spaces $\chL^2(\Sigma),\chL^2_{2,1}(\Sigma)$ as the completion of $C^\infty_0(\Sigma)$-functions under the following norms: \[
\|u\|^2_{\chL^2(\Sigma)}:=\int_{\Sigma}|u|^2dV_{\Sigma},
\] 
\[
\|u\|^2_{\chL^2_{2,1}(\Sigma)}:=\|^b\nabla^2_{\Sigma}u\|^2_{\chL^2}+\|^b\nabla_{\Sigma}u\|^2_{\chL^2}+\|\pa_w|_\Sigma(u)\|^2_{\chL^2}+\|u\|^2_{\chL^2},
\]
where $dV_{\Sigma}$ is an m-volume form on $\Sigma$.

\begin{lem}\label{lem_localisationapriori}
	We have an a-priori estimate for $u\in \chL^2_{2,1,\loc}(\Sigma)$,
	\begin{equation}\label{eq_xiestimateapriori}
		\|u\|_{\chL^2_{2,1}(\Sigma)}\leq C_{\xi}(\|N_\Sigma(\cT_\Sigma)(\xi)u\|_{\chL^2(\Sigma)}+\|u|_{\pa^{\ord}\Sigma^-}\|_{L^2_1(\pa^\ord\Sigma^-)}+\|u\|_{\chL^2(\Sigma)}),
	\end{equation}
	for all $\xi\in \bR$, for a constant $C_\xi$ dependent on $\xi$.
\end{lem}
\begin{proof}
	We consider a function on $\bM$ of the form $\tu(\sigma,\rho):= u(\sigma)\chi(\rho)$, for $\chi$ being a smooth function with rapidly decaying derivatives.  Note that the norms $L^2_{2,1}$ on $\bM$ Fourier transform in $\log(\rho)$ give us an estimate by Plancherel's theorem,\[
	\int_{\bR}\|u(\sigma)\|^2_{\chL^2_{2,1}(\Sigma)}\hat{\chi}^2d\xi \leq \|\tu\|^2_{L^2_{2,1}(\bM)}+C\int_{\bR}(|\xi|+|\xi|^2)^2\|u(\sigma)\|^2_{\chL^2_{1,0}(\Sigma)}\hat{\chi}^2d\xi,
	\]
	for $C$ a constant independent of $\xi$. Now we take the Fourier transform in $\log(\rho)$ of the a-priori estimate \eqref{eq_apriori} on $\bM$ and invoke Plancherel's theorem and the above estimate to get
	\begin{equation}
		\int_{\bR}\|u(\sigma)\|^2_{\chL^2_{2,1}}\hat{\chi}(\xi)^2d\xi\leq C_{\xi}\left(\int_{\bR}\|N_\Sigma(\cT_\Sigma)(\xi)u(\sigma)\|^2_{\chL^2}\hat{\chi}(\xi)^2d\xi+\int_{\bR}\|u|_{\pa^{\ord}\Sigma^-}\|^2_{L^2_1}\hat{\chi}^2(\xi)d\xi+\int_{\bR}\|u(\sigma)\|^2_{\chL^2}\hat{\chi}(\xi)^2d\xi\right),
	\end{equation}
	with $C_\xi$ being continuously varying in $\xi$. Now, we can take $\hat{\chi}$ to be a bump function centred at any $\xi_0\in \bR$, which integrates to $1$. Taking the support sufficiently small, we can obtain the estimate,\[
	\|u(\sigma)\|^2_{\chL^2_{2,1}}\leq 2C_{\xi_0}\left(\|N_\Sigma(\cT_\Sigma)(\xi_0)u(\sigma)\|^2_{\chL^2}+\|u|_{\pa^{\ord}\Sigma^-}\|^2_{L^2_1}+\|u(\sigma)\|^2_{\chL^2}\right).
	\]
\end{proof} 

\begin{prop}\label{prop_semifredholmaboundary}
	The mapping $\eqref{eq_indicialfamily}$ is semi-Fredholm for all $\xi\in \bR$ if and only if the family of normal operators at the a-boundaries of $\Sigma$ have no kernel for all parameters in $(\xi_1,\xi_2)\in\bR^2$.
\end{prop}
\begin{proof}
Near the a-boundary $\pa^{\ma}\Sigma$ of $\Sigma$, $N_\Sigma(\cT_\Sigma)(\xi)$ is a $b$-elliptic operator. Therefore, it follows that for a function $u$ supported in a sufficiently small neighbourhood of $\pa^{\ma}\Sigma$, we have an estimate,
	\begin{equation}\label{eq_estimateacornernormal}
	\|u\|_{\chL^2_{2,1}}\leq C_{ \xi}(\|N_\Sigma(\cT_\Sigma)(\xi)u\|_{\chL^2}),
	\end{equation}
	if and only if the normal operator at the a-boundaries are semi-Fredholm with no kernel for parameters $(\xi,\xi_2)$ for all $\xi_2\in \bR$. Since $\pa^{\ma}\Sigma$ is a compact manifold without  boundary (as $  X$ has at most codimension $2$ corners), the normal operator, being elliptic, is automatically semi-Fredholm, and we only require the kernel to be trivial to obtain the estimate \eqref{eq_estimateacornernormal}.

	Now by Lemma \ref{lem_partition}, we can write down a covering $U\cup V$ with a partition of unity $\chi_0,\chi_1$ and decompose $u=u_0+u_1$, with $u_0$ supported away from the a-boundaries and $u_1$ supported near the a-boundaries. Then, we can write \begin{equation}\label{eq_intermediate}
	\|u\|_{\chL^2}\leq \|u_0\|_{\chL^2}+\|u_1\|_{\chL^2} \leq C_\xi(\|N_\Sigma(\cT_\Sigma)(\xi)u_1\|_{\chL^2})+\|u_0\|_{\chL^2},
	\end{equation}
	by \eqref{eq_estimateacornernormal}.
	We note that for $|^\m\nabla^2\chi|+|^\m \nabla \chi|<\delta$, we have $|N_\Sigma(\cT_\Sigma)(\xi)(\chi u)| \leq |N_\Sigma(\cT_\Sigma)(\xi)(u)|+\delta |u|_{C^{1,0}}$. Therefore, we obtain that choosing $\delta$ small, and using \eqref{eq_xiestimateapriori} and \eqref{eq_intermediate},\begin{equation}\label{eq_interxi}
	\|u\|_{\chL^2_{2,1}}\leq C'_\xi(\|N_\Sigma(\cT_\Sigma)(\xi) u\|_{\chL^2}+\|u|_{\pa^{\ord}\Sigma^-}\|_{L^2_{1}}+\|u_0\|_{\chL^2}).
	\end{equation}
	The mapping $u\mapsto u_0=\chi_0 \cdot u$ from $\chL^2_{2,1}(\Sigma)\rt \chL^2(\Sigma)$ is compact, because $\chi_0$ is supported away from the a-boundaries. Therefore by Lemma \ref{lem_funcana} and the estimate \eqref{eq_interxi}, we obtain semi-Fredholmness of \eqref{eq_indicialfamily}. 
	
	Conversely, if the family of normal operators has some kernel for some parameter $(\xi,\xi_2)\in \bR^2$, it is easy to construct a family of functions $\{u_i\}$ with $N_\Sigma(\cT_\Sigma)(\xi)(u_i)\downarrow0$ and $\|u_i\|_{\chL^2_{2,1}(\Sigma)}=1$, but $u_i$ have pairwise disjoint supports and are therefore non-subsequentially-convergent, which implies the operator is not semi-Fredholm.
\end{proof}

Now we can write the following equation,\[
\beta(\sigma)\rho\pa_\rho=\alpha(\sigma)\pa_w+\pa_\rel,
\]
	for $\beta(\sigma)$ being an a-smooth function on $\Sigma$ which is positive away from the a-boundaries, $\alpha$ being an a-smooth function on $\Sigma$ and $\pa_\rel$ being a relative a-smooth vector field (as can easily be checked from writing in the local models in Definition \ref{defn_modeloperators}, crucially using that $\pa_w$ is non-vanishing at ordinary boundaries). As a result, we can estimate on $\bM$,\begin{equation}\label{eq_etaestimate}
	\|	\beta(\sigma)\cdot\rho\pa_\rho u\|_{L^2(\bM)}\leq C(\|\pa_v u\|_{L^2(\bM)}+\|^\m\nabla_\rel u\|_{L^2(\bM)}).
	\end{equation}

Now, we will prove the following:
\begin{prop}\label{prop_aboundnonvanfredholm}
 Let us assume that the global operator $\cT$ is injective near the a-corners corresponding to the a-boundaries of $\Sigma$ (in the sense of Definition \ref{defn_injective}). Then, the operator $\cT$ on $X$ is injective near $\Sigma$ if and only if the family of normal operators \eqref{eq_indicialfamily} (which are semi-Fredholm by Proposition \ref{prop_semifredholmaboundary} and the injectivity of $\cT$ near a-corners) are injective for all $\xi \in \bR$.
\end{prop}
\begin{proof}
	Let us assume the normal operators $(N_\Sigma(\cT_\Sigma)(\xi),|_{\pa^\ord\Sigma^-})$ are injective for all $\xi \in \bR$. Firstly, we have from the a-priori estimate \eqref{eq_apriori} on $\bM$ and \eqref{eq_etaestimate}, \[
	\|\beta(\sigma)\rho\pa_\rho u\|_{L^2(\bM)}\leq C(\|\cT_\Sigma u\|_{L^2(\bM)}+\|u|_{\pa^\ord\bM^-}\|_{L^2_1}+\|u\|_{L^2(\bM)}).
	\]
	Therefore, taking the Fourier transform in $\log(\rho)$ and invoking Plancherel gives us,
	\begin{equation}\label{eq_vectorfield}
	\int_\bR\|\beta(\sigma)\xi\hu(\xi)\|^2_{\chL^2}d\xi\leq C\left(\int_\bR\|N_\Sigma(\cT_\Sigma)(\xi)\hu(\xi))\|^2_{\chL^2}d\xi++\int_\bR\|u|_{\pa^\ord\Sigma^-}\|_{\chL^2_{2,1}}^2d\xi+\int_\bR\|\hu(\xi)\|^2_{\chL^2}d\xi\right).
	\end{equation}
	
	Since the operator $\cT$ is injective at the a-corners of the model space $\bM$, $\cT$ is injective at the a-boundary of $\bM$ given by $\pa^{\ma}\Sigma \times \llbracket0,\infty\rrbracket$ due to invariance of the operator in $\rho$. Now, consider a partition of unity $\{\chi_0,\chi_1\}$ associated to the a-boundaries of $\Sigma$, as in Lemma \ref{lem_partition}. Extend this to a $\rho$-invariant partition of unity on $\bM$ associated to the a-boundary $\pa^{\ma}\Sigma\times \llbracket0,\infty\rrbracket$, by pulling back along $\bM=\Sigma\times \llbracket0,\infty\rrbracket \rt \Sigma$. Then by the injectivity of $\cT$, we can write \[
	\|u\cdot\chi_1\|_{L^2}\leq C
	(\|\cT(u\cdot\chi_1)\|_{L^2}).
	\]
	Given $|^\m\nabla^2\chi_1|+|^\m \nabla \chi_1|<\delta \ll 1$, this gives us by the same method as estimate \eqref{eq_interxi} that\[
	\|u\|_{L^2(\bM)}\leq C'(\|\cT u\|_{L^2(\bM)}+\|u|_{\pa^{\ord}\bM^-}\|_{L^2_1})+C\|u_0\|_{L^2(\bM)}.
	\]
	Now, taking the Fourier transform of this in $\log(\rho)$ and invoking Plancherel along with \eqref{eq_vectorfield} gives us,
    \begin{equation}
        \begin{split}
            &\int_\bR\|\hu(\xi)\|^2_{\chL^2}d\xi+\int_\bR\|\beta(\sigma)\xi\hu(\xi)\|^2_{\chL^2}d\xi\\
            &
            \leq C\left(\int_\bR\|N_\Sigma(\cT_\Sigma)(\xi)\hu(\xi))\|^2_{\chL^2}d\xi+\int_\bR\|u|_{\pa^\ord\Sigma^-}\|_{\chL^2_{2,1}}^2d\xi+\int_\bR\|\chi_0\cdot\hu(\xi)\|^2_{\chL^2}d\xi\right).
        \end{split}
    \end{equation}
	Now we wish to get rid of the term $\|\chi_0\cdot \hu(\xi)\|^2_{\chL^2}$ on the RHS.	Since $\beta(\sigma)>0$ away from $\pa^{\ma}\Sigma$ and $\chi_0$ vanishes in a neighbourhood of $\pa^{\ma}\Sigma$, for $|\xi|$ large the $\|\chi_0\cdot \hu(\xi)\|_{\chL^2}$-term on the RHS is dominated by $\|\beta(\sigma)\xi \hu(\xi)\|^2_{\chL^2}$ on the LHS, and since $N_\Sigma(\cT_\Sigma)(\xi)$ is injective, we can get rid of the term for $|\xi|$ small. Therefore, we get estimate \eqref{eq_uniformindicialestimate} for all $\xi\in \bR$, proving the claim.
	
	Conversely, if \eqref{eq_indicialfamily} has a non-zero element in the kernel, say $u_\xi(\sigma)$ for some $\xi \in \bR$, we can construct a counterexample similar to the proof of Proposition \ref{prop_cornernonvan}.
\end{proof}

\subsection{Analysis near a vanishing a-boundary}\label{subsec_vanabound}
In this subsection, we will work with $\Sigma$ being a vanishing a-boundary of $X$, with the local time function being given by $t'=f(\sigma)\rho^{\nu}$, for $f:\Sigma\rt \llbracket0,\infty)$ a bounded a-smooth function vanishing at the boundaries of $\Sigma$ (this can be seen from the local models in Definition \ref{defn_modeloperators}). Then, we can form the model space $\bM=\Sigma\times \llbracket0,\infty\rrbracket\ni (\sigma,\ttil)$ and write the model operator as \[
\cT_{\Sigma,\deg}:=\mu(\sigma)\ttt-\cL_{\Sigma}.
\] 
As before in section \ref{subsubsec_vancorner}, it suffices to prove an $L^2$ estimate on $\bM$ for $u\in C^\infty_0(\bM)$ of the form \[
\|u\|_{L^2(\bM)}\leq C\|(\mu(\sigma)\ttt-\cL_{\Sigma})u\|_{L^2(\bM)},
\]
by the same proof as Lemma \ref{lem_equivmodel}. Since the operator is translation-invariant in $\log(\ttil)$ the above estimate is equivalent by a Fourier transform in $\log(\ttil)$ to proving estimates \[
\|\hu(\xi)\|_{\chL^2(\Sigma)}\leq C\|(\mu(\sigma)i\xi-\cL_\Sigma)u\|_{\chL^2(\Sigma)},
\]
for a constant $C$ independent of $\xi \in \bR$. So, we define the family of \textbf{degenerate normal operators} as:
  \begin{defn}
  	Let $\Sigma$ be a vanishing a-boundary of $  X$, with rate of vanishing $\eta>0$, and such that the local time $t'$ decays at rate $\rho^\nu$ with $\nu>0$. Given the heat operator $\cT$, we define the family of \textit{degenerate normal operators} at $\Sigma$ by
  	 \begin{equation}\label{eq_indicial}
  		N_{\Sigma,\deg}(\cT_{\Sigma,\deg})(\xi):=\mu(\sigma)i\xi+\cT|_{\Sigma}: \chL^2_{\xi;2,1}(\Sigma) \rt\chL^2(\Sigma), 
  	\end{equation}
  	where $\mu(\sigma)=-\frac{1}{\eta}\pa_v(t'^{-\eta/\nu})|_\Sigma$, and the Sobolev spaces are the completions of $C^\infty_0(\Sigma)$-functions under the norms \[
  	\|u\|^2_{\chL^2(\Sigma)}:= \int_{\Sigma}|u|^2dV_{\Sigma},
  	\]
  	\[
  	\|u\|^2_{\chL^2_{\xi;2,1}(\Sigma)}:= \|^b\nabla^2u\|^2_{\chL^2(\Sigma)}+\|^b\nabla u\|^2_{\chL^2(\Sigma)}+\|\xi\mu(\sigma)u\|^2_{\chL^2(\Sigma)}+\|u\|^2_{\chL^2(\Sigma)}.
  	\]
  \end{defn}

\begin{rem}
	Note that we need the domain Sobolev space $\chL^2_{\xi;2,1}(\Sigma)$ to be dependent on $\xi$, because if $\mu(\sigma)$ blows up near a boundary of $\Sigma$, the operator has drastically different behaviour at $\xi=0$ compared to $\xi\neq 0$.
\end{rem}

We can write down an a-priori estimate by transferring the a-priori estimate from $  X$:
\begin{equation}\label{eq_aprioridegen}
	\|u\|_{L^2_{2,1}(\bM)}\leq C(\|\cT_{\Sigma,\deg} u\|_{L^2(\bM)}+\|u\|_{L^2(\bM)}),
\end{equation}
where we note that\[
\|u\|_{L^2_{2,1}(\bM)}:= \|\mu(\sigma)\ttt u\|_{L^2}+\|\nabla^2_\rel u\|_{L^2}+\|\nabla_\rel u\|_{L^2}+\|u\|_{L^2}.
\]
 Now, we prove that:
\begin{lem}\label{lem_degenboundary}
	We have an a-priori estimate,
	\begin{equation}\label{eq_apriorisigmadeg}
	\|u\|_{\chL^2_{\xi;2,1}}\leq C(\|N_{\Sigma,\deg}(\cT_{\Sigma,\deg})(\xi)u\|_{\chL^2}+\|u\|_{\chL^2}),
	\end{equation}
	for a constant $C$ independent of $\xi$. Moreover, if we assume that the operator $\cT$ is injective at the a-corners $\pa^{\ma}\Sigma  =\pa\Sigma\subset   X$ in the sense of Definition \ref{defn_injective}, then we have an estimate \begin{equation}\label{eq_apriorisigmadegpos}
	\|u\|_{\chL^2_{\xi;2,1}}\leq C'(\|N_{\Sigma,\deg}(\cT_{\Sigma,\deg})(\xi)u\|_{\chL^2}+\|\chi_0u\|_{\chL^2}),
	\end{equation}
	for $\chi_0$ being a bump function supported away from the boundaries of $\Sigma$.
\end{lem}
\begin{proof}
	Firstly at $\xi=0$, the operator is uniformly $b$-elliptic, and the a-priori estimate is a consequence of elliptic estimates. Now let us fix $\xi_0\neq 0$. We consider a function on $\bM$ of the form $\tu(\sigma,\rho):= u(\sigma)\chi(\rho)$, for $\chi$ being a smooth function with rapidly decaying derivatives. Taking the Fourier transform of the a-priori estimate \eqref{eq_aprioridegen} on $\bM$ and invoking Plancherel's theorem gives us,\begin{equation}
		\int_{\bR}\|\mu(\sigma)u(\sigma)\|^2_{\chL^2}\xi^2\hat{\chi}(\xi)^2d\xi\leq C\left(\int_{\bR}\|(i\xi\mu(\sigma)-\cL)u\|^2_{\chL^2}\hat{\chi}(\xi)^2d\xi+\int_{\bR}\|u(\sigma)\|^2_{\chL^2}\hat{\chi}(\xi)^2d\xi\right).
	\end{equation}
	Now, we can take $\hat{\chi}$ to be a bump function centred at $\xi_0$, which integrates to $1$. Taking the support sufficiently small, we can obtain an estimate,\[
	\|\xi_0\mu(\sigma)u(\sigma)\|^2_{\chL^2}\leq 2C\left(\|(i\xi_0\mu(\sigma)-\cL)u\|^2_{\chL^2}+\|u(\sigma)\|^2_{\chL^2}\right).
	\]
	The estimates for other parts of the norm follow similarly, proving the first claim \eqref{eq_apriorisigmadeg}. 
	
	Now we assume that the operator $\cT$ is injective at $\pa\Sigma=\pa^{\ma}\Sigma$. Then, we can transfer it to an estimate on $\bM$ for functions supported sufficiently close to $\pa\bM=\pa\Sigma\times \llbracket0,\infty\rrbracket$. Consider a partition of unity $\{\chi_0,\chi_1\}$ associated to the a-boundaries of $\Sigma$, as in Lemma \ref{lem_partition}. Then, we can extend this to a $t$-invariant partition of unity on $\bM$ associated to the a-boundary $\pa^{\ma}\Sigma\times \llbracket0,\infty\rrbracket$, by pulling back along $\bM=\Sigma\times \llbracket0,\infty\rrbracket \rt \Sigma$. So, by the injectivity of $\cT$, we can write \[
	\|u\cdot\chi_1\|_{L^2(\bM)}\leq C(\|\cT_{\Sigma,\deg}(u\cdot\chi_1)\|_{L^2(\bM)}).
	\]
	This gives us on taking $|^\m\nabla^2\chi|+|^\m \nabla\chi| < \delta \ll1$ that\[
	\|u\|_{L^2_{2,1}(\bM)}\leq C'\|\cT_{\Sigma,\deg} u\|_{L^2(\bM)}+C\|\chi_0 u\|_{L^2(\bM)}.
	\]
	proving estimate \eqref{eq_apriorisigmadegpos}.
\end{proof}

Now, we will prove as before that:
\begin{prop}\label{prop_fredholmvanboundary}
\begin{enumerate}
	\item The $\xi$-dependent mapping \eqref{eq_indicial} is semi-Fredholm (with respect to the $\xi$-dependent norms) for all $\xi \in \bR$ if and only if the normal operators at the a-boundaries of $\Sigma$ (which are vanishing a-corners of $  X$) are injective for $\alpha\neq 0$ and the degenerate normal operators are injective for $\alpha = 0$ (as in Definition \ref{defn_cornerdegenop}). 
	\item Let us assume that $\cT$ is an injective operator near the a-corners corresponding to the a-boundaries of $\Sigma$, in the sense of Definition \ref{defn_injective}. Then the family of degenerate normal operators \eqref{eq_indicial} (which is semi-Fredholm by the previous point) is injective for all $\xi \in \bR$ if and only if we have a uniform bound
	\begin{equation}\label{eq_discontL^2}
		\|\hu(\xi)\|_{\chL^2}\leq C'\|N_{\Sigma,\deg}(\cT_{\Sigma,\deg})\hu(\xi)\|_{\chL^2}.
	\end{equation}
\end{enumerate}
\end{prop}
\begin{proof}
	The proof proceeds similarly as in the previous subsection. We only need to be careful about the possible blow-up of $\mu(\sigma)$.
	
	\begin{enumerate}
		\item Firstly, we have the a-priori estimate \eqref{eq_apriorisigmadeg}.	The semi-Fredholmness of the operator $\mu(\sigma)i\xi-\cL_\Sigma$ depends on the behaviour of $\mu(\sigma)$ near the a-boundaries of $\Sigma$ and on $\xi$. For $\xi=0$, it is clear that we have semi-Fredholmness of $N_{\Sigma,\deg}(\cT_{\Sigma,\deg})(0)=-\cL_\Sigma$ if and only if the normal operators at the a-corners are injective. Now we assume $\xi\neq0$. Let us write $\pa^{\ma}\Sigma =\pa^{\ma \uparrow}\Sigma \sqcup \pa^{\ma \downarrow}\Sigma$ as the boundaries where $\mu(\sigma)$ blows up and where it converges respectively. Near the boundaries $\pa^{\ma\downarrow}\Sigma$, the operator is uniformly $b$-elliptic, and so injectivity of normal operators implies an estimate,\[
		\|u\|_{\chL^2_{\xi;2,1}}\leq C_\xi(\|(i\xi\mu(\sigma)-\cL_\Sigma)u\|_{\chL^2}+\|\chi_\downarrow\cdot u\|_{\chL^2}),
		\]
		for $\chi_\downarrow$ a cutoff function that vanishes in a neighbourhood of $\pa^{\ma \downarrow}\Sigma$. Now, we note that the mapping $\chL^2_{\xi;2,1}\rt \chL^2$ given by $u\mapsto \chi_\downarrow u$ is compact. This proves the claim by Lemma \ref{lem_funcana}.
		\item Due to injectivity near the a-boundaries of $\Sigma$, we obtain from Lemma \ref{lem_degenboundary} the estimate \eqref{eq_apriorisigmadegpos}
		\begin{equation}\label{eq_lemma3}
		\|\hu(\xi)\|_{\chL^2_{\xi;2,1}}\leq C(\|N_{\Sigma,\deg}(\cT_{\Sigma,\deg})(\xi)\hu(\xi)\|_{\chL^2}+\|\chi_0\cdot\hu(\xi)\|_{\chL^2}).
		\end{equation}
		Thus as before for $|\xi|\gg1$, we can estimate
		\begin{equation}\label{eq_estimatebigxi}
		\|\hu(\xi)\|_{\chL^2}\leq C'\|N_{\Sigma,\deg}(\cT_{\Sigma,\deg})(\xi)\hu(\xi)\|_{\chL^2}.
		\end{equation}
		Therefore by injectivity of the normal operators and continuity of the operators away from $\xi=0$, we obtain \[
		\|\hu(\xi)\|_{\chL^2}\leq C_\delta\|N_{\Sigma,\deg}(\cT_{\Sigma,\deg})(\xi)\hu\|_{\chL^2},		
		\]
		for a constant $C_\delta$, for all $|\xi|>\delta$. However, because of the discontinuous behaviour of the operator at $\xi=0$, it is not a-priori clear that the above can be extended to $\xi=0$ without $C_\delta$ blowing up. So, we proceed as follows: first, we have by injectivity at $\xi=0$, \[
		\|\hu(0)\|_{\chL^2} \leq C_0\|N_{\Sigma,\deg}(\cT_{\Sigma,\deg})(0)\hu\|_{\chL^2}.
		\]
		Then, we have \[
		\|\chi_0 \hu(\xi)\|_{\chL^2}\leq C_0\|N_{\Sigma,\deg}(\cT_{\Sigma,\deg})(0)(\chi_0 u)\|_{\chL^2}\leq 2C_0\|N_{\Sigma,\deg}(\cT_{\Sigma,\deg})(\xi)(\chi_0 \hu)\|_{\chL^2},
		\]
		for all $|\xi|<\delta$ for some positive $\delta$ depending on $\chi_0$, because the normal operator is $i\xi\mu(\sigma)-\cL_\Sigma$ and $\chi_0$ is compactly supported. Therefore,\[
		\|\chi_0 \hu(\xi)\|_{\chL^2}\leq 2C_0\|N_{\Sigma,\deg}(\cT_{\Sigma,\deg})(\xi)\hu(\xi)\|_{\chL^2}+2C_0\cdot \ep \|u\|_{\chL^2_1},
		\] 
		for $\ep$ depending on $\chi_0$. Plugging into \eqref{eq_lemma3}, we get \[
		\|\hu(\xi)\|_{\chL^2_{\xi;2,1}}\leq C'\|N_{\Sigma,\deg}(\cT_{\Sigma,\deg})(\xi)\hu(\xi)\|_{\chL^2}+2C'\ep\|u\|_{\chL^2_{2,1}},
		\]
		for all $|\xi|<\delta$. Now we choose $\chi_0$ such that $\ep \ll 1$, and we can cancel the last term on the right to obtain \eqref{eq_estimatebigxi} for $|\xi|<\delta$,which when combined with \eqref{eq_estimatebigxi} proves \eqref{eq_discontL^2} for all $\xi\in \bR$.
	\end{enumerate} 
\end{proof}

\section{Fredholm theorems}\label{sec_fredholm}

As a result of the previous two sections, we have the following, noting that studying the mapping $\cT$ between weighted Sobolev spaces is equivalent to studying the conjugated operator $\cT^{\bla}$:

\begin{prop}\label{thm_mainfredholm}
	 The heat operator \eqref{eq_sobolevmappingL} between the weighted Sobolev spaces at rate $\bla$ is semi-Fredholm if and only if:
	\begin{enumerate}
		\item the conjugated operator $\cT^{\bla}$ is injective near each codimension $2$ a-corner, in the sense of Definition \ref{defn_injective} (in particular, it is necessary that the family of normal operators at each a-corner is injective).
		\item the family of normal operators for $\cT^{\bla}$ at each non-vanishing a-boundary as in \eqref{eq_indicialfamily}, and the family of degenerate normal operators for $\cT^{\bla}$ at each vanishing a-boundary as in \eqref{eq_indicial}, which are Fredholm by the previous point, are injective over all parameters $\xi \in \bR$.
	\end{enumerate}
\end{prop}
\begin{proof}
	By Theorem \ref{thm_fredholmlocal}, semi-Fredholmness is equivalent to estimates of the form \eqref{eq_L^2} near each a-boundary and near each a-corner. This immediately gives us the first point. The second point follows from Propositions \ref{prop_aboundnonvanfredholm} and \ref{prop_fredholmvanboundary}.
\end{proof}

We can write it succinctly as,\begin{equation}
	\begin{split}
\textrm{semi-Fredholmness} &\iff \textrm{injectivity near a-boundaries} \\
&\iff \textrm{injectivity near codimension 2 a-corners}\\
&+\textrm{injective normal operators at non-vanishing a-boundaries}\\
&+\textrm{injective degenerate normal operators at vanishing a-boundaries}
\end{split}
\end{equation}

Now, we will study the adjoint operator to infer Fredholmness of the global operator. Let us fix an m-volume form on $X$, say for example the volume form for the m-metric on $X$.
\begin{defn}
	The \textit{adjoint} operator to $\pa_v-\cL$ is the formal $L^2$-adjoint with respect to the m-volume form on $X$, and we take a restriction at the ordinary boundaries at time $t=b$. So, we obtain a mapping,
	\begin{equation}\label{eq_mappingadjoint}
	(\cT^*,|_{t=b}): \brh^{-\bla}L^2_{\be;k+2r,r}(  X)\rt \brh^{-\bla}L^2_{\be;k+2r-2,r-1}(  X)\oplus\brh^{-\bla}L^2_{k+2r-1}(\pa^\ord  X|_{t=b})  
	\end{equation}
\end{defn}
We note that $\cT^*$ is also a parabolic operator of heat type as described in the beginning, considered over the reversed interval $[-b,-a]$. We obtain the following $L^2$ identity:
\begin{equation}\label{eq_l2ibp}
\la \cT u,v\ra_{\brh^{\bla}L^2} +\la u|_{t=a},v|_{t=a}\ra_{\brh^{\bla}L^2}=\la u^*,\cT^*v^*\ra_{\brh^{-\bla}L^2}+\la u^*|_{t=b},v^*|_{t=b}\ra_{\brh^{-\bla}L^2},   	
\end{equation}
where $u^*:=\brh^{-2\bla}u\gamma^{1/2}, v^*:= \brh^{-2\bla}v\gamma^{1/2}$.  
\begin{rem}
    There will be a uniformly positive weight function arising from the vector field $\pa_v$ associated with the measure in the inner product $\la u^*|_{t=b},v^*|_{t=b}\ra$, which we absorb into the weight function $\brh$ without loss of generality. 
\end{rem}
\begin{lem}\label{lem_parabolicreg}
	Let $u\in L^2_{\loc}(X \backslash \pa^{\ma}X)$ be a locally integrable function such that for all $v\in C^\infty_0(X\backslash t^{-1}(a))$, we have $\la (\pa_v-\cL)v,u\ra=0$. Then, we have that $u\in C^\infty_\loc(X\backslash \pa^{\ma}X)$ and all derivatives of $u$ vanish at the final ordinary boundary $t^{-1}(b)$ of $X$.  
\end{lem}
\begin{proof}
	In the interior of $X$, the claim is standard parabolic regularity. It remains to show regularity at the ordinary boundary $\pa^{\ord}X$. By localisation, it suffices to show regularity for $\pa_v-\cL$ being the heat operator on the model space $\bR^m \times [0,\ep)$. Here, note that we may extend $u$ by $0$ on $t<0$, to obtain an extension $\tu$ that is still in $L^2_{\loc}$ and satisfies the orthogonality for any $v\in C^\infty_0(\bR^m\times(-\infty,\infty))$. Therefore, we are reduced to standard parabolic regularity and $\tu\in C^\infty_{\loc}$, which implies $u$ is smooth up to the boundary and vanishing at the boundary.
\end{proof}
\begin{lem}\label{lem_l2cokernel}
	Let $\bla$ be a rate at which the mapping \eqref{eq_sobolevmappingL} for $k=0,r=1$ is semi-Fredholm. Then, the $L^2$-cokernel of its image is naturally identified with the kernel of the mapping \eqref{eq_mappingadjoint}, via the correspondence $(u,u|_{t=a})\leftrightarrow u^*$, for $u^*= \brh^{-2\bla}u$.
\end{lem}
\begin{proof}
	Suppose that $(u,g)$ lies in the $L^2$-cokernel of $(\cT,|_{t=a})$. First, by taking test functions $v\in C^\infty_0(X\backslash t^{-1}(a))$, $\la \cT v, u\ra =0$, so $u$ is locally smooth on $  X\backslash t^{-1}(a)$ and vanishing at the final ordinary boundary $\pa^{\ord}X\cap t^{-1}(b)$ by Lemma \ref{lem_parabolicreg}. Then, it follows from \eqref{eq_l2ibp} that $\cT^*(u^*)=0$ in the interior (by taking test functions supported away from $t=a$ and $t=b$), and since $u^*|_{t=b}=0$ we get $g=u|_{t=a}$. Now we have by the a-priori estimate on $\cT^*$, that \[
	\|u^*\|_{\brh^{-\bla}L^2_{\be;2,1}}\leq C_{\bla}\|u^*\|_{\brh^{-\bla}L^2},
	\]
	from which the higher regularity of $u^*$ follows. Conversely, if an element lies in the kernel of $\cT^*$, it is in the $L^2$-cokernel of the image of $\cT$ by the integration-by-parts formula.
\end{proof}

Therefore, $\cT^{\bla}$ being semi-Fredholm implies that the adjoint controls the Fredholmness of $\cT^{\bla}$, so we can state the main Theorem, which follows from Lemma \ref{lem_l2cokernel} and Proposition \ref{thm_mainfredholm} applied to the conjugate operator:
\begin{theorem}\label{thm_mainthmfred}
	Let $X$ be a manifold with corners and a-corners, $t: X\rt I=[a,b]$  be a proper m-fibration, and $\cT=\pa_v-\cL$ a heat-type operator on $X$, satisfying Assumption \ref{ass_cond}. The mapping \eqref{eq_sobolevmappingL} between the weighted Sobolev spaces,\[
    (\cT, |_{ \pa^{\ord}  X^-}) : \brh^{\bla}L^2_{\be;k+2r,r}(  X) \rt \brh^{\bla}L^2_{\be;k+2r-2,r-1}(  X) \oplus \brh^{\bla}L^2_{k+2r-1}(\pa^\ord  X^-),
    \] is Fredholm if and only if:
	\begin{enumerate}
		\item the conjugate operator $\cT^{\bla}$ and its adjoint $\cT^{*-\bla}$ are injective near each codimension $2$ a-corner, in the sense of Definition \ref{defn_injective}. 
		\item the family of normal operators for $\cT^{\bla}$ at each non-vanishing a-boundary as in \eqref{eq_indicialfamily}, and the family of degenerate normal operators for $\cT^{\bla}$ at each vanishing a-boundary as in \eqref{eq_indicial} (which are Fredholm by the previous point and Propositions \ref{prop_semifredholmaboundary} and \ref{prop_fredholmvanboundary}) are isomorphisms over all parameters $\xi \in \bR$.
	\end{enumerate}
		By Propositions \ref{prop_cornernonvan} and \ref{prop_cornervan}, injectivity near a-corners is implied if the families of normal operators at each non-vanishing a-corner are isomorphisms and at each vanishing a-corner are positive. 
\end{theorem}
We can write it succinctly as,
	\begin{alignat*}{2}
		\textrm{Fredholmness} &\iff &&\textrm{ injective and surjective near a-boundaries} \\
		&\iff &&\textrm{ injective and surjective near codimension 2 a-corners,}\\
		&  &&\textrm{+ normal operators at non-vanishing a-boundaries are isomorphisms,}\\
		& &&\textrm{+ degenerate normal operators at vanishing a-boundaries are isomorphisms.}
	\end{alignat*}

Let us also note the following useful lemma:
\begin{lem}\label{lem_zeroindexfamily}
	Suppose that the conjugate operator $\cT^{\bla}$ and its adjoint $\cT^{*-\bla}$are injective near each codimension 2 a-corner. Then,
	\begin{enumerate}
		\item At the non-vanishing a-boundaries,the family of normal operators is Fredholm and zero-index for all $\xi\in \bR$. 
		\item  At the vanishing a-boundaries, the family of degenerate normal operators is Fredholm and zero index for $\xi\in \bR\backslash \{0\}$.
	\end{enumerate} 
\end{lem}
\begin{proof}
	Injectivity near the corners implies that the family of normal operators at the non-vanishing a-boundaries and the family of degenerate normal operators at the vanishing a-boundaries are Fredholm. Now we split into the two cases of vanishing and non-vanishing a-boundaries. At the non-vanishing a-boundaries, it is clear from the proof of Proposition \ref{prop_aboundnonvanfredholm} that the normal operator is injective for $|\xi|\gg1$ (since the term $\|\beta(\sigma)\xi\hu(\xi)\|^2_{\chL^2}$ dominates $\|\chi_0\cdot\hu(\xi)\|^2_{\chL^2}$). Similarly, the adjoint normal operator is also injective for $|\xi|\gg1$. Therefore the normal operator is an isomorphism for $|\xi|\gg1$, so it is index zero. Since the family of normal operators are all Fredholm by Propositions \ref{prop_semifredholmaboundary} and \ref{prop_fredholmvanboundary} and continuously varying in $\xi$, this proves the claim.
	
	At the vanishing a-boundaries, a similar argument as above gives us that the family of normal operators is zero index for all $\xi\neq0$ (which does not extend over because of the discontinuity of the degenerate normal operator at $\xi=0$).
\end{proof}

\subsection{Optimal regularity and H\"{o}lder estimates}\label{sec_optimalsobolev}
The $\be$-Sobolev and $\be$-H\"{o}lder spaces have rather weak regularity. Therefore, we define new Sobolev spaces that are in a sense intermediary between the ordinary Sobolev spaces and $\be$-Sobolev spaces, following ideas of Mazzeo \cite{Maz91a,Maz91b}. We define a modified Sobolev space, as the completion of $C^\infty_0$-functions under the norm: \begin{equation}\label{eq_mixednorm}
\|u\|^2_{L^2_{*;k+2r,r}}:= \|u\|^2_{L^2_{k+2r-2,r-1}}+\|\pa_vu\|^2_{L^2_{k+2r-2,r-1}}+\|^\m\nabla_\rel^2u\|^2_{L^2_{k+2r-2,r-1}},
\end{equation}
with the $L^2_{k+2r,r}$-norms as defined in \S \ref{subsec_ordinarysobolev}. Therefore, the main difference between the norms \eqref{eq_sobolevnormeta} and \eqref{eq_mixednorm} is that in the former, we weight the mixed derivatives with $\brh^{\be}$, while in the latter, we only weight the terms with the highest-order time derivative. In other words, we can write,
\begin{equation}\label{eq_starspace}
\|u\|^2_{L^2_{{*;k+2r,r}}} \equiv \sum_{\substack{\beta\leq r-1 \\ \alpha\leq k+2r-2\beta}} \|{^\m\nabla^\beta}({^\m\nabla}_\rel^\alpha u)\|^2_{L^2} + \sum_{\alpha\leq k} \|\brh^{\be}{^\m\nabla^r}(^\m \nabla_\rel ^\alpha u)\|^2_{L^2}.
\end{equation}
This norm makes $L^2_{*;k+2r,r}$ the largest Banach space for which the heat operator $\cT = \pa_v-\cL$ defines a mapping into $L^2_{k+2r-2,r-1}\oplus L^2_{k+2r-1}$. The linearised operator \eqref{eq_sobolevmappingL} gives a well-defined bounded map between weighted spaces,
\begin{equation}\label{eq_sobolevmappingL2}
	(\cT, |_{\pa^\ord X^-}) : \brh^{\bla}L^2_{*;k+2r,r}(X) \rt \brh^{\bla}L^2_{k+2r-2,r-1}(X) \oplus \brh^{\bla}L^2_{k+2r-1}(\pa^\ord  X^-). 
\end{equation}
Note that $L^2_{\be;k+2r,r}\equiv L^2_{*;k+2r,r}$ for $r\leq1$ by the representation \eqref{eq_starspace}, so we obtain an a-priori estimate from \eqref{eq_apriori},
\[
\|u\|_{\brh^{\bla}L^2_{*;2,1}}\leq C(\|\cT u\|_{\brh^{\bla}L^2}+\|u|_{\pa^\ord X^-}\|_{\brh^{\bla}L^2_1}+\|u\|_{\brh^{\bla}L^2}).
\]
And we can extend it to a higher order estimate by taking m-derivatives,\[
\|u\|_{\brh^{\bla}L^2_{*;k+2r,r}}\leq C(\|\cT u\|_{\brh^{\bla}L^2_{k+2r-2,r-1}}+\|u|_{\pa^\ord X^-}\|_{\brh^{\bla}L^2_{k+2r-1}}+\|u\|_{\brh^{\bla}L^2_{k+2r-2,r-1}}).
\]
However, unlike \eqref{eq_apriorihigher}, we do \textbf{not} have an a-priori estimate of the form \[
\|u\|_{\brh^{\bla}L^2_{*;k+2r,r}}\leq C(\|\cT u\|_{\brh^{\bla}L^2_{k+2r-2,r-1}}+\|u|_{\pa^\ord X^-}\|_{\brh^{\bla}L^2_{k+2r-1}}+\|u\|_{\brh^{\bla}L^2}),
\]
for a general $r>1$ and weight $\bla$. This is because we do not have higher order interior estimates near the vanishing a-boundaries for $r>1$, as in the following example:
\begin{ex}
	Consider an operator $T=t^2\pa_t+i$ on the interval $t\in \llbracket0,1]$, so that $T$ vanishes to order $1$ at $\{0\}$ and the weight $\be=\{1\}$. Then, considering $u= e^{i t^{-1}}$, we have that $Tu\equiv 0$, and over an interval $(a,b)$, $\|u|_{(a/2,4a)}\|_{L^2}=|\log(8)|^{1/2}$. However, for $r>1$, the norm $\|u|_{(a,2a)}\|_{L^2_{*;r}}\uparrow \infty$ as $a\downarrow 0$, while $\|u|_{(a,2a)}\|_{L^2_{\be;r}}=O(1)$ as $a\downarrow 0$.
\end{ex}
Nevertheless, we can show that:
\begin{theorem}\label{prop_optimall2Fredholm}
	If we are at a rate $\bla$ such that the mapping between $\be$-Sobolev spaces \eqref{eq_sobolevmappingL} is semi-Fredholm, we have an a-priori estimate for $u\in \brh^{\bla}L^2_{*;k+2r,r}(X)$,
	\begin{equation}\label{eq_star}
	\|u\|_{\brh^{\bla}L^2_{*;k+2r,r}}\leq C(\|\cT u\|_{\brh^{\bla}L^2_{k+2r-2,r-1}}+\|u|_{\pa^\ord X^-}\|_{\brh^{\bla}L^2_{k+2r-1}}+\|u\|_{\brh^{\bla}L^2}).
	\end{equation}
\end{theorem} 
\begin{proof}
	As mentioned before, the spaces $L^2_*$ are equivalent to the spaces $L^2_{\be}$ away from the vanishing a-boundaries. Therefore, by the general methods in the previous sections, it suffices to prove the estimate \eqref{eq_star} for $u$ being compactly supported in a small neighbourhood $U\cong \Sigma \times \llbracket0,\delta)$ of a vanishing a-boundary. Thus from the local models in Definition \ref{defn_modeloperators} we can assume without loss of generality that $\cT=\rho^{\eta}f(\sigma)(\rho\pa_\rho)-\cL_\Sigma$, with the coefficients frozen at $\Sigma$. Now we already know the estimate \eqref{eq_star} for $r=1$, since then the $L^2_*$-spaces and the $L^2_{\be}$-spaces coincide. Then semi-Fredholmness at rate $\bla$ implies that\begin{equation}\label{eq_optimalest}
		\|u\|_{\brh^{\bla}L^2}\leq C\| (\rho^{\eta}f(\sigma)(\rho\pa_\rho)-\cL_\Sigma)u\|_{\brh^{\bla}L^2}
	\end{equation}
	for $u$ supported $U$ with $\delta$ sufficiently small.We take derivatives in $\rho\pa_{\rho}$ to get\[
	\rho\pa_{\rho}\cT u = (\cT +\eta\rho^\eta f(\sigma))(\rho\pa_{\rho}u).
	\]
	This is because of the commutations $[\rho\pa_\rho,\rho^\eta\rho\pa_\rho f(\sigma)] = \eta \rho^{\eta}f(\sigma)\rho\pa_\rho$, $[\rho\pa_\rho,\cL_\Sigma]=0$. Clearly the mapping properties of $(\cT +O(\rho^\eta))$ are the same as that of $\cT$ for $\eta>0$. Therefore, we can obtain the same estimate \eqref{eq_optimalest} for one higher derivative, and continue inductively to obtain \eqref{eq_star} in a small neighbourhood of the a-boundary, then use partitions of unity as before to obtain a global estimate.
\end{proof}

Moreover, we can upgrade the regularity of solutions:
\begin{theorem}\label{prop_optimalreg}
	Suppose we are working at a rate $\bla$ such that the mapping \eqref{eq_sobolevmappingL} is semi-Fredholm. Given an equation\[
	\cT u=f,\;\; u|_{ \pa^\ord X^-}=g,
	\]
	for $f\in \we L^2_{k+2r-2,r-1}, g\in \brh^{\bla}L^2_{k+2r-1}( \pa^\ord X^-)$ and $u\in \brh^{\bla}L^2_{\be;2,1}(X)$, it follows that $u\in \brh^{\bla}L^2_{*;k+2r,r}$.
\end{theorem}
\begin{proof}
As before, the regularity is immediate away from the vanishing a-boundaries from standard a-priori estimates. Then, we work near a vanishing a-boundary as before, and by Fredholmness at $\bla$, we obtain the same $L^2$-estimate \eqref{eq_optimalest} near $\Sigma$, and we would like to claim that this implies a bound on $\rho\pa_\rho u$. However, this is not a valid proof, since we do not know a-priori that $\rho\pa_\rho u$ is integrable, while the estimate \eqref{eq_optimalest} assumes that $u\in L^2$.
	
	Therefore, we use the following trick: we approximate $u$ by $u_i:=u\cdot \chi_i$ for suitable cutoff functions, so that we know a-priori integrability of $u_i$, then show that the derivatives $\rho\pa_\rho u_i$ are bounded in $L^2$ and converge pointwise to $\rho\pa_\rho u$.

	Let us assume $u$ is supported compactly near a vanishing a-boundary which corresponds to a local time variable $t'=0$. Define cutoff functions $\chi_\ep:=\varphi(\ep^{-1}t)$, for $\varphi:[0,\infty)\rt [0,1]$ being a smooth cutoff function with $\varphi(x)\equiv 0$ for $x\in (0,1)$, $\varphi(x)\equiv 1$ for $x\in (2,\infty)$, and $\varphi$ is smooth, non-decreasing. Then, it is easy to see that all m-derivatives of $\chi_{\ep}$ are uniformly bounded in $\ep$ away from the ordinary boundaries. Now, we consider $u_\ep:=\chi_{\ep} \cdot u$, which have compact support away from the a-boundaries and by interior regularity, $\rho\pa_{\rho}u_\ep\in L^2_{\be;2,1}$. Therefore from \eqref{eq_optimalest} we obtain \[
	\|\de u_\ep\|_{L^2}\leq C\|\cT(\de u_\ep)\|_{L^2}\leq C(\|\de \cT(u_\ep)\|_{L^2}+\|[\de,\cT](u_\ep)\|_{L^2}).
	\]
	Then we can compute,\[
	\cT(u_\ep)=\chi_{\ep}\cdot f+\cT^0(\chi_{\ep})u+\nabla_\rel \chi_\ep \cdot \nabla_\rel u.
	\]
	Clearly, $\cT(\chi_\ep)=O(\brh^{\be})$, since $\chi_\ep$ is constant on time-slices. Then, we have\[
	\de \cT(u_\ep)=\de(\chi_\ep\cdot f)+O(\pa_vu)+O(u),
	\]
	and the commutator as before,\[
	[\rho\pa_\rho,\cT](u_\ep) = \eta\rho^\eta f(\sigma)\rho\pa_\rho(u_\ep) = O(\pa_v u)+O(u).
	\]
	As a result, we have a uniform bound in $\ep$,\[
	\|\de u_\ep\|_{L^2}\leq C(\|f\|_{L^2_{2,1}}+\|u\|_{L^2_{\be;2,1}}).
	\]
	Moreover, as the support of $\de \chi_\ep$ shrinks as $\ep\downarrow0$, the derivatives $\de u_i$ converge over an exhaustion of the neighbourhood of the a-boundary to $\de u$, and as the $L^2$ norms are bounded, we have $\de u_i \rt \de u$ in $L^2$, which provides an extra derivative of $u$. Continuing similarly, we obtain derivatives of all orders of $u$.
\end{proof}

\begin{coro}
	Suppose we are at a rate $\bla$ where the mapping \eqref{eq_sobolevmappingL} is Fredholm. Then, the mapping \eqref{eq_sobolevmappingL2} is Fredholm as well. Moreover, the kernel and cokernel of the mapping are independent of $k,r$ for all $k,r>1$, and are the same for $\be$-Sobolev spaces or $*$-Sobolev spaces.
\end{coro}
\begin{proof}
	It is clear that if an element lies in the kernel of the lowest order mapping (i.e.\ $k=0,r=1$) then it lies in all higher order mappings (i.e.\ $k,r\geq1$) by the above claim, while the converse inclusion is clear, so we have an equality of the kernel over $k,r\geq1$. The claim for the cokernel follows from Lemma \ref{lem_l2cokernel}.
\end{proof}
\subsubsection{H\"{o}lder estimates}
In order to improve regularity of solutions to nonlinear problems via bootstrapping, it will be useful to have optimal H\"{o}lder estimates at Fredholm rates, similar to Theorem \ref{prop_optimall2Fredholm}. First, we define the modified H\"{o}lder spaces, as the completion of $C^\infty_0$-functions under the norm: \begin{equation}\label{eq_holdernormoptimal}
	\|u\|_{C^{*;k+2r+\delta,r+\delta/2}}:= \|u\|_{C^{k+2r-2+\delta,r-1+\delta/2}}+\|\pa_vu\|_{C^{k+2r-2+\delta,r-1+\delta/2}}+\|^\m\nabla_\rel^2u\|_{C^{k+2r-2+\delta,r-1+\delta/2}}.
\end{equation}
Note that the $C^{k+2r+\delta,r+\delta/2}$-norms are as defined in \S \ref{subsec_ordinarysobolev}. Therefore, the main difference between the $C^{\be}$-norms \eqref{eq_holdernormeta} and $C^*$-norms \eqref{eq_holdernormoptimal} is that in the former, we weight each mixed derivative with $\brh^{\be}$, while in the latter, we weight only the terms with the highest order time derivatives. 
Then, it is easy to check that we have a well-defined bounded linear mapping,
\begin{equation}\label{eq_optimalholdermapping}
	(\cT, |_{\pa^\ord X^-}) : \brh^{\bla}C^{*;k+2r+\delta,r+\delta/2}(X) \rt \brh^{\bla}C^{k+2r-2+\delta,r-1+\delta/2} \oplus \brh^{\bla}C^{k+2r+\delta}(\pa^\ord  X^-). 
\end{equation}
Now, we will prove that at a Fredholm rate $\bla$, we have the following estimate on a sufficiently small open neighbourhood of each a-boundary, analogous to estimate \eqref{eq_L^2}:
\begin{equation}\label{eq_holderl2}
	\|u\|_{{*;2+\delta,1+ \delta/2}}\leq C(\|\cT^{\bla}u\|_{\delta,\delta/2}+\|u|_{\pa^\ord X^-}\|_{C^{2+\delta}}).
\end{equation}
We will prove this by contradiction, following a well-known blow-up argument, see for example Brendle--Kapouleas \cite[\S 6]{BK16} and Su--Tsai--Wood \cite[\S 6]{STW24}. First we have a Liouville theorem ruling out solutions to the model problems as defined in \S \ref{subsubsec_modelproblems}.
\begin{lem}\label{prop_C0kernel}
	Given a semi-Fredholm rate $\bla$ for the mapping \eqref{eq_sobolevmappingL}, the model operators on model spaces as defined in \S \ref{subsubsec_modelproblems} have trivial $C^{\be;2,1}$-kernel.
\end{lem} 
\begin{proof}
	Let us first consider model problems at the non-degenerate a-corners. Suppose to the contrary, that there exists a $C^{\be;2,1}$-solution $u$ to the model problem on the model space $\bM$ with $\|u\|_{\be;C^{2,1}}=1$, so that $\cT u=0$ and $u|_{\pa^{\ord}\bM^-}\equiv 0$. Then we can consider the mapping problem at a slightly weaker weight $\bla-\boldsymbol{\ep}$ (by abuse of notation, this refers to weighting at both ends $\rho\rt0$ and $\rho\rt \infty$ for each non-compact direction $\rho$), at which the mapping will still be Fredholm. From the $C^0$-estimate on $u$, it is clear that $\|u\|_{\brh^{\bla-\boldsymbol{\ep}}L^2(\bM)}<\infty$, a contradiction to the $L^2$-injectivity of the model problems at a semi-Fredholm rate. 
	
	However, this weighting does not immediately work for the model problems associated to a degenerate a-boundary or a-corner, because of the possible blow-up of coefficients of the operator in the time-direction. In particular, consider the model problem on $\bM=\Sigma \times \llbracket0,\infty\rrbracket$ for $\Sigma$ being an a-boundary or $\Sigma = \Xi\times \llbracket0,\infty\rrbracket \ni (\sigma,s)$, with a $C^{\be;2,1}$-solution \[
    \cT_{\Sigma,\deg}u=
	(\mu(\sigma)s^\alpha\ttt-\cL_\Sigma)u=0,
	\]
	(where the potential $\mu(\sigma)s^\alpha$ may blow-up at the boundaries of $\Sigma$). We can weight along the $s$-direction (again, by abuse of notation, we are referring to weighting in both directions $s\rt 0$ and $s\rt \infty$) to obtain integrability in the space-directions, however weighting in the time-direction introduces a possibly unbounded extra term given by $\epsilon \mu(\sigma)s^\alpha u$
	so the operator upon weighting in $\ttil$ may {not} be a small perturbation of the model problem. Instead, we define for a fixed $h>0$, \[
	u_h:=u(x,th)-u(x,t).
	\]
	Then, it is clear that $u_h$ satisfies the model problem $\cT_{\Sigma,\deg}u_h = (F(\sigma)\ttt-\cL_\Sigma)u_h \equiv 0$ and lies in $C^{\be;2,1}$, and moreover we can estimate
	\begin{equation}\label{eq_diffest}
	|F(\sigma)u_h|\leq \left|\int_t^{th}|F(\sigma)\ttt u| \ttil^{-1}d\ttil\right|\leq |\log(h)|,
	\end{equation}
	since $\|u\|_{C^{\be;2,1}}=1$. Now we can weight in the space direction to obtain integrability of $u_h$ in the space direction, and then weight in the time direction by a cutoff function $\chi(\ttil)$ to get \[
	(F(\sigma)\ttt-\cL_\Sigma)(\chi(\ttil)u_h)=\ttt\chi(\ttil)\cdot F(\sigma)u_h,
	\] 
	and thus by \eqref{eq_diffest} and the integrability of $u_h$, the $L^2$-norm of the above is uniformly bounded. Therefore, we obtain a uniform estimate on $\|\chi(\ttil)u_h\|_{L^2}$ by the $L^2$-injectivity of the model problem, and consequently we get $L^2$-convergence of $\chi(\ttil)u_h \rt u_h$ by taking cutoff functions $\chi$ whose supports exhaust the time domain. Thus $\|u_h\|_{L^2}$ is finite, and by $L^2$-injectivity of the model problem, it follows that $u_h\equiv 0$. Since $h$ is arbitrary, it follows that $u(\sigma,\ttil)\equiv u(\sigma)$, so the resulting model problem becomes\[
	\cT_{\Sigma,\deg}u = \cL_\Sigma u=0.
	\]
	Thus we are reduced to the p.d.e.\ $\cL_\Sigma u =0$ on $\Sigma$. Because we know $L^2$-injectivity of the normal operator at $\xi=0$, $\cL_\Sigma$ is $L^2$-injective. Weighting by a small amount in the space direction implies that $u\equiv 0$, proving the claim. 

\end{proof}

Now, suppose to the contrary that estimate \eqref{eq_holderl2} fails. Let $\{u_i\}_{i=1}^\infty$ be a sequence of functions supported in an open neighbourhood of the a-boundaries such that $\|u_i\|_{0,0}=1$ and $\|\cT^{\bla}u_i\|_{\delta,\delta/2}\downarrow 0$. By the Schauder estimate, $\|u_i\|_{C^{\be;2+\delta,1+\delta/2}}<C$. Now, we will attempt to find a subsequence of $u_i$ that converges to a non-zero function $u \in C^{\be;\delta,\delta/2}$ solving the model problem. Let us choose points $x_i\in U$ such that $|u_i(x_i)|>1/2$. Then, we have three possibilities:
\begin{enumerate}
	\item $x_i\rt x\in U\backslash\pa^{\ma}U$: In this case, $u_i\rt u \in C^{2+\delta,1+\delta/2}$ pointwise to a non-vanishing function. Therefore, $\cT^{\bla} u_i \rt {\cT}^{\bla}u $ in $C^{\delta,\delta/2}$, so $\cT^{\bla}u\equiv 0,u|_{\pa^\ord X^-}\equiv 0$. However, working at a slightly lower Fredholm rate, $u$ is integrable, so this contradicts the injectivity estimate \eqref{eq_L^2}. 
	\item $x_i\rt x$ for a point in the interior of an a-boundary. In this case, we reduce the problem to one of the local models in \S\ref{subsec_nonvanabound}, \S\ref{subsec_vanabound}, so that the operator is translation-invariant in $\rho$ or $\ttil$. Let us write the points in coordinates $x_i=(\sigma_i,\rho_i)$. Then, we can translate $u_i$ to $\tu_i(\sigma,\rho):=u_i(\sigma, \rho_i\cdot\rho)$ so that $\tu_i(\sigma_i,1)\rt u(\sigma,1)$ for a fixed point $\sigma\in \Sigma$. As a result, we can obtain a subsequential convergence of $\tu_i\rt u\neq 0$ in $C^{\be;\delta,\delta/2}$, and by the Schauder interior estimates \eqref{eq_schauder}, $\tu_i\rt u \in C^{2+\delta,1+\delta/2}$ pointwise to a non-vanishing function on the model space $\bM$ with $\cT_\Sigma^{\bla}u\equiv 0,u|_{\pa^\ord \bM^-}\equiv 0$, but this is ruled out by Lemma \ref{prop_C0kernel}.
	\item $x_i\rt x$ for a point at a codimension $2$ a-corner: Here, if the corner is non-vanishing, the problem is translation-invariant, and we can do the same trick as before by translating in two variables to obtain a contradiction. If the corner is vanishing, however, the problem is no longer translation-invariant. Instead, we firstly change coordinate so that we are in the model problem described in \S \eqref{subsubsec_vancorner}, on model space $\bM=\Xi\times \llbracket0,\infty\rrbracket \times \llbracket 0,\infty\rrbracket\ni (\sigma,\rho,t)$ and model operator $\mu_0(\sigma)\rho^\alpha t\pa_{t}-\cL'$. Then, suppose the points are given by $x_i=(\sigma_i,\rho_i,t_i)$ with $\rho_i,t_i\downarrow 0$. For a fixed $1>\ep>0$, consider balls centered at $(\sigma,\rho_i,t_i)$ given by $\Xi\times(\ep\rho_i,\ep^{-1}\rho_i
	)\times(\ep\rho_i^{\alpha} t_i,\ep^{-1}\rho_i^{\alpha}t_i)$. Then, we can define diffeomorphisms $\Psi_i$ from $(\sigma,x,s)\in\Xi \times (-|\log(\ep)|,|\log(\ep))\times(-|\log(\ep)|,|\log(\ep)|)$ to these balls given by $\Psi_i(\sigma,x,s):=(\sigma, \rho_ie^x,\rho_i^{\alpha}t_i e^{ s})$. Considering the pull-back functions $\Psi_i^{-1}(u_i)$, they will subsequentially converge on $\Xi \times (-|\log(\ep)|,|\log(\ep))\times(-|\log(\ep)|,|\log(\ep)|)$ to a solution of the p.d.e.\ $\mu_0(\sigma)e^{\alpha x}\pa_s-\cL$ in $C^{2+\delta,1+\delta/2}$ by Schauder interior estimates. By taking $\ep\downarrow 0$, we obtain convergence on larger balls and thus a $C^{2+\delta,1+\delta/2}$-solution to $\mu_0(\sigma)e^{\alpha x}\pa_s-\cL$ on the whole of $\Xi\times \bR\times \bR$, but this is a solution to the model problem via an exponential change of variable, which is ruled out by Lemma \ref{prop_C0kernel}. 
\end{enumerate}
Thus, we have that:

\begin{theorem}
	At a Fredholm rate $\bla$ for the $L^2$-problem, the H\"{o}lder mapping \eqref{eq_optimalholdermapping} is Fredholm as well. Moreover, the kernel and cokernel of the mapping are the same as that for the $L^2$-Sobolev mapping.
\end{theorem}
\begin{proof}
	Combined with the a-priori estimate \eqref{eq_schauder}, the above discussion proves that the mapping is semi-Fredholm for $k=0,r=1$. To extend it to higher orders, we follow exactly the same methods as in \S \ref{sec_optimalsobolev}, replacing $L^2$-norms with $C^0$-norms.
	
	Once we know the image of the mapping is closed, to prove Fredholmness it suffices to show by Theorem \ref{prop_optimalreg} that any equation of the form \[
	\cT^{\bla}u=f, u|_{\pa^\ord X^-} = g
	\]
	has a solution, for $f,g$ being compactly supported smooth functions orthogonal to the cokernel. But this can be solved in the $L^2$-spaces for $k,r \gg1$, and clearly by the standard Sobolev embedding \eqref{eq_parabolicholderembedding} this implies that the solutions are in the H\"{o}lder spaces as well. This proves the claim. 
\end{proof}
\begin{theorem}[A-priori estimate]
	Given an equation of the form\[
	\cT u=f,\; u|_{ \pa^{\ord}X^-}=g,
	\]
	with the coefficients of $\cT$ with regularity $C^{*;k+2r+\delta,r+\delta/2}$, the terms $f\in \brh^{\bla}C^{k+2r+\delta,r+\delta/2}$ and $g\in \brh^{\bla}C^{k+2r+2+\delta}(\pa^\ord X^-)$, and $\bla$ is a Fredholm rate, then given $u\in \brh^{\bla}C^{\be;2,1}$, we obtain $u\in \brh^{\bla}C^{*;k+2r+2+\delta,r+1+\delta/2}$ and we have an estimate\[
	\|u\|_{\brh^{\bla}C^{*;k+2r+2+\delta,r+1+\delta/2} } \leq C(\|f\|_{\brh^{\bla}C^{k+2r+\delta,r+\delta/2}}+\|g\|_{\brh^{\bla}C^{k+2r+2+\delta}} + \|u\|_{\brh^{\bla}C^{0}}).
	\]
\end{theorem}

\subsection{Change-of-index formula}\label{subsec_coi}
We will write down a change-of-index formula for the operator as the rate $\bla$ crosses some critical rate, following a similar method to Melrose--Mendoza \cite[\S 6]{MM83}. Let us assume that we are at a rate $\bla=(\lambda_1,\dots, \lambda_N)$ such that the conjugated operator $\cT^{\bla}$ and its adjoint $\cT^{*-\bla}$ are both injective near each codimension $2$ a-corner. Let $\Sigma_j$ be a non-vanishing a-boundary, and suppose that for a positive $0<\delta\ll 1$ we have Fredholmness of the operator at weights $\bla_-:=(\lambda_1,\dots, \lambda_j-\epsilon,\dots,\lambda_N)$ and $\bla_+:=(\lambda_1,\dots,\lambda_j+\epsilon,\dots, \lambda_N)$ for all $0<\epsilon<\delta$. Recall from Lemma \ref{lem_l2cokernel} that an element $(u,g)$ in the cokernel of the mapping at rate $\bla$ corresponds to an element $u^*= \brh^{-2\bla}u$ in the kernel of the adjoint mapping $(\cT^*,|_{t=b})$ at the rate $-\bla$. Therefore, as the rate $\bla$ increases from one Fredholm rate to another, the kernel of the operator is non-increasing while the cokernel is non-decreasing. First, we note the following Lemma about elements of the cokernel:
\begin{lem}\label{lem_compactsupport}
	We can choose finite-dimensional vector spaces of functions $V_{\bla^-}\subset V_{\bla^+}$ which are supported away from the a-boundaries and $V_{\bla^\pm}$ is complementary to the image of the mapping \eqref{eq_sobolevmappingL} for $k=0,r=1$ in the range $\brh^{\bla^\pm}L^2\oplus \brh^{\bla^\pm}L^2_1$ at rates $\bla^\pm$.
\end{lem}
\begin{proof}
	Let us write the kernel of the adjoint operator at rate $-\bla^-$ as spanned by $\la u_1^*,\dots, u_L^*\ra$ and at rate $-\bla^+$ spanned by $\la u_1^*,\dots, u_{L+M}^*\ra$. Consider a cutoff function $\chi$ that is compactly supported away from the a-boundaries such that the elements $\la u_1^*,\dots,u_{L+M}^*\ra$ are linearly independent on $\textrm{supp}(\chi)$. Define $V_{\bla^-}:=\la \chi \cdot \brh^{2\bla}u_1^*,\dots, \chi\cdot \brh^{2\bla}u_L^*\ra$, and $V_{\bla^+}:=\la \chi \cdot \brh^{2\bla}u_1^*,\dots, \chi\cdot \brh^{2\bla}u_{L+M}^*\ra$. These are integrable at all rates due to the vanishing of $\chi$. Now, it is clear that these are linearly independent, and they are not $L^2$-orthogonal to the cokernel of the mapping at rates $\bla^-$ and $\bla^+$, so they cannot lie in the image. 
\end{proof}
 Now, we note the following result analogous to Lemma 6.14 of \cite{MM83}:
\begin{lem}\label{lem_auxiliaryvect}
	For $\epsilon$ sufficiently small, there exists a finite-dimensional family of functions $V\subset \brh^{\bla^-}L^2_{2,1}$ such that $V\cap \brh^{\bla^+}L^2_{2,1}=\{0\}$ and $u\in \brh^{\bla^-}L^2_{2,1}$ satisfies $(\cT u , u|_{t=a})\in \brh^{\bla^+}L^2\oplus \brh^{\bla^+}L^2_{1}$ if and only if $u\in \brh^{\bla^+}L^2_{2,1}\oplus V$. Further, we have $\dim(V)=\ind(\bla^-)-\ind(\bla^+)$.
\end{lem}
\begin{proof}
	Let $W = \Ima_{\bla^-}\cap V_{\bla^+}\subset \brh^{\bla^-}L^2\oplus \brh^{\bla^-}L^2_1$ where $\Ima_{\bla^-}$ is the image of the mapping at rate $\bla^-$. It is easy to see that this is a complementary subspace to $V_{\bla^-}\subset V_{\bla^+}$ (as defined in Lemma \ref{lem_compactsupport}). Then its inverse image modulo the kernel is well-defined, $W^{-1}:=(\cT,|_{t=a})_{\bla^-}^{-1}W\cap \ker_{\bla^-}^\perp$. It is clear that $W^{-1}\cap \brh^{\bla^+}L^2_{2,1}=\{0\}$, since otherwise, the image of the mapping at rate $\bla^+$ would intersect $V_{\bla^+}$ (which is complementary to its image by Lemma \ref{lem_compactsupport}). We then define the vector space\[
	V:= K \oplus W^{-1},
	\]
	where $K$ is a complement of $\Ker_{\bla^+}$ in $\Ker_{\bla^-}$. Now, it is clear that $V\cap \brh^{\bla^+}L^2_{2,1}=\{0\}$. Given $u$ satisfying the equation, we can uniquely invert $u$ modulo a kernel at rate $\bla^-$, and then the claimed decomposition follows. Finally, the dimension of $V$ is immediate because\begin{equation}\begin{split}
			\ind(\bla^-)-\ind(\bla^+)&=\dim(\Ker_{\bla^-}/\Ker_{\bla^+})+\dim(\Coker_{\bla^+}/\Coker_{\bla^-})\\
			&=\dim(K)+\dim(W^{-1}).
		\end{split}
	\end{equation}
\end{proof}

\begin{defn}\label{def_pol}
	Consider the translation-invariant problem $(\cT_\Sigma u,|_{\pa^{\ord}\bM^-} u)=0$ on the model space $\bM=\Sigma_j\times \llbracket0,\infty\rrbracket$. Define the vector space of \textit{polyhomogeneous solutions} at $\Sigma_j$ at rate $\bla$, denoted $W_{\bla}$, to be the vector space of solutions to the problem spanned by finite linear combinations of functions of the form $u=\rho_j^{\lambda_j+i\xi}\log(\rj)^na(\sigma)$, for $\xi\in \bR$, $n\in \bZ_{\geq 0}$ and $a(\sigma)$ being an a-smooth function (with weight) on $\Sigma_j$.
\end{defn}

\begin{theorem}\label{prop_polyhomogeneous}
	 Let $\cT$ be a heat type operator with a-smooth coefficients. The elements that appear/disappear in the cokernel/kernel on going from $\bla^-$ to $\bla^+$ (i.e.\ $\Coker^+/\Coker^-$ and $\Ker^-/\Ker^+$) are leading-order polyhomogeneous functions of $\rho_j$ near the a-boundary $\Sigma_j$ with coefficients in $C^{\infty}(\Sigma_j)$. Further, the index of the operator changes according to the count of the total kernel of polyhomogeneous solutions,\[
		\ind(\bla^+)-\ind(\bla^-)=-\dim(W_{\bla}),
		\]
		as defined in Definition \ref{def_pol}.
		In particular, the space of polyhomogeneous solutions is finite.
\end{theorem}
\begin{proof}We already know that the elements in the kernel and cokernel of the operator at a Fredholm rate are a-smooth at a slightly lower rate, by Theorem \ref{prop_optimalreg}. Now, we consider germs of $C^\infty(X)$ functions at the a-boundary $\Sigma_j$, i.e.\ we consider pairs \[\{(U,f):U \textrm{ is an open neighbourhood of $\Sigma_j$ and $f$ is a-smooth on $U$}\},\]
	under the equivalence relation $(U_1,f_1)\sim(U_2,f_2)$ if there exists an open neighbourhood $U_{12}\subset U_1\cap U_2$ with $f|_{U_{12}}=g|_{U_{12}}$. We denote this vector space by $\mathscr{C}(\Sigma_j)$. Then, define weighted spaces as before, $\brh^{\bla^+}\mathscr{C}(\Sigma_j), \brh^{\bla}\mathscr{C}(\Sigma_j)$ and $\brh^{\bla^-}\mathscr{C}(\Sigma_j)$.

	Now, we define a vector space,\[
	\cB_{\bla}(\Sigma_j):= \{u\in \brh^{\bla^-}\sC/\brh^{\bla^+}\sC : (\cT u,u|_{\pa^\ord X^-}) \in \brh^{\bla^+}\sC\oplus \brh^{\bla^+}\sC(\pa^\ord)\}.
	\]
	
	\begin{lem}
		There exists a natural isomorphism between $B_{\bla}$ and $V$ as defined in Lemma \ref{lem_auxiliaryvect}.
	\end{lem}
	\begin{proof}
		Given an element $a\in V$, it is clear that we can restrict it to a small neighbourhood $U$ to obtain a nonzero element in $\cB_{\bla}$. Thus we have a natural injection $V\hookrightarrow \cB_{\bla}$. Now, given linearly independent elements $\{u_1,\dots,u_k\}$ of $\cB_{\bla}$, we can transfer it to $  X$ such that $(\cT u_i,u_i|_{\pa^\ord X^-}) \in \brh^{\bla^+}L^2\oplus \brh^{\bla^+}L^2_1$, so we can either i) solve away the error at a slightly higher Fredholm rate by adding lower-order terms, or ii) be unable to solve it, and hence lose an element in the cokernel. Since we are adding lower order terms, the leading order contribution comes from $u_i$ and hence it is a non-trivial element. As a result, we get $\dim(\cB_{\bla})\leq \ind(\bla^-)-\ind(\bla^+)$. Therefore, it follows that the injection is an isomorphism and $\dim(\cB_{\bla})= \dim(V) = \ind(\bla^-)-\ind(\bla^+)$.
	\end{proof}
	
	\begin{lem}
		The elements of $\cB_{\bla}$ are pairs $(U,f)$ where to leading order $f$ is a logarithmic polynomial, i.e.\ $f=\rj^{\lambda_j+i\xi}\sum_{k=0}^n\log(\rj)^ka_k(\sigma)+O(\rho_j^{\lambda_j+\delta})$, for $\xi\in \bR$ and $a_i(\sigma)$ being a-smooth functions on $\Sigma$.
	\end{lem}
	\begin{proof}
		If an element $u\in \brh^{\bla^-}\mathscr{C}$ on an open neighbourhood $U$ satisfies $\cT u \in \brh^{\bla^+}\mathscr{C}$, then we can take its derivative in $\rho\pa_\rho$, to get\[
		\cT(\rho\pa_\rho u) = \rho\pa_\rho(\cT u) + [\cT,\rho\pa_\rho]u =  \rho^{\delta_1} \rho\pa_\rho(u),
		\]
        modulo $\brh^{\bla^+}\mathscr{C}$,
		where $\delta_1> 0 $ depends on the rate of decay of the coefficients of $\cT$. Similarly, we can iterate derivatives and obtain \[
		\cT((\rho\pa_\rho)^ku)=\sum_{j=1}^k \rho^{\delta_j}(\rho\pa_\rho)^ju,
		\]
        modulo $\brh^{\bla^+}\mathscr{C}$.
		Now, by taking $\bla^-$ very close to $\bla$, we can ensure that the RHS above lies in $\brh^{\bla^+}\mathscr{C}$ for any given finite $k$. Therefore, taking $k= \ind(\bla^-)-\ind(\bla^+)$, we have that all the elements $u,\rho\pa_{\rho}u,\dots,(\rho\pa_\rho)^k u \in \cB_{\bla}$. 
		 Moreover, as $k=\dim(\cB_{\bla})$, a linear combination of these elements lies in $\brh^{\bla^+}\mathscr{C}$. Therefore, on a smaller neighbourhood, $u(\sigma,\rj)$ satisfies a pointwise o.d.e.\ in $\rj$ on an interval $\llbracket0,\kappa)$, given by \[
		 \sum_{i=0}^ka_i(\rho\pa_\rho)^iu(\sigma,\cdot)=O(\rho_j^{\lambda_j+\delta}).
		 \] Since the initial conditions at $\rj=\kappa$ vary smoothly due to interior regularity of the kernel, we immediately obtain the claimed decomposition.
	\end{proof}
	This completes the proof of Theorem \ref{prop_polyhomogeneous}.
\end{proof}

\begin{rem}
	If $\Sigma_j$ is a \textit{vanishing} a-boundary, then there are no critical rates $\bla$ satisfying the condition of the theorem. The reason is that writing out the weights in the local model gives us $\log(t')^\lambda s^{-\lambda}$ near the corners up to bounded factors. The $\log$-weight does not change the operator, and weighting by $s$ implies that an element that appears in the kernel at rate $\lambda$ will also appear in the kernel at rate $\lambda^-$, and conversely, anything in the cokernel will appear at rate $\lambda^+$. In particular, the set of non-critical rates at $\Sigma_j$ for a fixed set of rates at all other a-boundaries forms a (possibly empty) open interval.
\end{rem}

\subsection{Extension to vector bundles}\label{subsubsec_vecval}
	Now we formulate our p.d.e.\ for a vector bundle $E$ over $  X$, considering an operator of the form $\pa_v-\cL$ where $\cL:\Gamma^\infty(E)\rt \Gamma^\infty(E)$ is a certain relative \textit{strongly elliptic} operator, and our results will go through. We begin with a an m-connection $\nabla^E$ on $E$, and $\nabla^E_\rel$ the induced relative connection. We define a heat operator on $E$ as follows:\[
	\cT = \pa_v-\cL: \Gamma^\infty(E)\rt \Gamma^\infty(E),\;\; \cT(u) :=  \nabla^E_{\pa_v}u-\cL(u),
	\] 
	for $\cL$ being a relative strongly elliptic operator on $E$ and $u$ being a section of $E$. Now, as before we can note down local models for this problem, which allow us to prove interior estimates. Then, we can define normal operators as before, and obtain Fredholm theorems.
		
	\subsubsection{Sobolev spaces}
	To define Sobolev spaces, we begin with a smooth metric on the fibres of $E$. Then, we define Sobolev spaces as usual, \[
	\|u\|^p_{L^p_{\be;k+2r,r}}:=\sum \|(\nabla^E_{\pa_v})^j(\nabla^E_\rel)^k u\|^p_{L^p}.
	\]
	Note that to differentiate sections $\nabla^E_\rel u  \in \Gamma^\infty(E\otimes {^\m}T_\rel^*  X)$, we need to transport it to sections of $E\otimes {^\m}T^*  X$, which we can by dualizing the inclusion $0\rt {^\m}T_\rel  X\rt {^\m}T  X$ with respect to the m-metric.
	This gives us a well-defined mapping, \[
	(\cT, |_{t=a}): L^2_{\be;k+2r,r}\rt L^2_{\be;k+2r-2,r-1}\oplus L^2_{k+2r-1}.
	\]
	Due to $\cL$ being strongly elliptic, we have the vector-bundle interior estimates of \S \ref{subsec_vectorvalued}, and as a result, the proofs of interior estimates in the scalar case immediately carries over to give us a-priori estimates. 
	\subsubsection{Normal operators}
	We define normal operators as before, with the only difference being that they will be defined on sections of $E|_\Xi$ for an a-corner $\Xi$. We also define the degenerate normal operators at a vanishing a-boundary as before, defined on sections of $E|_\Sigma$. Finally, we define the notion of positivity via the $L^2$-inner product arising from a choice of metric on the restricted vector bundles.

    \begin{defn}\label{defn_normaloperatorvec}
	Let $M$ be a codimension $k\geq 1$ a-corner component of $  X$. Let the boundary-defining functions corresponding to $M$ be $\rho_{i_1},\dots, \rho_{i_k}$, so that $M\subseteq \Sigma_{i_1}\cap\cdots\cap \Sigma_{i_k}$. Define the \textit{family of normal operators} associated to $\cT$ as follows: for a parameter $\xi \in \bR^k$, given a section $u\in C^{\infty}(E|_M)$, define
	\begin{equation}\label{eq_normaloperatordefnvec}
	N_M(\cT)(\xi)(u):= \brh^{-i\xi}\cT(\brh^{i\xi}\tilde{u})|_{M},
	\end{equation}
	for $\tilde{u}$ being a smooth extension of $u$ on a small neighbourhood $U$ of $M$. As a result, we obtain a family of operators dependent on the parameter $\xi\in \bR^k$.
\end{defn}
Similarly we define degenerate normal operators:
 \begin{defn}
  	Let $\Sigma$ be a vanishing a-boundary of $  X$, with rate of vanishing $\eta>0$, and such that the local time $t'$ decays at rate $\rho^\nu$ with $\nu>0$. Given the heat operator $\cT$, we define the family of \textit{degenerate normal operators} at $\Sigma$ by
  	 \begin{equation}\label{eq_indicialvec}
  		N_{\Sigma,\deg}(\cT_{\Sigma,\deg})(\xi):=\mu(\sigma)i\xi+\cT|_{\Sigma}: \chL^2_{\xi;2,1}(E|_\Sigma) \rt\chL^2(E|_\Sigma), 
  	\end{equation}
  	where $\mu(\sigma)=-\frac{1}{\eta}\pa_v(t'^{-\eta/\nu})|_\Sigma$, and the Sobolev spaces are the completions of $C^\infty_0(\Sigma)$-functions under the norms \[
  	\|u\|^2_{\chL^2(E|_\Sigma)}:= \int_{\Sigma}|u|^2dV_{\Sigma},
  	\]
  	\[
  	\|u\|^2_{\chL^2_{\xi;2,1}(E|_\Sigma)}:= \|^b\nabla^2u\|^2_{\chL^2(E|_\Sigma)}+\|^b\nabla u\|^2_{\chL^2(E|_\Sigma)}+\|\xi\mu(\sigma)u\|^2_{\chL^2(E|_\Sigma)}+\|u\|^2_{\chL^2(E|_\Sigma)}.
  	\]
  \end{defn}

	\subsubsection{Fredholm theorems}
	Given these definitions, it is straightforward to prove by the exact same methods that the same Fredholm theorems hold for the operator $\cT$. 
    \begin{theorem}\label{thm_mainthmfredvec}
	Let $X$ be a manifold with corners and a-corners, $t: X\rt I=[a,b]$  be a proper m-fibration, and $E\downarrow X$ be a vector bundle. Let $\cT=\nabla^E_{\pa_v}-\cL$ a heat-type operator on sections of $E$, satisfying Assumptions \ref{ass_cond}. The mapping between the weighted Sobolev spaces,\[
    (\cT, |_{ {\pa^{\ord}  X^-}}) : \brh^{\bla}L^2_{\be;k+2r,r}(E) \rt \brh^{\bla}L^2_{\be;k+2r-2,r-1}(E) \oplus \brh^{\bla}L^2_{k+2r-1}(E|_{\pa^\ord  X^-}),
    \] is Fredholm if and only if:
	\begin{enumerate}
		\item the conjugate operator $\cT^{\bla}$ and its adjoint $\cT^{*-\bla}$ are injective near each codimension $2$ a-corner, in the sense of Definition \ref{defn_injective}. 
		\item the family of normal operators for $\cT^{\bla}$ at each non-vanishing a-boundary as in \eqref{eq_indicialfamily}, and the family of degenerate normal operators for $\cT^{\bla}$ at each vanishing a-boundary as in \eqref{eq_indicial} (which are Fredholm by the previous point and Propositions \ref{prop_semifredholmaboundary} and \ref{prop_fredholmvanboundary}) are isomorphisms over all parameters $\xi \in \bR$.
	\end{enumerate}
		By Propositions \ref{prop_cornernonvan} and \ref{prop_cornervan}, injectivity near a-corners is implied if the families of normal operators at each non-vanishing a-corner are isomorphisms and at each vanishing a-corner are positive. 
\end{theorem}

\bibliographystyle{alphaurl}
\bibliography{bibliography}

@article{AN63,
  author    = {S. Agmon and L. Nirenberg},
  title     = {Properties of solutions of ordinary differential equations in {Banach} space},
  journal   = {Communications on Pure and Applied Mathematics},
  volume    = {16},
  pages     = {121--239},
  year      = {1963}
}

@article{BK16,
author = {S. Brendle and N. Kapouleas},
title = {Gluing {Eguchi-Hanson} Metrics and a Question of {Page}},
journal = {Communications on Pure and Applied Mathematics},
volume = {70},
number = {7},
pages = {1366-1401},
year = {2017}
}

@book {BH18,
    AUTHOR = {B\"uhler, T. and Salamon, D.},
     TITLE = {Functional analysis},
    SERIES = {Graduate Studies in Mathematics},
    VOLUME = {191},
 PUBLISHER = {American Mathematical Society, Providence, RI},
      YEAR = {2018},
     PAGES = {xiv+466},
      ISBN = {978-1-4704-4190-6},
}

@book{Eva10,
  author    = {L. C. Evans},
  title     = {Partial Differential Equations},
  edition   = {2nd},
  publisher = {American Mathematical Society},
  year      = {2010}
}

@article{Gho26b,
  author    = {S. Ghosh},
  title     = {{Neck pinching in Lagrangian mean curvature flow: the non-generic case}},
  journal      = {in preparation}
}

@article{GS18,
author = {P. Gianniotis and F. Schulze},
title = {{Ricci} flow from spaces with isolated conical singularities},
journal = {Geometry and Topology},
volume = {22},
number = {7},
pages = {3925-3977},
year = {2018}
}

@article{Hui90,
author = {G. Huisken},
title = {Asymptotic behaviour for singularities of the mean curvature flow},
journal = {Journal of Differential Geometry},
volume = {31},
number = {1},
pages = {285--299},
year = {1990}
}

@article{Joy15,
  author    = {D. Joyce},
  title     = {Conjectures on {Bridgeland stability for Fukaya categories of Calabi--Yau manifolds, special Lagrangians, and Lagrangian} mean curvature flow},
  journal   = {EMS Surveys in Mathematical Sciences},
  volume    = {2},
  number    = {1},
  pages     = {1--62},
  year      = {2015}
}

@misc{Joy16,
  author    = {D. Joyce},
  title     = {Manifolds with analytic corners},
  year      = {2016},
  eprint    = {arXiv:1605.05913},
  archivePrefix = {arXiv},
  primaryClass = {math.DG},
}

@book{Kry96,
  author    = {N. V. Krylov},
  title     = {Lectures on Elliptic and Parabolic Equations in H\"older Spaces},
  series    = {Graduate Studies in Mathematics},
  volume    = {12},
  publisher = {American Mathematical Society},
  year      = {1996}
}

@book{Kry08,
  author    = {N. V. Krylov},
  title     = {Lectures on Elliptic and Parabolic Equations in Sobolev Spaces},
  series    = {Graduate Studies in Mathematics},
  volume    = {96},
  publisher = {American Mathematical Society},
  year      = {2008}
}

@book{LSU68,
  author    = {O. A. Ladyzhenskaya and V. A. Solonnikov and N. N. Ural'tseva},
  title     = {Linear and Quasilinear Equations of Parabolic Type},
  series    = {Translations of Mathematical Monographs},
  volume    = {23},
  publisher = {American Mathematical Society},
  year      = {1968}
}

@article{LMPS23,
  author    = {J. Lira and R. Mazzeo and A. Pluda and M. S\'aez},
  title     = {Short-time existence for the network flow},
  journal   = {Communications on Pure and Applied Mathematics},
  volume    = {76},
  pages     = {3968--4021},
  year      = {2023}
}

@article{Loy04, 
author= {P. Loya}, 
 title= {Index theory of Dirac operators on manifolds with corners up to codimension two}, 
 journal={Aspects of boundary problems in analysis and geometry},
  pages ={ 131--166}, 
   editor    = {J. Gil et al.},
  publisher = {Birkh\"auser},
  address   = {Basel}, 
  year={2004}
  }

@article{Loy05, 
author= {P. Loya}, 
 title= {The index of b-pseudodifferential operators on manifolds with corners}, 
 journal={Ann. Global Anal. Geom.},
  pages ={101--133}, 
  volume={27},
  year={2005}
  }

@article{LM03, 
author= {P. Loya and R.B. Melrose}, 
 title= {Fredholm perturbations of elliptic operators on manifolds with corners}, 
 journal={preprint},
  year={2003}
  }

@article{Maz91a,
author = {R. Mazzeo},
title = {{Elliptic theory of differential edge operators I}},
journal = {Communications in Partial Differential Equations},
volume = {16},
number = {10},
pages = {1615--1664},
year = {1991},
publisher = {Taylor \& Francis},
}

@article{Maz91b,
 author = {R. Mazzeo},
 journal = {Indiana University Mathematics Journal},
 number = {4},
 pages = {1277--1299},
 publisher = {Indiana University Mathematics Department},
 title = {{Regularity for the singular Yamabe problem}},
 urldate = {2026-07-15},
 volume = {40},
 year = {1991}
}

@inproceedings{Mel90,
  author    = {R. Melrose},
  title     = {Pseudodifferential operators, corners and singular limits},
  booktitle = {Proceedings of the International Congress of Mathematicians, Kyoto},
  pages     = {217--234},
  year      = {1990},
  publisher = {Springer-Verlag}
}

@book{Mel93,
  author    = {R. Melrose},
  title     = {The Atiyah–Patodi–Singer Index Theorem},
  publisher = {A. K. Peters},
  address   = {Wellesley, MA},
  year      = {1993}
}

@misc{Mel96,
  author    = {R. Melrose},
  title     = {Differential Analysis on Manifolds with Corners},
  year      = {1996},
  note      = {unfinished book available at \url{https://math.mit.edu/~rbm/book.html}}
}

@inproceedings{MM83,
  title={Elliptic Operators of Totally Characteristic Type},
  author={R. Melrose and G. Mendoza},
  year={1983},
  url={https://api.semanticscholar.org/CorpusID:46475850}
}

@article{Per02,
  author    = {G. Perelman},
  title     = {{The entropy formula for the Ricci flow and its geometric applications}},
  year      = {2002},
  eprint    = {/math/0211159v1},
  primaryClass = {math.DG}
}

@article{Per03a,
  author    = {G. Perelman},
  title     = {{Ricci flow with surgery on three-manifolds}},
  year      = {2003},
   eprint    = {/math/0303109},
  primaryClass = {math.DG}
}

@article{Per03b,
  author    = {G. Perelman},
  title     = {{ Finite extinction time for the solutions to the Ricci flow on certain three-manifolds}},
  year      = {2003},
   eprint    = {/math/0307245},
  primaryClass = {math.DG}
}

@article{STW24,
  author    = {W.-B. Su and C.-J. Tsai and A. Wood},
  title     = {Infinite Time Singularities of the {Lagrangian} Mean Curvature Flow},
  year      = {2024},
  eprint    = {arXiv:2401.02228},
  archivePrefix = {arXiv},
  primaryClass = {math.DG}
}

@article{Sto26,
  author    = {M. Stolarski},
  title     = {Closed {Ricci} flows with singularities modeled on asymptotically conical shrinkers},
  journal   = {Geometry and Topology},
  volume    = {30},
  pages     = {1277-1370},
  year      = {2026}
}

@article{LZ24,
  author    = {T.-K. Lee and X. Zhao},
  title     = {Closed mean curvature flows with asymptotically conical singularities},
  year      = {2024},
   eprint   ={arXiv:2405.15577},
  primaryClass = {math.DG}
}

@book{Top06,
  title={{Lectures on the Ricci Flow}},
  author={Topping, P.},
  series={London Mathematical Society Lecture Note Series},
  volume={325},
  year={2006},
  publisher={Cambridge University Press}
}

@article{Gho25,
  author        = {S. Ghosh},
  title         = {Short-time existence of {Lagrangian} mean curvature flow},
  journal   ={The Journal of Geometric Analysis},
  volume  ={292},
    year  = {2026}
}

\end{document}